\documentclass{article}

\usepackage{arxiv}

\usepackage[utf8]{inputenc} 
\usepackage[T1]{fontenc}    
\usepackage{booktabs}       
\usepackage{amsfonts}       
\usepackage{nicefrac}       
\usepackage{microtype}      
\usepackage{lipsum}
\usepackage{graphicx}

\usepackage{subcaption}
\usepackage{float}

\usepackage{amsfonts,amsthm}
\usepackage{amssymb,amsmath,latexsym,tensor}
\usepackage{mathrsfs}
\usepackage{enumerate}
\usepackage{savesym}
\usepackage{graphicx}
\usepackage{adjustbox}
\usepackage[pdftex]{color}
\usepackage[font=footnotesize,labelfont=bf]{caption}
\savesymbol{AND}
\usepackage{epstopdf}
\usepackage{algorithm}
\usepackage{algorithmic}

\newcommand{\R}{\mathbb{R}}

\DeclareMathOperator{\E}{\mathbb{E}}
\DeclareMathOperator{\Tr}{Tr}
\DeclareMathOperator*{\argmin}{arg\,min}
\DeclareMathOperator*{\argmax}{arg\,max}

\usepackage{mathtools}
\DeclarePairedDelimiterX{\norm}[1]{\lVert}{\rVert}{#1}
\DeclarePairedDelimiterX{\inp}[2]{\langle}{\rangle}{#1, #2}

\usepackage{cleveref}

\newtheorem{theorem}{Theorem}[section]
\newtheorem{corollary}{Corollary}[theorem]
\newtheorem{lemma}[theorem]{Lemma}
\newtheorem{remark}[theorem]{Remark}

\newtheorem{definition}[theorem]{Definition}

\usepackage{xcolor}
\usepackage{mdframed}

\usepackage{array,multirow}

\title{Sketch \& Project Methods for Linear Feasibility Problems:  Greedy Sampling \& Momentum}

\author{
 Md Sarowar Morshed \\
Department of Mechanical $\&$ Industrial Engineering\\
Northeastern University\\
Boston, MA \\
  \texttt{morshed.m@northeastern.edu} \\
  \And
 Md. Noor-E-Alam \\
Department of Mechanical $\&$ Industrial Engineering\\
Northeastern University\\
Boston, MA \\
Corresponding author: \texttt{mnalam@neu.edu} \\
}

\begin{document}
\maketitle
\begin{abstract}
We develop two greedy sampling rules for the \textit{Sketch \& Project} method for solving linear feasibility problems. The proposed greedy sampling rules generalize the existing max-distance sampling rule and uniform sampling rule and generate faster variants of \textit{Sketch \& Project} methods. We also introduce greedy capped sampling rules that improve the existing capped sampling rules. Moreover, we incorporate the so-called heavy ball momentum technique to the proposed greedy \textit{Sketch \& Project} method. By varying the parameters such as sampling rules, sketching vectors \textemdash we recover several well-known algorithms as special cases, including \textit{Randomized Kaczmarz} (RK), \textit{Motzkin Relaxation} (MR), \textit{Sampling Kaczmarz Motzkin} (SKM). We also obtain several new methods such as \textit{Randomized Coordinate Descent}, \textit{Sampling Coordinate Descent}, \textit{Capped Coordinate Descent} etc. for solving linear feasibility problems. We provide global linear convergence results for both the basic greedy method and the greedy method with momentum. Under weaker conditions, we prove $\mathcal{O}(\frac{1}{k})$ convergence rate for the Cesaro average of sequences generated by both methods. We extend the so-called certificate of feasibility result for the proposed momentum method that generalizes several existing results. To back up the proposed theoretical results, we carry out comprehensive numerical experiments on randomly generated test instances as well as sparse real-world test instances. The proposed greedy sampling methods significantly outperform the existing sampling methods. And finally, the momentum variants designed in this work extend the computational performance of the \textit{Sketch \& Project} methods for all of the sampling rules.

\end{abstract}



\section{Introduction}
\label{sec:intro}
In this work, we consider the problem of solving the following linear feasibility problem:
\begin{align}
\label{eq:1}
Ax \leq b, \ \ b \in \R^m, \ A \in \R^{m\times n}, \ m \gg n.
\end{align}
In the last decade, projection-based iterative methods gain a considerable amount of traction for solving problem \eqref{eq:1}. Recent advances in the area of projection methods suggest that one can interpret most projection-based methods under one big family of methods widely known as \textit{Sketch \& Project} (SP) methods \cite{gower:2015}. The SP framework connects several projection-based methods such as \textit{Randomized Newton}, \textit{Randomized Kaczmarz},  and \textit{Randomized Coordinate Descent}, \textit{Random Gaussian Pursuit} and \textit{Randomized Block Kaczmarz} etc. Although the original setup of SP methods was proposed for solving a system of linear equations, recently this has been extended to a wide array of methods such as \textit{Quasi-Newton} methods \cite{gower:2017}, \textit{Matrix Pseudo-inverse} \cite{gower2016linearly}, \textit{Randomized Subspace Newton} \cite{gower2019rsn} \textit{Newton-Raphson} \cite{yuan2020sketched} and the references therein. Interestingly in \cite{gower2016linearly}, the authors discussed how one can recover most \textit{Quasi-Newton}-type methods such as \textit{Davidon–Fletcher–Powell}, \textit{Powell-Symmetric-Broyden}, \textit{Bad Broyden}, and \textit{Broyden–Fletcher–Goldfarb–Shanno} from the proposed SP methods by choosing different sketching matrices and positive definite matrices. In \cite{necoara:2019}, the authors proposed a weaker version of SP methods for solving a wide range of convex feasibility problems. Recently, several generalized and accelerated variants of the SP method have been proposed in \cite{richtrik2017stochastic, loizou:2017,NIPS:2018, gower2019adaptive} for solving a system of linear equations.

One key ingredient for the sketching methods is the choice of selecting sketching matrices at each iteration. Several important sampling strategies used widely in the broader sense of Kaczmarz method are \textit{Uniform Sampling} \cite{strohmer:2008, lewis:2010,wright:2016}, \textit{Maximum Distance Sampling} \cite{motzkin}, \textit{Kaczmarz Motzkin Sampling} \cite{haddock:2017, haddock:2019, Morshed2019, morshed2020generalization,morshed:momentum}, Capped Sampling \cite{bai:2018,gower2019adaptive}. Especially in their recent work, Gower \textit{et. al} discussed the above-mentioned sampling strategies in the context of SP methods for solving a system of linear equations. They showed that by sampling indices based on the \textit{Sketched Loss} of the current iterate, one can design efficient algorithms. Before we delved into the contributions of our work, first let us provide some background information of some classical and modern algorithmic developments over the years for solving problem \eqref{eq:1}.

The most well-known and simplest of the projection methods is the Kaczmarz method \cite{kaczmarz:1937}. Kaczmarz method is a variant of the SP method which updates the next iterate as: $x_{k+1} = \mathcal{P}_{\mathcal{X}_i}(x_k)$ \footnote{$\mathcal{P}_{\mathcal{X}_i}(x_k)$ denotes the orthogonal projection of $x_k$ onto the hyper-plane $\mathcal{X}_i$. Note that, Kaczmarz method uses unit coordinate vector as the sketching vector at each iteration.}. The research into Kaczmarz-type methods boomed in the last decade after Strohmer \textit{et. al} \cite{strohmer:2008} proposed the RK method which
significantly improves the theoretical and practical efficiency \footnote{Instead of selecting the hyper-plane $\mathcal{X}_i$ by cyclic method \cite{kaczmarz:1937}, they proposed to select hyper-plane randomly.}. Another classical method is the so-called MR method \cite{agamon,motzkin} that chooses hyper-plane $\mathcal{X}_i$ with the maximum positive residual \footnote{The perceptron algorithm in machine learning \cite{ramdas:2014,ramdas:2016} can be sought as a variant of the MR-type method.}. The work Strohmer \textit{et. al} \cite{strohmer:2008} spurred various extensions of the RK method for solving various type of problems such as linear system, linear feasibility, least square etc. (see \cite{lewis:2010,needell:2010,eldar:2011,zouzias:2013,lee:2013,NEEDELL:2014,agaskar:2014,ma:2015,NEEDELL:2015, blockneddel:2015,needell:2016, wright:2016,needell:2016,haddock:2017, greedbai:2018,haddock:2019, Morshed2019,razaviyayn:2019, haddock:2019,needell2019block, morshed2020generalization,morshed:momentum} and the references therein).

In the last decade, a huge amount of optimization and machine learning works have been associated with developing efficient \& accelerated iterative methods. The two standout acceleration techniques are the so-called \textit{Polyak Momentum} \cite{polyak1964some} and \textit{Nesterov Accelerated Gradient} \cite{nesterov:1983} that are widely used in training deep neural network learning. These methods roots back to developing an efficient version of the \textit{Gradient Descent} (GD) method for solving the unconstrained minimization problem. In the context of projection-based iterative methods, these acceleration techniques have been incorporated to various method such as \textit{Coordinate Descent} \cite{nesterov:2012}, \textit{Randomized Kaczmarz} \cite{wright:2016}, \textit{Sketch \& Project} \cite{loizou:2017}, \textit{Affine Scaling} \cite{morshed:2018}, \textit{Quasi-Newton} \cite{NIPS:2018}, \textit{Randomized Gossip} \cite{peter:2019}, \textit{Sampling Kaczmarz Motzkin} \cite{Morshed2019, morshed2020generalization, morshed:momentum} etc. (for more details please see the references therein).


In this work, we first propose two greedy sampling techniques that generalize the available sampling strategies and generate efficient algorithmic variants of the SP method for solving the linear feasibility problem. We extend available greedy techniques such as Kaczmarz-Motzkin sampling \cite{haddock:2017,haddock:2019,morshed2020generalization,morshed:momentum} and \textit{Capped Kaczmarz} \cite{greedbai:2018,gower2019adaptive} to the SP framework for solving the LF problem of \eqref{eq:1}. Furthermore, we introduce the heavy ball momentum scheme to the proposed greedy SP method to accelerate the efficiency. The proposed greedy rule-based SP methods outperform the available SP methods. Moreover, the momentum variants significantly outperform the proposed greedy algorithms on a wide variety of test instances in terms of CPU consumption time and solution quality \footnote{Note that with some modifications to the proposed methods one can design efficient algorithms for solving linear feasibility problems with both equality and inequality equations.}.

\subsection{Outline}
The paper is organized as follows. In section \ref{sec:contr}, we provide a brief summary of existing projection-based methods that deal with solving the feasibility problem of \eqref{eq:1}. At the end of section \ref{sec:contr}, we list a summary of the important contributions of this work. In section \ref{sec:tech}, we provide some technical backgrounds \& tools to handle the analysis of the proposed methods. In section \ref{sec:algorithms}, we provide the proposed algorithms. At the end of section \ref{sec:algorithms}, we try to provide a visual representation of the proposed methods. In section \ref{sec:sr}, we discuss the proposed greedy sampling rules along with their algorithmic influence on the function $f(x)$. The main convergence results for the basic method and the momentum variants are provided in section \ref{sec:conv}. In section \ref{sec:num}, we carry out extensive numerical experiments to measure the performance of the proposed greedy sampling rules and momentum variants. The paper is concluded in section \ref{sec:conclusions} with remarks and future research directions. In Appendix \ref{sec:prel}, we mention some preliminary results we borrow from the literature. In Appendix \ref{sec:prf}, we provide the necessary proofs of the proposed technical results. In Appendix \ref{appendix:exp}, we provide some extra experimental figures.

\subsection{Notation}
We follow the standard linear algebra notation throughout the paper. The notation $\mathbb{R}^{m\times n}$ will be used to denote the set of $m \times n$ real-valued matrices. Similarly, $\mathbb{R}^{m\times n}_+$ will be used to denote the set of $m \times n$ real-valued non-negative matrices. The feasible region of the LF problem defined in \eqref{eq:1} is given by $\mathcal{X} = \{x \in \R^n | \ Ax \leq b\}$. Similarly, the $i^{th}$ hyper-plane $\mathcal{X}_i$ of the feasibility problem is defined as $\mathcal{X}_i = \{x \in \R^n | \ a_i^Tx \leq b\}$, where the notation $a_i^T$ denotes the rows of matrix $A$. For a real-valued matrix $A$, the notation $\|A\|$  and $\|A\|_F$ denotes respectively the spectral and the Frobenius norm. The notation $x^+$ will be used to denote the positive part of any real number $x \in \R$, i.e., $x^+ = \max \{x,0\}$. For any two arbitrary matrices $M,\ N$, the notation $M \succ N$ defines the positive definiteness of the matrix $M-N$. Given any sampling rule $\mathcal{R}$, by which the index $i$ will be chosen, we use the notation $\E[\cdot \ | \ i \sim \mathcal{R} ]$ to denote the expectation with respect to the sampling rule $\mathbb{R}$. Let, $B \in \R^{n \times n}$ be any positive definite matrix. We denote the inner product equipped with the $B$ matrix as $\langle x, B x \rangle = x^TBx= \|x\|_B $.  For a closed convex set $\emptyset \neq \mathcal{X} \subseteq \R^n$, the notation $\mathcal{P}_{\mathcal{X}}^{B}(x)$ denotes the projection operator onto $\mathcal{X}$, in the $B-$norm, i.e., $\mathcal{P}_{\mathcal{X}}^{B}(x) = \argmin_{y \in \R^n} \{\|y-x\|_B : \ y \in \mathcal{X}\} $. Also, we define the distance from $x$ to the set $\mathcal{X}$ in the $B-$norm as $ d_B(x,\mathcal{X}) = \inf_{y \in \mathcal{X}} \|y-x\|_B  =  \| x- \mathcal{P}_{\mathcal{X}}^{B}( x) \big \|_B$. Similarly, the notation $d(x,\mathcal{X})$ denotes the Euclidean distance from $x$ to the set $\mathcal{X}$, i.e., $ d(x,\mathcal{X}) = \inf_{y \in \mathcal{X}} \|y-x\|  =  \| x- \mathcal{P}_{\mathcal{X}}( x) \big \|$, where the notation $\mathcal{P}_{\mathcal{X}}( x) $ denotes the orthogonal projection of $x \in \R^n $ onto the feasible region $\mathcal{X}$.

\section{Preliminaries \& Our Contributions}
\label{sec:contr}
In this section, we first discuss some preliminary works that deal with solving the linear feasibility problem of \eqref{eq:1}.

\paragraph{Sketch \& Project Methods} In \cite{necoara:2019}, the following method is proposed for solving the linear feasibility problem of \eqref{eq:1}: given a random iterate $x_k$, the goal of \textit{Sketch \& Project Methods} is to seek the closest point $x_{k+1}$ such that $x_{k+1}$ solves the following sketched feasibility problem:
\begin{align}
\label{eq:sk_siam}
x_{k+1} =     \argmin_x \|x-x_k\|^2 \quad \text{subject to} \quad  S^TAx \leq  S^T b,
\end{align}
where, $S \in \R^m_+$ is selected randomly from distribution $\mathcal{D}$. The solution of the sketching problem of \eqref{eq:sk_siam} is given by
\begin{align}
\label{eq:sk1_siam}
   x_{k+1}  = x_k -  \frac{\left[S^T(Ax_k-b)\right]^+}{\|A^TS\|^2} A^T S.
\end{align}
The authors specifically discussed the following case: choose $S = e_i$ with $i$ chosen with probability $\frac{\|a_i\|^2}{\|A\|^2_F}$. Furthermore, they showed that the proposed method converges linearly given that the so-called regularity condition holds \footnote{Please see Lemma \ref{neco}.}.

\paragraph{Momentum Sampling Kaczmarz Motzkin (MSKM)} In \cite{morshed:momentum}, the authors proposed the following update formula for solving the feasibility problem \eqref{eq:1}:
\begin{align}
\label{mskm:1}
& x_{k+1}  = x_k - \delta  \frac{\left(a_{i^*}^Tx_k -b_{i^*}\right)^+}{\|a_{i^*}\|^2} a_{i^*} + \gamma(x_k -x_{k-1}),
\end{align}
where $ \delta > 0$ is the projection parameter and $\gamma \geq 0$ is the momentum parameter \footnote{Note that, this method can be sought as momentum extension of the methods proposed in \cite{haddock:2017,haddock:2019}.}. The index $i^*$ at iteration $k$ is selected by the following rule: $i^* = \argmax_{i \in \phi_k(\tau)} \{a_i^Tx_k-b_i, 0\}$. Here, the set $\phi_k(\tau)$ denotes the set consisting of $\tau$ rows uniformly sampled from the rows of matrix $A$. One can easily recover momentum variants of the RK and MR methods by choosing sample size $\tau$ as $\tau =1$ and $\tau = m$ respectively. Note that the momentum induced RK method takes too many iterations too converge (cheaper per iteration cost) and the momentum induced MR method has a higher per iteration cost (takes few iterations). By introducing this specific sampling rule MSKM method enjoys better performance as shown in \cite{morshed:momentum}.

\paragraph{Heavy Ball Momentum.} Polyak momentum, popularly known as heavy ball momentum is one of the most oldest and important acceleration techniques for solving unconstrained minimization problem: $x^{*} = \argmin_{x \in \mathbb{R}^n} \mathcal{F}(x)$. The heavy ball update for solving the above problem is given by:
\begin{align*}
    x_{k+1} = x_k - \alpha_k \nabla \mathcal{F}(x_k) + \gamma (x_k-x_{k-1}),
\end{align*}
where, $\gamma$ is the momentum parameter. When, $\gamma = 0$, this method resolves into the so-called GD method. Polyak \cite{polyak1964some} proved that for twice continuously differentiable function $\mathcal{F}(x)$ with $\mu$ strong convexity constant and $L-$Lipschitz gradient, the momentum GD method achieves accelerated rate (with appropriate step-size parameters $\alpha_k$ and momentum parameter $\gamma$). Building on the above-mentioned works, in this paper, we propose two greedy sampling strategy based Sketch \& Project methods and the corresponding momentum variants for solving the linear feasibility problem of \eqref{eq:1}. In the following, we provide a brief summary of the contributions of this work.

\subsection{Summary of Our Contributions}
\paragraph{Adaptive \textit{Sketch \& Project} method with greedy sampling strategies.} We generalize the \textit{Sketch \& Project} methods by introducing a new parameter \footnote{This is standard in the linear system framework \cite{gower:2015,richtrik2017stochastic}.}, i.e., positive definite matrix $B \in \R^{n \times n}$. We introduce the \textit{Greedy Sampling} rule that generalizes several available sampling strategies such as uniform sampling, maximum distance sampling. Moreover, we show that the \textit{Greedy Sampling} rule produces more efficient algorithms than both the above-mentioned sampling rules. Furthermore, we introduce the \textit{Greedy Capped Sampling} rule that extends the scope of the so-called capped sampling strategy to make the resulting algorithm much more efficient.

\paragraph{Adaptive \textit{Sketch \& Project} method with momentum.}
We propose heavy ball momentum techniques to the developed adaptive sketching methods. The proposed momentum algorithms outperform the basic sketching methods for the majority of test instances. Furthermore, one can recover a variety of momentum methods and their convergence results from our convergence results for solving linear feasibility problems. In Table \ref{tab:AK}, we provide three variants of momentum algorithms that can be obtained form our proposed methods \footnote{Note that for simplicity, we decide to provide variants of Kaczmarz and CD methods with orthogonal projection. One can derive a wide variety of methods by choosing a different combination of sketching vector, matrix $B$, projection parameter $\delta$, and sampling rule.}.
\begin{table}[H]
\centering
\caption{\textit{Adaptive Sketch \& Project} methods with momentum for solving problem \eqref{eq:1}.}
\adjustbox{max width=\textwidth}{
\begin{tabular}{|c|c|c|c|c|}
\hline
$S_i $ & $B$ & Sampling Rule, $q = m$ & $x_{k+1}$ & Algorithm \\ \hline
\multirow{4}{*}{\begin{tabular}[c]{@{}c@{}} $e_i$\end{tabular}} & \multirow{4}{*}{$I$} & $\mathbb{P}(i) = \frac{\|a_i\|^2}{\|A\|^2_F}$ & \multirow{4}{*}{\begin{tabular}[c]{@{}c@{}}$x_k + \gamma(x_k-x_{k-1}) $\\ $- \delta \frac{\left(a_{i}^Tx_k -b_{i}\right)^+}{\|a_{i}\|^2} a_{i}  $ \end{tabular}} & MRK \cite{morshed:momentum} \\ \cline{3-3} \cline{5-5} 
 &  & $i = \argmax_{j} \frac{|(a_j^Tx_k-b_j)^+|^2}{\|a_j\|^2} $ &  & MMR \cite{morshed:momentum} \\ \cline{3-3} \cline{5-5} 
 &  & $i = \argmax_{j \in \phi_k(\tau)} \frac{|(a_j^Tx_k-b_j)^+|^2}{\|a_j\|^2} $ &  & SKM, MSKM \cite{haddock:2017,morshed:momentum} \\ \cline{3-3} \cline{5-5}
 &  & Capped Sampling ($i \in \mathcal{W}_k$) &  & MCK, New \\ \hline
 $S_i$ & $B$ & Sampling Rule, $q = m$ & $x_{k+1}$ & Algorithm  \\ \hline
\multirow{4}{*}{\begin{tabular}[c]{@{}c@{}} $e_i$\end{tabular}} & \multirow{4}{*}{$A \succ 0$} &$\mathbb{P}(i) = \frac{A_{ii}}{\textbf{Tr}(A)}$ & \multirow{4}{*}{\begin{tabular}[c]{@{}c@{}}$x_k + \gamma(x_k-x_{k-1}) $\\ $- \delta \frac{\left(a_{i}^Tx_k -b_{i}\right)^+}{A_{ii}} e_{i} $ \end{tabular}} & MRCD, New \\ \cline{3-3} \cline{5-5} 
 &  & $i = \argmax_{j} \frac{| (a_{j}^Tx_k-b_{j})^+ |^2}{A_{jj}^2} $ &  & MMCD, New \\ \cline{3-3} \cline{5-5} 
 &  & $i = \argmax_{j \in \phi_k(\tau)} \frac{| (a_{j}^Tx_k-b_{j})^+ |^2}{A_{jj}^2}  $ &  & MSCD, New \\ \cline{3-3} \cline{5-5} 
 &  & Capped Sampling ($i \in \mathcal{W}_k$) &  & MCCD, New \\ \hline
 $S_i$ & $B$ & Sampling Rule, $q = n$, $A \geq 0$ & $x_{k+1}$ & Algorithm  \\ \hline
\multirow{4}{*}{\begin{tabular}[c]{@{}c@{}} $A^+_i$\end{tabular}} & \multirow{4}{*}{$A^TA^+ $} &$\mathbb{P}(i) = \frac{\|A_i^+\|^2}{\|A^+\|^2_F}$ & \multirow{4}{*}{\begin{tabular}[c]{@{}c@{}}$x_k + \gamma(x_k-x_{k-1}) $\\ $- \delta \frac{\left(A_i^{+T}(Ax_k -b)\right)^+}{\|A^+_{i}\|^2} e_{i}  $  \end{tabular}} & MRCD-LS, New \\ \cline{3-3} \cline{5-5} 
 &  & $i = \argmax_{j} \frac{|\left(A_j^{+T}(Ax_k -b)\right)^+|^2}{\|A^+_{j}\|^2} $ &  & MMCD-LS, New \\ \cline{3-3} \cline{5-5} 
 &  & $i = \argmax_{j \in \phi_k(\tau)}\frac{|\left(A_j^{+T}(Ax_k -b)\right)^+|^2}{\|A^+_{j}\|^2}  $ &  & MSCD-LS, New \\ \cline{3-3} \cline{5-5} 
 &  & Capped Sampling ($i \in \mathcal{W}_k$) &  & MCCD-LS, New \\ \hline
\end{tabular}}
\label{tab:AK}
\end{table}

\paragraph{Global linear rate:}
We study the convergence behavior of the proposed adaptive sketching method as well as the momentum induced adaptive sketching method in great detail. We establish a global linear rate for both methods. We show that the terms $\E[d_B(x_{k}, \mathcal{X})^2]$ and $\E[f(x_k)]$ converge for a wide range of projection parameters $0 < \delta < 2$ and momentum parameter $\gamma \geq 0$. Our result connects several well-known convergence results with respect to sampling rules, positive definite matrix $B$, and sketching vectors $S_i \geq 0$. 

\paragraph{Sub-linear rate:}
For a fair understanding, we show that the Cesaro average of iterate, i.e., $\Tilde{x}_k = \frac{1}{k} \sum \limits_{i =0}^k x_i $ generated by the basic and momentum methods enjoys $\mathcal{O}(\frac{1}{k})$ sub-linear rate. One can obtain various well-known Cesaro average results from our proposed result.

\paragraph{Certificate of feasibility}
If there exists a point $x^*$ such that $\theta(x_k) < 2^{1-\sigma}$, then this point will be called a certificate of feasibility for the rational system $Ax \leq b$ (see Lemma \ref{lem:skm4}). When the feasibility problem $Ax \leq b$ is feasible, it is of practical benefit to find a certificate of feasibility after finitely many iterations. Moreover, if one fails to obtain a feasibility certificate after finitely many iterations, one needs to provide a lower bound on the probability that the system is infeasible. Assuming the system if feasible, we obtain an upper bound on the probability of finding a certificate of feasibility for the momentum induced adaptive sketching method (see Theorem \ref{th:4}). Our certificate of feasibility result extends the results obtained in \cite{haddock:2017} for the SKM method and MSKM method of \cite{morshed:momentum}. Moreover, from our result, one can show certificate of feasibility results for many new methods.

\section{Technical Tools}
\label{sec:tech}

In this section, we discuss some preliminary results that we will frequently use throughout the paper. First, we start by providing the assumptions of this work. Then we introduce function $f(x)$, variants of which are frequently used in the literature for analyzing the behavior of \textit{Sketch \& Project} methods \cite{richtrik2017stochastic, loizou:2017,necoara:2019}. For completeness, we discuss some useful results in Appendix \ref{sec:prel} that we borrowed from the literature \cite{haddock:2017,morshed:momentum,morshed2020generalization}. These results are instrumental for the convergence analysis of projection-based iterative methods in solving the linear feasibility problem of \eqref{eq:1}

\subsection{Assumptions}
Throughout the paper, we will assume that the following assumptions hold: (1) the system $Ax \leq b$ is consistent, and (2) matrix $A$ has no zero rows. 

\subsection{Function $f(x)$}
In this subsection, we introduce the function $f(x)$. Before we delved into the definition, first let us formalize the sampling rule. At iteration $k$, given any random iterate $x_k$, in deriving the next update $x_{k+1}$ the proposed algorithms have to choose an index $i$. The index $i$ is chosen following some sampling rule \footnote{Throughout the paper, we use $\mathcal{R}$ to denote any generic sampling rule. Later, in section \ref{sec:sr}, we will discuss the proposed greedy sampling rules.} $\mathcal{R}$, i.e., $i \sim \mathcal{R}$. We use $\E[\cdot \ | \ i \sim \mathcal{R}]$ to denote the resulting expectation. 

\begin{mdframed}[backgroundcolor=gray!20,   topline=false,   bottomline =false,   rightline=false,   leftline=false] \begin{definition}
Let's define function $f(x)$ as follows:
\begin{align}
\label{eq:sk3}
 f(x) = \E[f_i(x) \ | \ i \sim \mathcal{R}], \quad   f_{i}(x) = \frac{\big |\left[S_{i}^T(Ax-b)\right]^+\big |^2}{2\|A^TS_{i}\|^2_{B^{-1}}},
\end{align}
where, $B \in \R^{n \times n}$ is a positive definite matrix and $S_i \in \R^m_+$ is the $i^{th}$ sketching vector from the sketched vector set $\mathcal{S}(q) $, i.e., $\mathcal{S}(q) = \{S_1,S_2,...,S_q\}$.
\end{definition}
\end{mdframed}

The function $f_i(x)$ is the so-called sketched loss \cite{gower2019adaptive}. Note, that, the gradients $\nabla f_i$ and $\nabla^B f_i$ of function $f_i$ are given by
\begin{align}
\label{eq:sk4}
  \nabla f_i(x) = \frac{\left[S_{i}^T(Ax-b)\right]^+}{\|A^TS_{i}\|^2_{B^{-1}}}A^T S_{i}, \quad   \nabla^{B} f_i(x) = \frac{\left[S_{i}^T(Ax-b)\right]^+}{\|A^TS_{i}\|^2_{B^{-1}}}  B^{-1} A^T S_{i},
\end{align}
where, $\nabla^B f_i$ denotes the gradient of $f_i$ with respect to the $B-$ norm. Denote, $\mathcal{X}_{S_i} = \{x \ | \ S_i^TAx \leq S_i^Tb\}$, this is the so-called sketched feasible region. Consider the following feasibility problem:
\begin{align}
    \label{reform:feas}
    \text{Find} \quad x \in \mathcal{X}' = \bigcap\limits_{{i \in \{1,2,...,q\}}} \mathcal{X}_{S_i}.
\end{align}
It can be easily show that $\mathcal{X} \subseteq \mathcal{X}' $. the above convex feasibility problem can be reformulated as the following stochastic optimization problem:
\begin{align}
    \label{prob:stoc}
    x = \argmin f(x) = \argmin \E[f_i(x) \ | \ i \sim \mathcal{R}].
\end{align}
In the following Lemma, we will show that the problems \eqref{eq:1} and \eqref{reform:feas} are equivalent, i.e., $\mathcal{X} = \mathcal{X}' $ if the following property holds \footnote{This equivalence has been shown in \cite{necoara:2019} for the case of $B=I$.}.

\begin{mdframed}[backgroundcolor=gray!20,   topline=false,   bottomline =false,   rightline=false,   leftline=false]  \begin{lemma}
\label{neco}
Exactness holds, i.e., $\mathcal{X} = \mathcal{X}' $ provided that there exists a constant $\mu > 0$ such that the following identity holds:
\begin{align}
\label{stoc:reg}
    \mu \ d_B(x,\mathcal{X})^2 \leq \E [d_B(x,\mathcal{X}_{S_i})^2 \ | \ i \sim \mathcal{R}].
\end{align}
\end{lemma} \end{mdframed}

\begin{proof}
See Appendix \ref{pstoc:reg}
\end{proof}
This is the so-called \textit{Stochastic Linear Regularity Property} defined in \cite{necoara:2019}. Later, in section \ref{sec:sr}, we will show that the proposed greedy sampling rules always enjoy this property. In other words, in our proposed greedy sampling setting the problems \eqref{1}, \eqref{reform:feas}, and \eqref{prob:stoc} are equivalent.

\section{Algorithms}
\label{sec:algorithms}
In this section, we first propose a \textit{Sketch \& Project} framework that is equipped with specific sampling rule $\mathcal{R}$ and a positive definite matrix $B \in \R^{n \times n}$. The proposed method generalizes the method proposed in \cite{necoara:2019} and allow us to design efficient algorithms based on greedy sampling strategies. Then, we propose a momentum variant of the proposed \textit{Sketch \& Project} method that generalizes several existing works for solving linear feasibility problems.

\subsection{Adaptive Sketch \& Project (ASP)}
Now, we discuss the proposed \textit{Adaptive Sketch \& Project} (ASP) method for solving the linear feasibility problem. The main ingredients for the ASP algorithm include: 1) A positive definite matrix $B$, 2) Sketching vector set $ \mathcal{S}(q)$, 3) Sketching Rule $\mathcal{R}$, and 4) Projection parameter $\delta$. At iteration $k$ given a current iterate $x_k$, we seek to find the closest point $x_{k+1}$ such that $x_{k+1}$ satisfies the following sketched feasibility problem:
\begin{align}
\label{eq:sk}
x_{k+1} =     \argmin_x \|x-x_k\|^2_B \quad \text{subject to} \quad  S_{i}^TAx \leq  S_{i}^T b,
\end{align}
where, the sketched vector $S_i \in \R^m_+$ is selected from the set of sketching vectors $ \mathcal{S}(q)$ following the sampling rule $\mathcal{R}$, i.e., $i \sim \mathcal{R}$ \footnote{Note that index $i$ is depended on iteration $k$, i.e., $i = i_k$. For ease of presentation, we denote the index as $i$.}. The solution of the sketching problem of \eqref{eq:sk} is given by
\begin{align}
\label{eq:sk1}
   x_{k+1}  = x_k -  \frac{\left[S_{i}^T(Ax_k-b)\right]^+}{\|A^TS_{i}\|^2_{B^{-1}}} B^{-1} A^T S_{i}.
\end{align}
If we allow projection parameter $\delta \in (0,2)$, this update formula becomes the following
\begin{align}
\label{eq:sk2}
   x_{k+1}  = x_k - \delta \frac{\left[S_{i}^T(Ax_k-b)\right]^+}{\|A^TS_{i}\|^2_{B^{-1}}} B^{-1} A^T S_{i} = x_k - \delta \ \nabla^{B} f_i(x_k).
\end{align}
Then, we get the following algorithm:

\begin{algorithm}
\caption{ASP Algorithm: $x_{k+1} = \textbf{ASP}(A,b,x_0,\delta, B, \mathcal{S}(q),  \mathcal{R}, K)$}
\label{alg:arsp}
\begin{algorithmic}
\STATE{Choose initial point $x_0 \in \R^n$}
\WHILE{$k \leq K$}
\STATE{From the sketched vector set $\mathcal{S}(q)$, select index $i$ based on sampling rule $\mathcal{R}$, i.e.,  $i \in \mathcal{R}$. Then update
\begin{align*}
& x_{k+1}  = x_k - \delta \ \nabla^{B} f_{i}(x_k);
\end{align*}
$k \leftarrow k+1$;}
\ENDWHILE
\RETURN $x$
\end{algorithmic}
\end{algorithm}

\begin{mdframed}[backgroundcolor=gray!20,   topline=false,   bottomline =false,   rightline=false,   leftline=false]  \begin{remark}
\label{rem:hist}
Algorithm \ref{alg:arsp} can be interpreted as the \textit{Stochastic Gradient Descent} (SGD) method with fixed learning rate $\delta$ for solving the stochastic optimization problem of \eqref{prob:stoc}. Note that, if we take take $B = I$ without any specific sampling rule $\mathcal{R}$. Then the proposed ASP algorithm resolves into the algorithm proposed in \cite{necoara:2019} (see equations \eqref{eq:sk_siam} and \eqref{eq:sk1_siam}). This work can also be sought as an extension of the work \cite{gower2019adaptive} for solving linear feasibility problems.
\end{remark} \end{mdframed}

\subsection{Adaptive Sketch \& Project with Momentum (ASPM)}
In this subsection, we discuss the ASP algorithm with momentum. The main ingredients for the ASPM algorithm include: 1) A positive definite matrix $B$, 2) Sketching vector set $ \mathcal{S}(q)$, 3) Sketching Rule $\mathcal{R}$, 4) Projection parameter $\delta$, and 5) Momentum parameter $\gamma$. Introducing the momentum scheme into the ASP algorithm, we can derive the following update formula:
\begin{align}
\label{eq:momupdate}
   x_{k+1} = x_k - \delta \ \nabla^{B} f_i(x_k) + \gamma (x_k-x_{k-1}).
\end{align}
Then, we get the following algorithm:

\begin{algorithm}
\caption{ASPM Algorithm: $x_{k+1} = \textbf{ASPM}(A,b,x_0,\delta, \gamma, B, \mathcal{S}(q),  \mathcal{R}, K)$}
\label{alg:marsp}
\begin{algorithmic}
\STATE{Choose initial point $x_0 \in \R^n$}
\WHILE{$k \leq K$}
\STATE{From the sketched vector set $\mathcal{S}(q)$, select index $i$ based on sampling rule $\mathcal{R}$, i.e.,  $i \in \mathcal{R}$. Then update
\begin{align*}
& x_{k+1}  = x_k - \delta \ \nabla^{B} f_{i}(x_k) + \gamma (x_k-x_{k-1});
\end{align*}
$k \leftarrow k+1$;}
\ENDWHILE
\RETURN $x$
\end{algorithmic}
\end{algorithm}

\subsection{Visualization of the momentum}
In this subsection, we provide a visual illustration of the momentum mechanism with respect to uniform and maximum sampling rules. We provide pictorial representations and compare the momentum update with the corresponding basic update for the \textit{Greedy Kaczmarz} (GK) method.

\begin{figure}[htbp]
\centering
    \includegraphics[scale = 0.65]{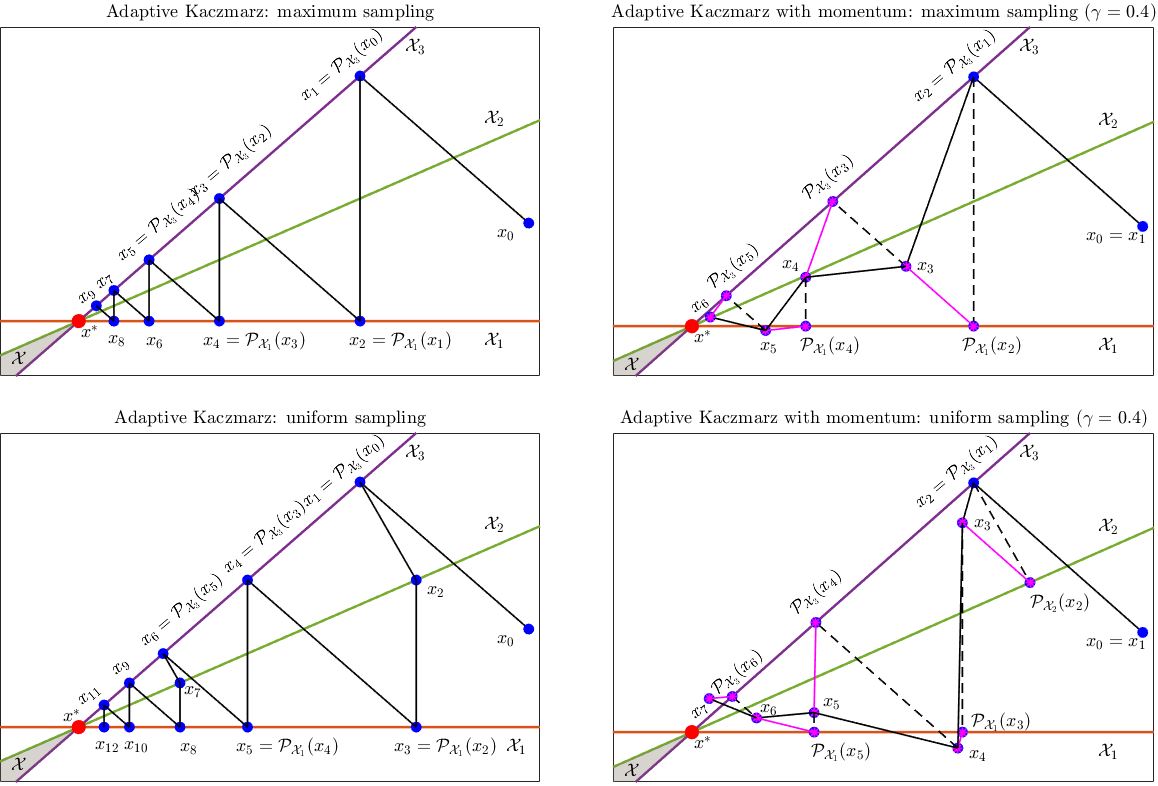}
    \caption{Graphical interpretation of the basic method and the momentum method with three hyper-planes $\mathcal{X}_j = \{x | a_j^Tx \leq b_j\}, \ j =1, 2, 3$. Shaded region $\mathcal{X}$ is the feasible region, top panel: GK with the maximum distance rule, bottom panel: GK with the uniform rule.}
    \label{fig:g1}
\end{figure}

In Figure \ref{fig:g1}, we draw several updates of the proposed methods in a $\R^2$ plane starting with the same initial point $x_0$. For ease of illustration, we select three hyper-planes $\mathcal{X}_1$, $\mathcal{X}_2$, and $\mathcal{X}_3$ with orthogonal projection, i.e., $\delta =1$. Given a random point $x_k$, the basic method finds the next update $x_{k+1}$ by the projection step that projects $x_k$ onto one of the hyper-planes to find the next update $x_{k+1}$. The notation $\mathcal{P}_{\mathcal{X}_1}(x_k)$ denotes the orthogonal projection of $x_k$ onto the hyper-plane $\mathcal{X}_1$. For the momentum variants, the extra momentum term, $\gamma(x_k-x_{k-1})$ is added to the projection step to find the next update. From Figure \ref{fig:g1}, we find that the momentum induced update $x_{k+1}$ is moves closer to the feasible region $\mathcal{X}$ then the update without momentum for both sampling rules \footnote{We will validate this statement later in the numerical experiments section by performing this comparison for a wide variety of large test instances.}. Furthermore, it can be noted that no matter which sampling rule we use the vector $x_{k+1}-\mathcal{P}_{\mathcal{X}_i}(x_k)$ is always parallel to the vector $x_{k}-x_{k-1}$ at any iteration $k \geq 1$ (please see Figure \ref{fig:g1}, the magenta-colored lines in the right panel sub-figures denote the vector $x_{k+1}-\mathcal{P}_{\mathcal{X}_i}(x_k)$, and the black colored lines in the right panel sub-figures denote the vector $x_{k}-x_{k-1}$).


\section{Function $f(x)$ \& Greedy Sampling Rules}
\label{sec:sr}
In this section, we discuss the proposed sampling rules and the corresponding results of the accompanying function $f(x)$. Previously, we defined $f(x) = \E[f_i(x) \ | \ i \sim \mathcal{R}]$, where $\mathcal{R}$ is a generic sampling rule. In this section, we will specifically discuss the implications of specific sampling rules on the properties of $f(x)$. In the first subsection, we provide some generic properties of function $f(x)$ that are true irrespective of sampling rules. In the second subsection, we discuss the proposed greedy sketched loss sampling rule and the corresponding properties of function $f(x)$. In the third subsection, we discuss the proposed greedy capped sketched loss sampling rule and the corresponding properties of function $f(x)$. We also discuss some well-known special cases of the proposed sampling rules.

\subsection{Properties of function $f(x)$}
In this subsection, we will discuss properties of the function $f(x)$. The following Lemma holds:
\begin{mdframed}[backgroundcolor=gray!20,   topline=false,   bottomline =false,   rightline=false,   leftline=false]  \begin{lemma}
\label{1}
Assume, the index $i$ is selected as $i \sim \mathcal{R}$, where $\mathcal{R}$ is a generic sampling rule. Then the following identities hold:
\begin{enumerate}
    \item $   f_i(x) =  \frac{1}{2} \big \| \nabla^{B} f_i(x) \big \|^2_B $.
    \item $   f(x) =  \frac{1}{2} \E \left[\|\nabla^B f_i(x)\|_B^2 \ | \ i \sim \mathcal{R} \right]   =  \frac{1}{2} \E \left[d_B(x,\mathcal{X}_{S_i})^2 \ | \ i \sim \mathcal{R} \right]$.
    \item $\nabla f(x) = \E[\nabla f_i(x) \ | \ i \sim \mathcal{R} ], \ \ \nabla^B f(x) = \E[\nabla^B f_i(x) \ | \ i \sim \mathcal{R} ]$.
\end{enumerate} 
\end{lemma} \end{mdframed}

\begin{proof}
  See Appendix \ref{p1}.
\end{proof}

\begin{mdframed}[backgroundcolor=gray!20,   topline=false,   bottomline =false,   rightline=false,   leftline=false]  \begin{lemma}
\label{2}
Assume, $\bar{x} \in \mathcal{X}$ (i.e., $A\bar{x} \leq b$), then for any $x \in \R^n$, we have
\begin{enumerate}
    \item $ \big \langle  \bar{x}-x, \nabla^{B} f_i(x) \big \rangle_{B} \leq - 2 f_{i}(x) $.
    \item $ 2 f(x)  \leq  d_B(x,\mathcal{X}) \ \|\E[\nabla^{B} f_i(x)\ | \ i \sim \mathcal{R} ]\|_B$.
\end{enumerate}
\end{lemma} \end{mdframed}

\begin{proof}
  See Appendix \ref{p2}.
\end{proof}

\begin{mdframed}[backgroundcolor=gray!20,   topline=false,   bottomline =false,   rightline=false,   leftline=false]  \begin{lemma}
\label{lem:distance}
For any $x \in \R^n$ and $\bar{x} \in \mathcal{X}$, the following identity holds,
\begin{align*}
  d_B(x,\mathcal{X})^2 \ = \ \big \|x- \mathcal{P}^B_{\mathcal{X}}(x)\big\|^2_B \ \leq \ \|x-\bar{x} \|^2_B.
\end{align*}
\end{lemma} \end{mdframed}

\begin{mdframed}[backgroundcolor=gray!20,   topline=false,   bottomline =false,   rightline=false,   leftline=false]  \begin{lemma}
\label{th:upper}
Assume, any sampling rule $\mathbb{R}$, where index $i$ is chosen with probability $p_i$. Let, $\mu_2 \geq 0$ be the smallest constant satisfying the inequality:
\begin{align*}
  \big \| \E [\nabla^B f_i(x) \ | \ i \sim \mathcal{R}]\big \|_B^2 & \leq 2 \mu_2 \ \E[f_i(x) \ | \ i \sim \mathcal{R}].
\end{align*}
Then, with the definition $Z = \E \left[\frac{S_iS_i^T}{\|A^TS_i\|^2_{B^{-1}}}  \ | \ i \sim \mathcal{R}\right]$, we have
\begin{align}
    \mu_2 = \lambda_{\max} \left(B^{-\frac{1}{2}} A^T Z A B^{-\frac{1}{2}}\right) \leq 1 \quad \text{and} \quad f(x) \leq \frac{\mu_2}{2} d_B(x,\mathcal{X})^2.
\end{align}
\end{lemma} \end{mdframed}

\begin{proof}
See Appendix \ref{pupper}.
\end{proof}

\begin{mdframed}[backgroundcolor=gray!20,   topline=false,   bottomline =false,   rightline=false,   leftline=false]  \begin{lemma}
\label{lem:GRsketching2}
For any $x,y \in \R^n$, we have
\begin{align}
\langle  \E[\nabla^B f_i(x) \ | \ i \sim \mathcal{R}] , y-x \rangle_B  = \langle \nabla f(x), y-x \rangle \leq f(y) - f(x),
\end{align}
for any random iterate $x$.
\end{lemma} \end{mdframed}

\begin{proof}
  The function $f(x) = \E[ f_i(x) \ | \ i \sim \mathcal{R}] $ is convex. Therefore, the identity of the above Lemma follows from the convexity of $f(x)$.
\end{proof}
We note that the condition of Lemma \ref{lem:GRsketching2} is weaker than the traditional strong convexity, and it is also weaker than the essentially strong convexity condition defined in \cite{karimi:2016}. For instance, the essentially strong convexity requires the following identity:
\begin{align*}
    f(x)-f(y)  \leq \langle \nabla f(x), x-y \rangle - \frac{\epsilon}{2} d_B(x,y)^2, \quad \forall x,y, \ \text{s.t.} \ \mathcal{P}^B_{\mathcal{X}}(x) = \mathcal{P}^B_{\mathcal{X}}(y),  
\end{align*}
for some $\epsilon>0$. The above condition clearly implies Lemma \ref{lem:GRsketching2}. In our convergence analysis, we need to derive constant $\mu_1 \geq 0$ such that the following identity holds:
\begin{align}
\label{r:mu1}
 f(x) =\E[ f_i(x) \ | \ i \sim \mathcal{R}]    \geq \frac{\mu_1}{2} d_B(x,\mathcal{X})^2,
\end{align}
for any $x \in \R^n$ with any sampling rule $\mathcal{R}$. This is the so-called restricted secant inequality condition defined in \cite{karimi:2016} and is weaker than the essentially strong convexity. Next, we discuss the proposed greedy sampling rules and the respective spectral constants $\mu_1$ and $\mu_2$.

\subsection{Greedy Sketched Loss Sampling}
\label{ss}
Now, we will discuss a special kind of sampling which is frequently used in the literature to develop better performing Kaczmarz-type methods for solving linear feasibility problems. Choose a sample of $\tau$ sketching vectors, uniformly at random from the sketched vector set $\mathcal{S}(q)$. Denote the index set generated by the above sampling process as, $\phi(\tau)$. From these $\tau$ sketched vectors, choose $i = \argmax_{i \in \phi(\tau)} f_i(x)$, i.e.,
\begin{align}
\label{def:i1}
i = \argmax_{i \in \phi(\tau)} f_i(x) = \argmax_{i \in \phi(\tau)} \frac{\big |\left[S_{i}^T(Ax-b)\right]^+\big |^2}{2\|A^TS_{i}\|^2_{B^{-1}}}.
\end{align}
Now, we will discuss the expectation calculation with respect to this greedy sampling rule. First, let us sort the sketched losses $f_{i}(x)$ from smallest to largest for any random iterate $x$. Denote, $f_{\underline{\mathbf{i_j}}}(x)$ as the $(\tau+j)^{th}$ entry on the sorted list, i.e.,
\begin{align}
\label{eq:sampling}
  \underbrace{f_{\underline{\mathbf{i_0}}}(x)}_{\tau^{th}} \ \leq ... \leq \ \underbrace{f_{\underline{\mathbf{i_j}}}(x)}_{(\tau+j)^{th}} \ \leq ... \leq \ \underbrace{f_{\underline{\mathbf{i_{q-\tau}}}}(x)}_{q^{th}}.
\end{align}
Now, from the sorted sketched losses $f_{i}(x)$, if we randomly select any entry of the residual vector at any given iteration $k$ the probability that any sample is selected is $\frac{1}{\binom{m}{\tau}}$. Also, each sketched loss has an equal probability of selection. Using the above discussion with the list provided in equation \eqref{eq:sampling}, we have the following:
{\allowdisplaybreaks
\begin{align}
\label{def:exp}
 \E[f_i(x) \ | \ i \sim \mathcal{G}(\tau)] & = \frac{1}{ \binom{q}{\tau}} \sum\limits_{j = 0}^{q-\tau} \binom{\tau-1+j}{\tau-1} f_{\underline{\mathbf{i_j}}}(x) \nonumber \\
 & = \frac{1}{2\binom{q}{\tau}} \sum\limits_{j = 0}^{q-\tau} \binom{\tau-1+j}{\tau-1} \frac{\big | \left[S_{\underline{\mathbf{i_j}}}^T(Ax-b)\right]^+ \big |^2}{\|A^TS_{\underline{\mathbf{i_j}}}\|^2_{B^{-1}}},
 \end{align}}
with $\E[f_i(x) \ | \ i \sim \mathcal{G}(\tau)]$ denotes the required expectation corresponding to this specific sampling rule. This greedy approach allows us to combine two well-known adaptive sketching rules \footnote{Note that, with the choice $S= e_i, \ B = I, \ q =m$ this specific sampling resolves into the sampling related the so-called SKM method \cite{haddock:2017,haddock:2019,Morshed2019,morshed2020generalization,morshed:momentum}}. For instance, take $\tau = 1$ in \eqref{def:exp}, then we have
{\allowdisplaybreaks
\begin{align}
\label{def:expK}
 \E[f_i(x) \ | \ i \sim \mathcal{G}(\tau)] & = \frac{1}{ \binom{q}{1}} \sum\limits_{j = 0}^{q-1}  f_{\underline{\mathbf{i_j}}}(x) = \frac{1}{q} \sum\limits_{i = 1}^{q}  f_i(x).
 \end{align}}
This is the so-called uniform sketching rule. Similarly, take $\tau = q$ in \eqref{def:exp}, then we have
{\allowdisplaybreaks
\begin{align}
\label{def:expM}
 \E[f_i(x) \ | \ i \sim \mathcal{G}(\tau)] & = \frac{1}{ \binom{q}{q}} \sum\limits_{j = 0}^{q-q} \binom{q-1+j}{q-1} f_{\underline{\mathbf{i_j}}}(x) = f_{\underline{\mathbf{i_0}}}(x) = \max_{i \in \{1,2,...,q\}} f_i(x).
 \end{align}}
This is the so-called maximum distance sketching rule. Note that, for the choice $S_i = e_i, B = I, q = m$ this rule resolves into the so-called MR rule for the linear feasibility problem.

\begin{mdframed}[backgroundcolor=gray!20,   topline=false,   bottomline =false,   rightline=false,   leftline=false]  \begin{lemma}
\label{lem:mu}
Let, $x$ be a random iterate generated by the greedy sampling rule defined above, then
\begin{align*}
    \frac{\mu_1}{ 2} \ d_B(x,\mathcal{X})^2 \ \leq \ \E[f_i(x) \ | \ i \sim \mathcal{G}(\tau)] \ \leq \  \frac{\mu_2}{2} \ d_B(x,\mathcal{X})^2,
\end{align*}
with 
\begin{align}
\label{ssmu}
\mu_1 = \frac{1}{\sigma \omega_2}\min \left\{\frac{1}{q-\tau+1}, \frac{1}{q-s}\right\}, \    \mu_2 = \min \left\{1, \frac{\tau}{\omega_1 q} \lambda_{\max} \left(B^{-\frac{1}{2}} A^T R^TR A B^{-\frac{1}{2}}\right) \right\},
\end{align}
where $\omega_1 = \min_{i \in \{1,...,q\}} \|A^TS_{i}\|^2_{B^{-1}}, \ \omega_2 = \max_{i \in \{1,...,q\}} \|A^TS_{i}\|^2_{B^{-1}}$, $\sigma $ is the Hoffman constant, and $R = [S_1,...,S_q]^T \in \R^{q \times m}$. The quantity $s$ denotes the number of zero entries in the sketched residual vector $\left[R(Ax-b)\right]^+$ (i.e., $s = q- \|\left[R(Ax-b)\right]^+\|_0$, where $\|\cdot\|_0$ denotes the zero norm).
\end{lemma} \end{mdframed}

\begin{proof}
See Appendix \ref{pmu}.
\end{proof}

\subsection{Greedy Capped Sketched Loss Sampling}
\label{gcs}
Now, we will discuss a greedy version of capped sampling \cite{gower2019adaptive}. Assume, $x$ is any random iterate. Take, $0 \leq \theta \leq 1$ and two sampled sketching vectors of sizes $\tau_1$ and $\tau_2$ respectively uniformly at random \footnote{Sampling with replacement. Note that, one can extend this method to multiple sampled sketched vectors.}. Let,
\begin{align}
    \mathcal{W} = \left\{i  | \ f_i(x) \geq \theta  \E[f_i(x) \ | \ i \sim \mathcal{G}(\tau_1)] + (1-\theta) \E[f_i(x) \ | \ i \sim \mathcal{G}(\tau_2)] \right\}.
\end{align}
Then, select index $i \in \mathcal{W}$ with probability $p_i$ to update the next iterate \footnote{Note that, an obvious generalization of this sampling rule can be: $i \in \mathcal{W} = \{i \ | \ f_i(x) \geq \sum \limits_{j=1}^{N} \theta_j \E[f_i(x) \ | \ i \sim \mathcal{G}(\tau_j)]\}$ with $0 \leq \theta_j \leq 1, \ 1 \leq \tau_j \leq q$.}. With the choice $\tau_1 = q, \ \tau_2 = 1$ this resolves into the Capped sampling proposed in (\cite{bai:2018}). This set in not empty \footnote{The computation of the quantity $\E[f_i(x) \ | \ i \sim \mathcal{G}(\tau)]\}$ is not of practical choice for the implementation of the proposed Capped sampling. Instead, we suggest to use a reasonable lower bound, i.e., $\E[f_i(x) \ | \ i \sim \mathcal{G}(\tau)]\} \geq \frac{1}{\omega_2 q} \|[R(Ax-b)]^+\|^2 $, please see the proof of Lemma \ref{lem:mu} for details.} as 
\begin{align*}
    \theta  \E[f_i(x) \ | \ i \sim \mathcal{G}(\tau_1)] + (1-\theta) \E[f_i(x) \ | \ i \sim \mathcal{G}(\tau_2)] \leq \max_{i \in \{1,2,..,q\}} f_i(x),
\end{align*}
in other words, $\max_{i \in \{1,2,..,q\}} f_i(x) \in \mathcal{W}$. The resulting expectation can be calculated as follows:
\begin{align}
    \E[f_i(x) \ | \ i \sim \mathcal{C}(\theta, \tau_1,\tau_2)] & = \sum \limits_{j \in \mathcal{W}} p_j f_j(x).
\end{align}

\begin{mdframed}[backgroundcolor=gray!20,   topline=false,   bottomline =false,   rightline=false,   leftline=false]  \begin{lemma}
\label{lem:c}
Let, $x$ be a random iterate generated by the greedy capped sampling rule defined above, then
\begin{align*}
    \frac{\mu_1}{ 2} \ d_B(x,\mathcal{X})^2 \ \leq \ \E[f_i(x) \ | \ i \sim \mathcal{C}(\theta, \tau_1,\tau_2)] \ \leq \  \frac{\mu_2}{2} \ d_B(x,\mathcal{X})^2,
\end{align*}
with $\mu_1 = \theta \mu_1(\tau_1)+ (1-\theta) \mu_1 (\tau_2), \    \mu_2 = \mu_2(q)$. Here, $\mu_1(\tau)$ and $\mu_2(\tau)$ are the spectral constants obtained from greedy sampling rule with $\tau$ sample sketched vectors.
\end{lemma} \end{mdframed}

\begin{proof}
See Appendix \ref{pc}.
\end{proof}

\begin{mdframed}[backgroundcolor=gray!20,   topline=false,   bottomline =false,   rightline=false,   leftline=false]  \begin{remark}
\label{rem:strongconvex}
Lemmas \ref{lem:mu} and \ref{lem:c} state that the function $ f(x) = \E[f_i(x) \ | \ i \sim \mathcal{R}]$ is strongly convex and has Lipschitz continuous gradient when restricted along the segment $[x,\mathcal{P}^B_{\mathcal{X}}(x)]$ under both the greedy and greedy capped sampling rules. To show this, first note that, $f^* = \min_{x} \E[f_i(x) \ | \ i \sim \mathcal{R}] = 0$ for any $x$ such that $x \in \mathcal{X}$. Similarly, we have $\nabla f(\mathcal{P}^B_{\mathcal{X}}(x)) = 0$. Then the results from the preceding Lemmas can be rewritten as
\begin{align}
    & \frac{\mu_1}{ 2}  \|x-\mathcal{P}^B_{\mathcal{X}}(x)\|_B^2 + \big \langle \nabla f(\mathcal{P}^B_{\mathcal{X}}(x)),x-\mathcal{P}^B_{\mathcal{X}}(x) \big \rangle  \ \leq \ f(x) - f^*, \label{eq:strongl}\\
    & f(x) - f^* \ \leq \  \big \langle \nabla f(\mathcal{P}^B_{\mathcal{X}}(x)),x-\mathcal{P}^B_{\mathcal{X}}(x) \big \rangle + \frac{\mu_2}{2}\  \|x-\mathcal{P}^B_{\mathcal{X}}(x)\|_B^2. \label{eq:strongm}
\end{align}
\end{remark} \end{mdframed}
Without loss of generality, we denoted $\mu_1$ and $\mu_2$ as the respective spectral constants for the corresponding adaptive sampling rules provided earlier. Equation \eqref{eq:strongl} represent the Lipschitz continuity condition and equation \eqref{eq:strongm} represent the strong convexity condition along the line segment $[x,\mathcal{P}^B_{\mathcal{X}}(x)]$.

\section{Main Results}
\label{sec:conv}

In this section, we derive the convergence results for the proposed ASP and ASPM algorithms. For ease of analysis, we prove the results for any sampling rule $\mathcal{R}$ with generic spectral constants $\mu_1$ and $\mu_2$. In subsection \ref{sec:res}, we provide convergence results for both ASP and ASPM methods. In subsection \ref{sec:certificate}, we provide a probabilistic estimation of the certificate of feasibility for the proposed ASPM method from which we can recover several well-known results. In subsection \ref{sec:cesaro}, we provide the convergence result for the ASPM method with respect to the so-called Cesaro average. Finally, in subsection \ref{sec:special}, we discuss some well-known algorithms and their convergence results that can be derived from the proposed Theorems.

\subsection{Convergence Results for ASP \& ASPM Methods}
\label{sec:res}
In Theorem \ref{th:b1}, we prove the convergence results for the decay of $\E[ d_B(x_{k},\mathcal{X})^2]$ and $\E[f(x_{k})]$ generated by the ASP method. Similarly Theorems \ref{th:1} and \ref{th:mom2} are the respective results for the ASPM method. Before we delved into the technical results, first let us define the following sets:
\begin{align}
\label{eq:s}
& Q = \left\{(\delta, \gamma) \ | \ 0 < \delta < 2, \ 0 \leq \gamma <  \frac{1-\sqrt{h_{\mathcal{R}}(\delta)}}{1-\sqrt{h_{\mathcal{R}}(\delta)} + \delta \sqrt{\mu_2}}\right\}, \nonumber \\
  & R  = \left\{(\delta, \gamma,\zeta) \ | \ 0 < \delta < 2, \ \zeta \geq 0, \ 0 \leq \gamma < \frac{ \zeta}{1+ \zeta} \right\},  \\
  & S = \left\{(\delta, \gamma, \zeta) \ |  \ \frac{\gamma \mu_2}{\mu_1} <  \frac{2}{1+ \zeta}-\delta+\gamma\leq \frac{1+ \gamma}{\delta \mu_1(1+ \zeta)} \right\}. \nonumber 
\end{align}
where, $ \eta = 2 \delta -\delta^2$ and $h_{\mathcal{R}}(\delta) = 1- \eta \mu_1 < 1$. The constant $\mu_1$ is dependent on the choice of sampling strategy $\mathcal{R}$. These sets were first introduced in \cite{morshed:momentum} to analyze the MSKM algorithm. We note that these sets are of crucial importance for proving convergence results related to projection-based momentum methods for the linear feasibility problem. 

\begin{mdframed}[backgroundcolor=gray!20,   topline=false,   bottomline =false,   rightline=false,   leftline=false]  \begin{theorem}
\label{th:b1}
Let, $x_k$ is the random iterate generated by the basic method with $0 < \delta < 2$. Then, the following identities
\begin{align*}
 \E[ d_B(x_{k+1},\mathcal{X})^2]  \leq  [h_{\mathcal{R}}(\delta)]^{k+1} d_B(x_0,\mathcal{X})^2  \ \text{and} \   \E[f(x_{k+1})] \leq \frac{\mu_2}{2} [h_{\mathcal{R}}(\delta)]^{k+1} d_B(x_0,\mathcal{X})^2
\end{align*}
hold. Also the average iterate $\Tilde{x}_k = \sum \limits_{l=0}^{k-1} x_l$ for all $k \geq 1$ satisfies 
\begin{align*}
    \E[d_B(\Tilde{x}_k,\mathcal{X})^2]  \leq \frac{d_B(x_0,\mathcal{X})^2 }{2\delta k(2-\delta) \mu_1} \quad \text{and} \quad \E[f(\Tilde{x}_k)] \leq \frac{d_B(x_0,\mathcal{X})^2}{2\delta k(2-\delta)},
\end{align*}
where, the constant $\mu_1$ is depended on the choice of sampling rule $\mathcal{R}$.
\end{theorem} \end{mdframed}

\begin{proof}
See Appendix \ref{pb1}.
\end{proof}

\begin{mdframed}[backgroundcolor=gray!20,   topline=false,   bottomline =false,   rightline=false,   leftline=false]  
\begin{remark}
\label{rem:specialresult}
Note that, Theorem \ref{th:b1} is a generalized result. As the constant $\mu_1$ varies form rules to rules, for different choices of sampling rules, we get the corresponding convergence results. In subsection \ref{sec:special}, we discuss two special algorithms, and their respective convergence results obtained from Theorem \ref{th:b1}.
\end{remark} \end{mdframed}

\allowdisplaybreaks{\begin{mdframed}[backgroundcolor=gray!20,   topline=false,   bottomline =false,   rightline=false,   leftline=false]  \begin{theorem}
\label{th:1}
Let $\{x_k\}$ be the sequence of random iterates generated by algorithm \ref{alg:marsp} and let $0 \leq \gamma < 1$ such that $(\delta, \gamma) \in Q$. Let's denote, $ \Pi_1 =  \sqrt{h_{\mathcal{R}}(\delta)}, \ \Pi_2 = \Pi_4 = \gamma, \ \Pi_3 =\delta \sqrt{\mu_2 }$ and $\Gamma_1, \Gamma_2, \Gamma_3, \rho_1, \rho_2$ as in \eqref{t0} with the above parameter choice. Then the sequence of iterates $\{x_k\}$ converges and the following result holds:
\allowdisplaybreaks{\begin{align*}
  \E \begin{bmatrix}
d_B(x_{k+1}, \mathcal{X})  \\[6pt]
\|x_{k+1}-x_k\|_B 
\end{bmatrix} & \leq \begin{bmatrix}
-\Gamma_2 \Gamma_3 \ \rho_1^{k}+ \Gamma_1 \Gamma_3 \ \rho_2^{k} \\[6pt]
- \Gamma_3 \ \rho_1^{k}+ \Gamma_3 \ \rho_2^{k} 
\end{bmatrix} \ d_B(x_0,\mathcal{X}) \leq \begin{bmatrix}
 1 \\[6pt]
2 \Gamma_3
\end{bmatrix} \  \rho_2^{k}  \ d_B(x_0,\mathcal{X}),
\end{align*}}
where $ \Gamma_3 \geq 0$ and $ 0 \leq |\rho_1| \leq  \rho_2 < 1$.
\end{theorem} \end{mdframed}}

\begin{proof}
See Appendix \ref{pt1}.
\end{proof}

\begin{mdframed}[backgroundcolor=gray!20,   topline=false,   bottomline =false,   rightline=false,   leftline=false]  \begin{remark}
\label{rem:mom1bound}
Note that, we can simplify Theorem \ref{th:1} to develop working bounds for the momentum parameter $\gamma$. Indeed, if the parameter pair $(\delta, \gamma)$ belongs to $Q$, i.e., $0 \leq \gamma <  \frac{1-\sqrt{h_{\mathcal{R}}(\delta)}}{1-\sqrt{h_{\mathcal{R}}(\delta)} + \delta \sqrt{\mu_2}}$ holds then the proposed momentum algorithm converges. Define, $\Tilde{\mu}_1 = \frac{\mu_1}{\mu_1+\sqrt{\mu_2}}$ and $\Tilde{\mu}_2 =\frac{1-\sqrt{1-\mu_1}}{1-\sqrt{1-\mu_1}+\sqrt{\mu_2}}$. Now, the function $\mathcal{H}(\delta) = \frac{1-\sqrt{h_{\mathcal{R}}(\delta)}}{1-\sqrt{h_{\mathcal{R}}(\delta)} + \delta \sqrt{\mu_2}}$ is decreasing in the interval $(0,2)$. That implies the function $\mathcal{H}$ attains maximum value at $\delta \rightarrow 0$, i.e,  
\begin{align}
\label{param1}
     \max_{\delta \in (0,2)} \mathcal{H}(\delta) = \lim_{\delta \rightarrow 0} \frac{1-\sqrt{h_{\mathcal{R}}(\delta)}}{1-\sqrt{h_{\mathcal{R}}(\delta)} + \delta \sqrt{\mu_2}} = \frac{\mu_1}{\mu_1+\sqrt{\mu_2}} \leq 0.5.
\end{align}
That implies the allowable range of $\gamma$ values for which the convergence result of Theorem \ref{th:1} holds is $0 \leq \gamma < 0.5$. More specifically, one can easily show that the following piece-wise conditions are valid:
\begin{align}
\label{param2}
 0 < \delta < 1 :\rightarrow  \gamma <\Tilde{\mu_1}  - (\Tilde{\mu_1} -\Tilde{\mu_2} )\delta,  \quad  1 < \delta < 2 :\rightarrow   \gamma < 2\Tilde{\mu_2} -\Tilde{\mu_2} \delta.
\end{align}
Furthermore, if $(\gamma, \delta) \in \{0 < \delta < 2, \ 0 < \gamma < 0.5, \  \gamma \ \leq \ 0.5 \Tilde{\mu_1} (2-\delta) \}$, then they must reside inside $Q$, i.e., $(\gamma, \delta) \in Q$.
\end{remark} \end{mdframed}

Note that, the following relation holds
\begin{align}
\big | \E [d_B(x_{k}, \mathcal{X})] \big |^2  \leq \E \left[d_B(x_k,\mathcal{X})^2\right],
\end{align}
for any random vector $x \in \R^n$. Therefore, the convergence result of Theorem \ref{th:1} is weaker compared to the usual $L_2$ convergence (the decay of the term $\E \left[d_B(x_k,\mathcal{X})^2\right]$). In Theorem \ref{th:mom2}, we derive the necessary decay bounds to show the convergence of the term $\E \left[d_B(x_k,\mathcal{X})^2\right]$.

\begin{mdframed}[backgroundcolor=gray!20,   topline=false,   bottomline =false,   rightline=false,   leftline=false]  \begin{lemma}
\label{lem:skm2}
The sequence $\{x_k\}$ generated by the ASPM algorithm are point-wise closer to the feasible region $\mathcal{X}$ with respect to the $B-$norm, i.e., for all $\bar{x} \in \mathcal{X}$ and $k \geq 0$, we have
\begin{align*}
    \|x_{k+1}-\bar{x}\|_B \ \leq \  \|x_{k}-\bar{x}\|_B.
\end{align*}
\end{lemma} \end{mdframed}

\begin{proof}
See Appendix \ref{pskm2}.
\end{proof}

\begin{mdframed}[backgroundcolor=gray!20,   topline=false,   bottomline =false,   rightline=false,   leftline=false]  \begin{lemma}
\label{lem:mom1}
Let $x_{k+1}$ is generated by the momentum algorithm, then we have
\begin{align*}
   \E[d_B(x_{k+1} & ,\mathcal{X})^2  \ | \ i \sim \mathcal{R}] + \zeta \E[\|x_{k+1}-x_{k}\|_B^2\ | \ i \sim \mathcal{R}]  \nonumber \\
   & \leq (1+\gamma) \ d_B(x_{k},\mathcal{X})^2 - \gamma \ d_B(x_{k-1},\mathcal{X})^2  + (\gamma^2+\gamma+ \zeta \gamma^2) \ \|x_{k}-x_{k-1}\|_B^2 \nonumber \\
    & + 2 \gamma \delta (1+\zeta) f(x_{k-1}) - 2 \delta [2-(\delta -\gamma)(1+\zeta)] f(x_{k}),
\end{align*}
where, $\zeta \geq 0$.
\end{lemma} \end{mdframed}

\begin{proof}
See Appendix \ref{pmom1}.
\end{proof}

\allowdisplaybreaks{\begin{mdframed}[backgroundcolor=gray!20,   topline=false,   bottomline =false,   rightline=false,   leftline=false]  \begin{theorem}
\label{th:mom2}
Let $\{x_k\}$ be the sequence of random iterates generated by algorithm \ref{alg:marsp}. Let $0 \leq \gamma <1$ and $\zeta \geq 0$ such that $(\delta, \gamma, \zeta) \in R \cap S$. Then the sequence of iterates $\{x_k\}$ converges and the following result holds.
\begin{align*}
  \E [d_B(x_{k+1},\mathcal{X})^2]  \leq  \rho^k (1+\alpha)   d_B(x_0,\mathcal{X})^2  \ \text{and} \  \E [f(x_{k+1})]  \leq \frac{\mu_2(1+\alpha)}{2} \rho^k  d_B(x_0,\mathcal{X})^2.
\end{align*}
Also the average iterate $\Tilde{x}_k = \sum \limits_{l=1}^{k} x_l$ for all $k \geq 0$ satisfies
\begin{align*}
    \E[d_B(\Tilde{x}_k,\mathcal{X})^2]  \leq \frac{(1+\alpha) \ d_B(x_0,\mathcal{X})^2}{ k(1-\rho)} \quad \text{and} \quad \E[f(\Tilde{x}_k)] \leq \frac{ \mu_2 (1+\alpha) }{2 k(1-\rho)} \ d_B(x_0,\mathcal{X})^2,
\end{align*}
where, $\alpha \geq 0$,  $0 < \rho < 1$.
\end{theorem} \end{mdframed}}

\begin{proof}
See Appendix \ref{pmom2}.
\end{proof}

\subsection{Certificate of Feasibility}
\label{sec:certificate}
In this subsection, we propose a generic Theorem related to the feasibility certification for the halting of the proposed momentum algorithm after finitely many iterations. It is a generalization of the feasibility certification results obtained in \cite{haddock:2017}, \cite{morshed:momentum}, and to a certain extent, it can be said to be an extension of Telgen's result \cite{telgen:1982}. Before we delve into the main Theorem, first let us define the following quantities:
\begin{align}
    \psi = \max_{j \in \{1,2,...,m\}} \|a_j\|_2, \ \lambda_2 = \lambda_{\max}(B), \ \lambda_1 = \lambda_{\min}(B), \ \xi = \frac{\lambda_2}{\lambda_1}.
\end{align}

\begin{mdframed}[backgroundcolor=gray!20,   topline=false,   bottomline =false,   rightline=false,   leftline=false]  \begin{theorem}
\label{th:4}
Assume $A, b$ are rational matrices with binary encoding length, $\sigma$. Also assume parameters $0 < \delta < 2$ and $\gamma, \zeta \geq 0$ satisfy the condition $(\delta, \gamma, \zeta) \in Q \cup \left(R \cap S\right)$. Suppose we run the momentum algorithm on the system $ Ax \leq b$ with $x_0=0$, and the number of iterations $k$ satisfies the following lower bound:
\begin{align*}
   \frac{4 \sigma - 4 -\log n + \log (1+\alpha)+ \log \xi + 2 \log \psi}{\log \left(\frac{1}{\bar{\rho}}\right)} < k-1.
\end{align*}
Define, $\bar{\rho} = \max\{\rho^2_2, \rho\} < 1$, where $\rho_2$ and $\rho$ are defined in Theorem \ref{th:1} and Theorem \ref{th:mom2} for the choice $(\delta,\gamma) \in Q$ and $(\delta, \gamma, \zeta)  \in R \cap S$, respectively. If the system $Ax \leq b$ is feasible, then,  
\begin{align*}
    p \ \leq H(\sigma, \xi, \psi, \alpha, k, \bar{\rho}) =   \sqrt{\frac{\xi(1+\alpha)}{n}}  \psi
  \ 2^{2\sigma -2} \ \bar{\rho}^{\frac{k-1}{2}},
\end{align*}
where $p$ is the probability that the current iterate $x_{k}$ is not a certificate of feasibility. Note that, with respect to $k$, function $H(\sigma, \xi, \psi, \alpha, k, \bar{\rho})$ is a decreasing function.
\end{theorem} \end{mdframed}

\begin{proof}
See Appendix \ref{p4}.
\end{proof}

\begin{mdframed}[backgroundcolor=gray!20,   topline=false,   bottomline =false,   rightline=false,   leftline=false]  \begin{corollary}
\label{cor:2}
(Theorem 4.7, Remark 4.8 in \cite{morshed:momentum}) Suppose $A, b$ are rational matrices with binary encoding length, $\sigma$, and that we run the momentum induced SKM method \footnote{Note that, by choosing $\gamma = 0$ in Corollary \ref{cor:2}, one can recover the certificate of feasibility for the SKM method proved in \cite{haddock:2017}.} on the system (MSKM method in \cite{morshed:momentum} with $0 < \delta < 2$ and $\gamma, \zeta \geq 0$ such that $(\delta, \gamma, \zeta) \in Q \cup \left(R \cap S\right)$). Suppose the number of iterations $k$ satisfies the following lower bound:
\begin{align*}
   \frac{4 \sigma - 4 -\log n + \log (1+\alpha) + 2 \log \psi}{\log \left(\frac{1}{\bar{\rho}}\right)} < k-1,
\end{align*}
If the system $Ax \leq b$ is feasible, then,  
\begin{align*}
    p \ \leq \  \sqrt{\frac{1+\alpha}{n}} \ 2^{2\sigma -2} \ \psi \ \bar{\rho}^{\frac{k-1}{2}},
\end{align*}
where $p = $the probability that the current update $x_k$ is not a certificate of feasibility.
\end{corollary} \end{mdframed}

\begin{proof}
Take $B = I, \ S_i = e_i$  along with the greedy sampling rule of subsection \ref{ss} in Theorem \ref{th:4}. Then we get, $\xi = 1$. Now, considering Theorem \ref{th:4}, we get the result of Corollary \ref{cor:2}.  
\end{proof}

\subsection{Cesaro Average}
\label{sec:cesaro}
Next, we discuss the convergence of the function decay, i.e., $\E[f(x)]$ for the average iterate $\Tilde{x_k} = \frac{1}{k} \sum \limits_{l =1}^{k}x_l$ \footnote{This is widely known as convergence with respect to the Cesaro average.}. We derive a $\mathcal{O}(\frac{1}{k})$ convergence for the proposed momentum algorithm with respect to the Cesaro average that is better than the rate obtained in Theorem \ref{th:mom2}. Furthermore, the convergence result holds under a somewhat weaker condition on the parameter pair $(\delta, \gamma)$. Moreover, several well-known results can be obtained as special cases from our proposed Theorem.

\begin{mdframed}[backgroundcolor=gray!20,   topline=false,   bottomline =false,   rightline=false,   leftline=false]  \begin{theorem}
\label{th:cesaro}
Let $\{x_k\}$ be the random sequence generated by the momentum algorithm. Let, $0 < \delta <2$ and $0 \leq \gamma < 1 $ such that the condition $0 < \delta < 2(1-\gamma)$ holds. Define $\Tilde{x_k} = \frac{1}{k} \sum \limits_{l =1}^{k}x_l$, then
\begin{align*}
    \E \left[f(\bar{x}_k)\right] \leq \frac{ (1-\gamma)^2 \ d_B(x_0,\mathcal{X})^2+ 2 \delta \gamma f(x_0)}{2 \delta k \left(2- 2 \gamma -\delta\right)}.
\end{align*}
\end{theorem} \end{mdframed}

\begin{proof}
See Appendix \ref{pcesaro}.
\end{proof}

\begin{mdframed}[backgroundcolor=gray!20,   topline=false,   bottomline =false,   rightline=false,   leftline=false]  \begin{remark}
\label{rem:ces2}
Theorem \ref{th:cesaro} holds for a wide range of projection and momentum parameter pairs (i.e., $(\delta, \gamma)$) compared to Theorem \ref{th:mom2}. Moreover, the convergence rate obtained in Theorem \ref{th:cesaro} is substantially better than the one obtained in Theorem \ref{th:mom2}. Finally, one can obtain several well-known results as special cases by choosing different parameter matrix $B$ and different sampling rules $\mathcal{R}$. For instance, the following result can be obtained for the MSKM method.
\end{remark} \end{mdframed}

\begin{mdframed}[backgroundcolor=gray!20,   topline=false,   bottomline =false,   rightline=false,   leftline=false]  \begin{corollary} 
\label{cor:cesaro}
(Theorem 4.9 in \cite{morshed:momentum}) Let $\{x_k\}$ be the random sequence generated by the MSKM algorithm \footnote{Note that, by choosing $\gamma = 0$ in Corollary \ref{cor:cesaro}, one can recover the Cesaro average Theorem of the SKM method proposed in \cite{morshed:momentum}.} Take, $0 \leq \gamma < 1 $ and $0 < \delta < 2(1-\gamma)$. Define $\Tilde{x_k} = \frac{1}{k} \sum \limits_{l =1}^{k}x_l$ and $f(x)$ as in \eqref{eq:sk3}, then
\begin{align*}
    \E \left[f(\bar{x}_k)\right] \leq \frac{ (1-\gamma)^2 \ d(x_0,\mathcal{X})^2+ 2 \delta \gamma f(x_0)}{2 \delta k \left(2- 2 \gamma -\delta\right)}.
\end{align*}
\end{corollary} \end{mdframed}
\begin{proof}
Consider the special sketched loss sampling in subsection \ref{ss} along with $B =I$.  $S_i = e_i$. Then it can be easy to check that Theorem \ref{th:cesaro} resolves into Corollary \ref{cor:cesaro}.
\end{proof}

\subsection{Special Cases}
\label{sec:special}
In this subsection, we discuss some special algorithms and the corresponding convergence results that can recover from our proposed methods. For simplification, we consider the following simple cases: $S_i = e_i$ and $B = I, A$. For completeness, in Appendix \ref{extracase}, we provide another variant of algorithms and the corresponding convergence results. 

\paragraph{Momentum Sampling Kaczmarz Motzkin} Consider, $q = m, \ B =I, \ S= e_i$ in the ASPM method. Assume, $\|a_i\|^2 = 1, i = 1,2,...,m$. Then the update formula of the ASPM method with the greedy sketched loss sampling (i.e., $i \sim \mathcal{G}(\tau)$ subsection \ref{ss}) resolves into the following update formula:
\begin{align}
\label{skm}
  x_{k+1} = x_k- \delta \frac{\left(a_{i}^Tx_k -b_{i}\right)^+}{\|a_{i}\|_2^2} a_{i} + \gamma (x_k-x_{k-1}).   
\end{align}
where, $i = \argmax_{i \in \phi_k(\tau)} \{a_i^Tx_k-b_i, 0\}$ and $\phi_k(\tau)$ denotes the collection of $\tau$ rows chosen uniformly at random out of $m$ rows of the constraint matrix $A$. Using the above parameter choice we get, $R= I_{m \times m}$. Since, $\omega_1 = \min_{i \in \{1,...,m\}} \|A^Te_i\|^2 = \min_{i \in \tau} \|a_i\|^2 = 1 = \omega_2 $  and $\sigma = L^2$ (Lemma \ref{lem0}).  Now, we can calculate the constants $\mu_1$, $\mu_2$ as follows:
\begin{align}
\label{mu:skm}
  \mu_1 =  \frac{1}{L^2}\min \left\{\frac{1}{m-\tau+1}, \frac{1}{m-s}\right\} \geq \frac{1}{m L^2}, \ \   \mu_2 = \min \left\{1, \frac{\tau}{m} \lambda_{\max}(A^TA)  \right\}.
\end{align}
Here, $s$ is the number of zero entries in the residual $(Ax-b)^+$. Using these parameter values, in the following Corollaries, we derived the convergence results for both the SKM ($\gamma = 0$ in \eqref{skm}) and MSKM methods. Note that, with the choice $\tau =1$ and $\tau = m$ in SKM, we can recover the RK method \cite{lewis:2010} and MR method \cite{motzkin} respectively.

\begin{mdframed}[backgroundcolor=gray!20,   topline=false,   bottomline =false,   rightline=false,   leftline=false]  \begin{corollary}
\label{cor:b1}
(Theorem 1 in \cite{haddock:2017}, Lemma 9 in \cite{morshed2020generalization}, Theorem 1 in \cite{morshed:momentum}) Let, $x_k$ be the random iterate generated by the SKM method with $0 < \delta < 2$, $\eta = 2\delta- \delta^2$. Then, the following identities
\begin{align*}
 \E[ d(x_{k},\mathcal{X})^2]  \leq  \left(1-\frac{\eta}{mL^2}\right)^{k} d(x_0,\mathcal{X})^2,
\end{align*}
and 
\begin{align*}
   \E\left[\big | (a_{i}^Tx_k-b_{i})^+ \big |^2\right] \leq  \min \left\{1, \frac{\tau}{m} \lambda_{\max}(A^TA) \right\} \left(1-\frac{\eta}{mL^2}\right)^{k} d(x_0,\mathcal{X})^2,
\end{align*}
hold. Also the average iterate $\Tilde{x}_k = \sum \nolimits_{l=0}^{k-1} x_l$ for all $k \geq 1$ satisfies
\begin{align*}
    \E[d(\Tilde{x}_k,\mathcal{X})^2]  \leq \frac{ m L^2  }{2\delta k(2-\delta)}  d(x_0,\mathcal{X})^2 \quad \text{and} \quad  \E\left[\big | (a_{i}^T\Tilde{x}_k-b_{i})^+ \big |^2\right] \leq \frac{d(x_0,\mathcal{X})^2}{2\delta k(2-\delta)}.
\end{align*}
\end{corollary} \end{mdframed}

\begin{proof}
Consider, $\mu_1$ and $\mu_2$ values from \eqref{mu:skm} in Theorem \ref{th:b1}. Then, with simplification, we get the result of Corollary \ref{cor:b1}.
\end{proof}

\begin{mdframed}[backgroundcolor=gray!20,   topline=false,   bottomline =false,   rightline=false,   leftline=false]  \begin{corollary} 
\label{cor:mskm2}
(Theorem 1 in \cite{morshed:momentum}) Let $\{x_k\}$ be the sequence of random iterates generated by the MSKM algorithm starting with $x_0 = x_1 \in \R^n$. With $0 < \delta < 2$, the sequence of iterates $\{x_k\}$ converges and the following result holds:
{\allowdisplaybreaks
\begin{align*}
\E \begin{bmatrix}
d(x_{k+1}, \mathcal{X})  \\[6pt]
\|x_{k+1}-x_k\|
\end{bmatrix} & \leq \begin{bmatrix}
-\Gamma_2 \Gamma_3 \ \rho_1^{k}+ \Gamma_1 \Gamma_3 \ \rho_2^{k} \\[6pt]
- \Gamma_3 \ \rho_1^{k}+ \Gamma_3 \ \rho_2^{k} 
\end{bmatrix} \ d(x_{0}, \mathcal{X})  \leq \begin{bmatrix}
 \rho_2^{k} \\[6pt]
2 \Gamma_3 \ \rho_2^{k} 
\end{bmatrix} \ d(x_{0}, \mathcal{X}).
\end{align*}}
\end{corollary} \end{mdframed}
\begin{proof}
Consider, $\mu_1$ and $\mu_2$ values from \eqref{mu:skm} in Theorem \ref{th:1}. Then, with simplification, we get the result of Corollary \ref{cor:mskm2}. 
\end{proof}

\begin{mdframed}[backgroundcolor=gray!20,   topline=false,   bottomline =false,   rightline=false,   leftline=false]  \begin{corollary} 
\label{cor:skm5}
(Theorem 4.6 in \cite{morshed:momentum}) Let $\{x_k\}$ be the sequence of random iterates generated by the MSKM algorithm with $x_0 = x_1 \in \R^n$. With $0 < \delta < 2$, the sequence of iterates $\{x_k\}$ converges and the following results hold:
\begin{align*}
  \E [d(x_{k+1},\mathcal{X})^2]  \leq  \rho^k (1+\alpha) \  d(x_0,\mathcal{X})^2 \ \ \text{and} \ \ \E [f(x_{k+1})]  \leq \frac{\mu_2(1+\alpha)}{2} \rho^k \ d(x_0,\mathcal{X})^2,
\end{align*}
where, $\alpha \geq 0$,  $0 < \rho < 1$ are provided in \eqref{mom28} \eqref{mom29} respectively.
\end{corollary} \end{mdframed}
\begin{proof}
Consider, $\mu_1$ and $\mu_2$ values from \eqref{mu:skm} in Theorem \ref{th:mom2}. Then, with simplification, we get the result of Corollary \ref{cor:skm5}.
\end{proof}

\paragraph{Momentum Sampling Co-ordinate Descent (MSCD)} Take, $q = m =n, \ S= e_i, \ B = A$ ($A$ is positive definite). Assume, $A_{ii} = 1, i = 1,2,...,m$. Then the ASPM method with the greedy sketched loss sampling (i.e., $i \sim \mathcal{G}(\tau)$ subsection \ref{ss}) resolves into the following update formula:
\begin{align}
\label{rcd-pd}
  x_{k+1} = x_k- \delta \frac{\left(a_{i}^Tx_k -b_{i}\right)^+}{A_{ii}} e_{i} + \gamma(x_k-x_{k-1}).   
\end{align}
where, $A_{ii}$ is the $i^{th}$ diagonal entry of matrix $A \in \R^{n \times n}$ and $i = \argmax_{i \in \phi_k(\tau)} \{a_i^Tx_k-b_i, 0\}$ and $\phi_k(\tau)$ denotes the collection of $\tau$ rows chosen uniformly at random out of $m$ rows of the constraint matrix $A$. Using the above parameter choice in subsection \ref{ss}, we get $R= I_{m \times m}$, $\|B^{-\frac{1}{2}}A^TR^T\|^2_F = \Tr(A) = m$, $\lambda_{\max}(B^{-\frac{1}{2}}A^TR^TRA B^{-\frac{1}{2}}) = \lambda_{\max}(A)$, $\|A^TS_i\|^2_{A} = A_{ii} = 1$ and $\sigma = L^2$ (Lemma \ref{lem0}). Now, we can calculate the constants $\mu_1$, $\mu_2$ as follows:
\begin{align}
\label{mu:scd}
  \mu_1 =  \frac{1}{L^2}\min \left\{\frac{1}{m-\tau+1}, \frac{1}{m-s}\right\} \geq \frac{1}{m L^2}, \ \   \mu_2 = \min \left\{1, \frac{\tau}{m} \lambda_{\max}(A) \right\}.
\end{align}
Here, $s$ is the number of zero entries in the residual $(Ax-b)^+$. Using these parameter values, in the following Corollaries, we derive the convergence result for \textit{Sampling Co-ordinate Descent} (SCD) and MSCD methods. Note that, with the choice $\tau =1$ in SCD, we can recover the \textit{Randomized Co-ordinate Descent} (RCD) method proposed in \cite{lewis:2010} for solving a system of linear equations. Then, we derive the convergence results for the MSCD method.

\begin{mdframed}[backgroundcolor=gray!20,   topline=false,   bottomline =false,   rightline=false,   leftline=false]  \begin{corollary}
\label{cor:rcd-pd}
(New Theorem) Let, $x_k$ be the random iterate generated by the SCD method with $0 < \delta < 2$. Then, the following identities
\begin{align*}
 \E[ d_A(x_{k},\mathcal{X})^2]  \leq  \left(1-\frac{\eta}{m L^2}\right)^{k}  d_A(x_0,\mathcal{X})^2,
\end{align*}
and 
\begin{align*}
   \E\left[\big | (a_{i}^Tx_k-b_{i})^+ \big |^2\right] \leq  \frac{\tau \lambda_{\max}(A)}{m} \left(1-\frac{\eta}{m L^2}\right)^{k} d(x_0,\mathcal{X})^2,
\end{align*}
hold. Also the average iterate $\Tilde{x}_k = \sum \nolimits_{l=0}^{k-1} x_l$ for all $k \geq 1$ satisfies
\begin{align*}
    \E[d_A(\Tilde{x}_k,\mathcal{X})^2]  \leq \frac{m L^2  }{2\delta k(2-\delta)} d_A(x_0,\mathcal{X})^2 \quad \text{and} \quad  \E\left[\big | (a_{i}^T\Tilde{x}_k-b_{i})^+ \big |^2\right]  \leq \frac{d_A(x_0,\mathcal{X})^2}{2\delta k(2-\delta)}.
\end{align*}
\end{corollary} \end{mdframed}

\begin{proof}
Consider, $\mu_1$ and $\mu_2$ values from \eqref{mu:scd} in Theorem \ref{th:b1}. Then, with simplification, we get the result of Corollary \ref{cor:rcd-pd}.
\end{proof}

\begin{mdframed}[backgroundcolor=gray!20,   topline=false,   bottomline =false,   rightline=false,   leftline=false]  \begin{corollary} 
\label{cor:mscd1}
(New Theorem) Let $\{x_k\}$ be the sequence of random iterates generated by the MSCD algorithm starting with $x_0 = x_1 \in \R^n$. With $0 < \delta < 2$, the sequence of iterates $\{x_k\}$ converges and the following result holds:
{\allowdisplaybreaks
\begin{align*}
\E \begin{bmatrix}
d_A(x_{k+1}, \mathcal{X})  \\[6pt]
\|x_{k+1}-x_k\|_A
\end{bmatrix} & \leq \begin{bmatrix}
-\Gamma_2 \Gamma_3 \ \rho_1^{k}+ \Gamma_1 \Gamma_3 \ \rho_2^{k} \\[6pt]
- \Gamma_3 \ \rho_1^{k}+ \Gamma_3 \ \rho_2^{k} 
\end{bmatrix} \ d_A(x_{0}, \mathcal{X})  \leq \begin{bmatrix}
 \rho_2^{k} \\[6pt]
2 \Gamma_3 \ \rho_2^{k} 
\end{bmatrix} \ d_A(x_{0}, \mathcal{X}).
\end{align*}}
\end{corollary} \end{mdframed}

\begin{proof}
Consider, $\mu_1$ and $\mu_2$ values from \eqref{mu:scd} in Theorem \ref{th:1}. Then, with simplification, we get the result of Corollary \ref{cor:mscd1}.
\end{proof}

\begin{mdframed}[backgroundcolor=gray!20,   topline=false,   bottomline =false,   rightline=false,   leftline=false]  \begin{corollary} 
\label{cor:mscd2}
(New Theorem) Let $\{x_k\}$ be the sequence of random iterates generated by the MSCD algorithm with $x_0 = x_1 \in \R^n$. With $0 < \delta < 2$, the sequence of iterates $\{x_k\}$ converges and the following results hold:
\begin{align*}
  \E [d_A(x_{k+1},\mathcal{X})^2]  \leq  \rho^k (1+\alpha) \  d_A(x_0,\mathcal{X})^2 \ \text{and}  \ \E [f(x_{k+1})]  \leq \frac{\mu_2(1+\alpha)}{2} \rho^k \ d_A(x_0,\mathcal{X})^2,
\end{align*}
where, $\alpha \geq 0$,  $0 < \rho < 1$ are provided in \eqref{mom28} \eqref{mom29} respectively.
\end{corollary} \end{mdframed}
\begin{proof}
Consider, $\mu_1$ and $\mu_2$ values from \eqref{mu:scd} in Theorem \ref{th:mom2}. Then, with simplification, we get the result of Corollary \ref{cor:mscd2}.
\end{proof}

\begin{mdframed}[backgroundcolor=gray!20,   topline=false,   bottomline =false,   rightline=false,   leftline=false]  
\begin{remark}
\label{special remark}
Note that, in the MSKM algorithm we assumed matrix $A$ has normalized rows. Similarly, for the MSCD method, we assumed $A_{ii} = 1$. We assumed this to show the equivalency with existing algorithms. It can be noted that these assumptions are not required for computational performance. Indeed, we find that irrespective of the assumptions the proposed algorithms perform the same. 
\end{remark}
 \end{mdframed}

\section{Numerical Experiments}
\label{sec:num}

In this section, we study the computational performance of the proposed momentum methods. We implement the above-mentioned methods in \textit{MATLAB R2020a} and carry out the experiments in a workstation with 64GB RAM, Intel(R) Xeon(R) CPU E5-2670, two processors running at 2.30 GHz. Throughout the experiment, we fixed the projection parameter to $\delta = 1$ \footnote{Empirically this specific choice has the best computational performance for both linear systems \cite{richtrik2017stochastic,loizou:2017} and linear feasibility problems \cite{haddock:2017, morshed2020generalization,morshed:momentum}.}. We test the proposed algorithms on two types of datasets: 1) synthetic data (Gaussian system), and 2) real-world data (Netlib LP instances). First, we select two types of projection algorithms and their respective momentum variants, i.e., GK and Adaptive Co-ordinate Descent for the computational experiments. Second, we fix the momentum parameter as $\gamma = 0.1, 0.2, 0.3, 0.4, 0.5$. Note that, we don't need to calculate the constants $\mu_1$ and $\mu_2$ for the selection of momentum parameter $\gamma$. Until otherwise mentioned, throughout the computational section, we fixed the following parameters: 1) initial point $x_0$ is fixed as $1000*[1,1,...,1]^T$ which is very far away from the feasible region of the considered test instances, 2) stopping criteria is selected as either $\|\left(Ax-b\right)^+\|_2 \leq 10^{-05}$ or number of iterations $300,000$, and 3) all experiments were run for $10$ times and the averaged performance was reported. This section is divided into four subsections. In the first, we discuss the test instances and the convergence measures we are going to test. Furthermore, at the end of the first subsection, we summarise the findings in brief from our numerical experiments. In the second subsection, we discuss the performance comparison of the proposed adaptive methods with no momentum under different sampling rules. In the third and fourth subsection, we discuss the effect of momentum on both GK and GCD methods.

\subsection{Test Instances and convergence parameters}
For the GK  method, the synthetic data is generated as follows: matrix $A\in \R^{m \times n}$ and vector $x \in \R^{n}$ are chosen to be i.i.d $\mathcal{N}(0,1)$, then the right-hand side $b$ of the feasibility problem is set to $b = Ax+ |\epsilon|$, where $\epsilon$ is a Gaussian random vector. It can be noted that, with this specific selection process, we maintain the consistency of the linear feasibility problem. For the GCD method, the synthetic data is generated as follows: matrix $A$ is set as $A = G^TG$, where $G \in \R^{m \times n}$ is a Gaussian matrix, and the right-hand side $b$ is generated by the same procedure as in the GK  method. From our convergence results, we have that the quantities $\E[d_B(x_k,\mathcal{X})^2] = \E[\|x_k - \mathcal{P}_{\mathcal{X}}^{B}(x_k) \|_B^2]$ and $E[f(x_k)]$ converges linearly to zero. For this reason, we select the decay of these quantities with respect to CPU time and number of iterations \footnote{Note that, at iteration $k$ we don't calculate the value of $\mathcal{P}_{\mathcal{X}}^{B}(x_k)$, instead we use the $x$ value used initially to generate the matrix $A$. This is not an obvious measure but for illustration purpose, we show the decay of the quantity $\|x_k-x_{int}\|^2_B$.}. Moreover, for simplification of implementation, for the GK  method we assumed $\|a_i\|^2 =1$ for all $i$, and for the GCD method we assumed $a_{ii} =1$ for all $i$. With this simplification, we get $\E[f(x_k)] = \E\left[|\left(a_i^Tx_k-b_i\right)^+|^2 | \ i \sim \mathcal{R}\right] \propto \|\left(Ax_k-b\right)^+\|_2$. To investigate the solution quality of each algorithm as they progress, we measure the number of violated constraints after each iteration. Let's define the following quantity:
\begin{align*}
 \text{Fraction of Satisfied Constraints (FSC)} = \frac{ \# \ \text{of satisfied constraints at iteration} \ k}{\text{Total number of constraints ($m$)}}.   
\end{align*}
Note that, $0 \leq \text{FSC} \leq 1$ holds for each $k$. Based on the above discussion, we get the following convergence measures: 1) Positive residual error ($\|\left(Ax_k-b\right)^+\|_2$), 2) Relative error ($\|x_k-x_{int}\|_B/\|x_0-x_{int}\|_B$), and 3) FSC. In the horizontal axis, we use either the number of iterations or the CPU time measured using the MATLAB tic-toc function. Note that, with the choice $\gamma = 0$, the momentum variants resolves into the basic methods. In the following, we summarise the findings of our numerical experiments:

\begin{itemize}
    \item From our experiments, we find that the proposed greedy sampling rules heavily outperform the traditional sampling rules in terms of CPU time and solution quality. For the greedy sketched loss sampling, the choice $1 < \tau \ll m$ leads to the best-performing methods whereas the choice $\theta \in [0.5, 1)$ leads to better algorithms for the greedy capped sampling rules.
    
    \item Throughout the experiment, we choose momentum parameter $\gamma$ arbitrarily. We find that for the majority of the test instances the choice $\gamma \in [0.3,0.4]$ leads to the best momentum variants. The choice $\gamma = 0.5$ generates good results for a handful of test instances but in generals this choice fails to converge faster in most cases and has worse performance than the basic method with no momentum.
    
    \item We test our methods for a wide range of random datasets with varying condition number of matrix $A$. For the case of ill-conditioned feasibility problems (condition number $A$ is large), the momentum variants outperform the basic method heavily. For the case of small condition numbers, the performance improvement is marginal. 
    
    \item For both GK and GCD methods, with the choice $\gamma \in (0,0.5)$, the momentum variants always converge faster and generate better solutions as they progress. The choice of momentum parameter is small compared to the choice $\gamma = 0.9$ which is used for the SGD method for deep neural network training.
\end{itemize}

\subsection{Comparison among sampling rules without momentum}
In this subsection, we carry out comparison experiments for both GK and GCD methods with respect to the proposed sampling rules. For a fair understanding we choose the following six sampling rules: Uniform ($\tau =1$), $\tau = 5$, $\tau = 50$, $\tau = 100$, Maximum distance ($\tau = m$), Capped ($\tau_1 = 1, \tau_2 = m, \theta = 0.5$). We plot positive residual error, relative error, FSC vs time and number of iterations graphs.

\begin{figure}[htbp]
\centering
    \includegraphics[scale = 0.79]{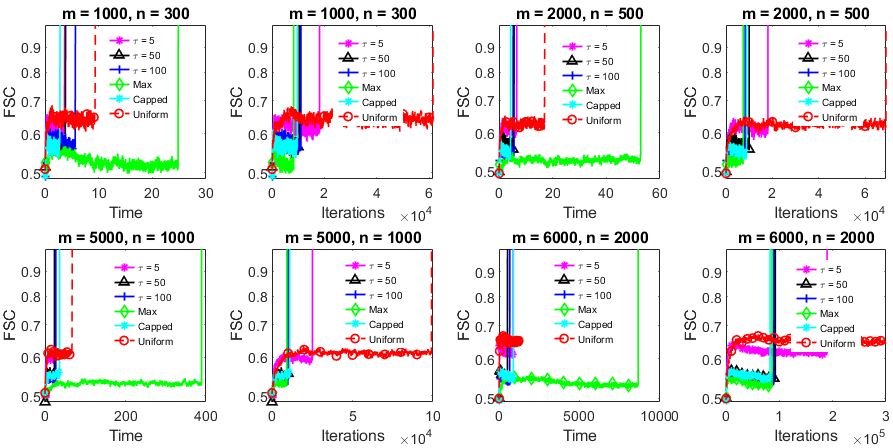}
    \caption{GK: comparison among sampling rules on Gaussian data, FSC vs time and No. of iterations.}
    \label{fig:2}
\end{figure}

\begin{figure}[htbp]
\centering
    \includegraphics[scale = 0.75]{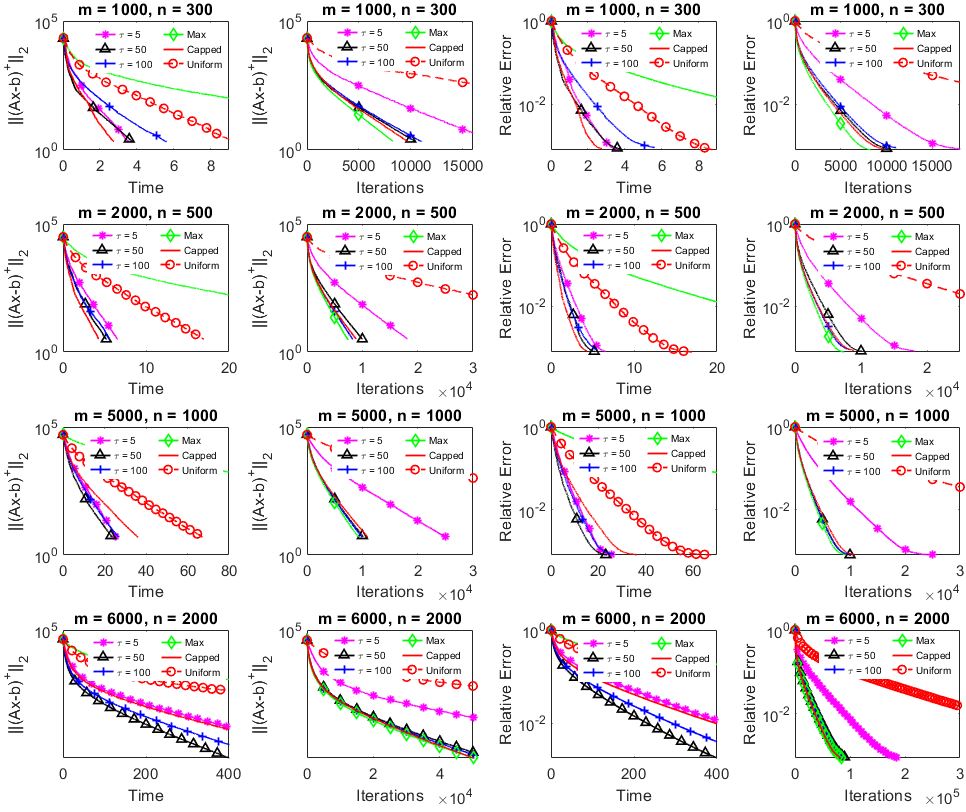}
    \caption{GK: comparison among sampling rules on Gaussian data, left 2 panels: Positive residual error $\|\left(Ax-b\right)^+\|_2$ vs time and No. of iterations, right 2 panels: relative error $\|x_k-x_{int}\|_B/\|x_0-x_{int}\|_B$ vs time and No. of iterations.}
    \label{fig:1}
\end{figure}

\paragraph{Comparison between sampling rules: GK on synthetic data.} In this subsection, we compare the performance of different sampling rules on the GK method. We will test these GK variants on randomly generated Gaussian test instances as well as four real-world sparse instances from Netlib LP test instances. For the Gaussian datasets, we consider six problems of sizes $1000 \times 300,\ 2000 \times 500, \ 5000 \times 1000, \ 6000 \times 2000$ respectively. In Figures \ref{fig:1} and \ref{fig:2}, we plot the performance measures of the selected variants of AK.

From the Figures, it is evident that the proposed greedy sampling rules heavily outperform both the uniform and maximum distance sampling rules. Furthermore, the performance of the proposed variants with greedy sampling and greedy capped sampling rules perform equality compared to each other. Another interesting point can be noted that whereas the uniform sampling rule takes the most number of iterations, on the other hand, the maximum distance rule takes the most time. It can be concluded from Figures \ref{fig:1} and \ref{fig:2} that the maximum distance rules perform poorly compared to other sampling rules.  

\begin{figure}[htbp]
\centering
    \includegraphics[scale = 0.76]{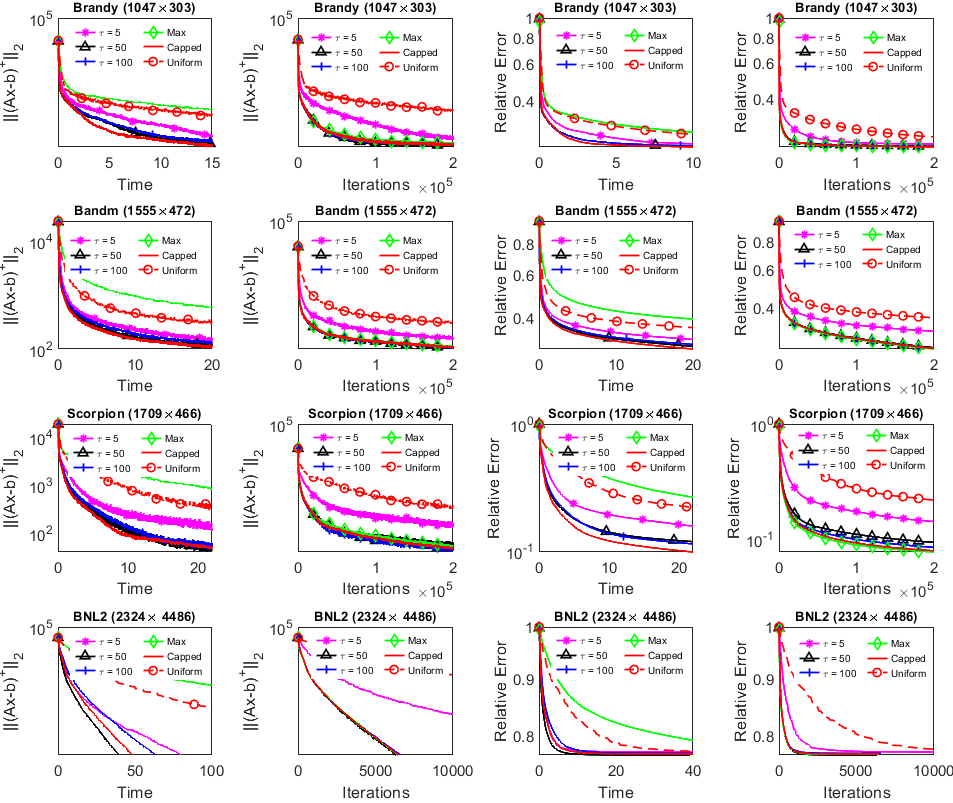}
    \caption{GK: comparison among sampling rules on Netlib LP instances, left 2 panels: Positive residual error $\|\left(Ax-b\right)^+\|_2$ vs time and No. of iterations, right 2 panels: relative error $\|x_k-x_{int}\|_B/\|x_0-x_{int}\|_B$ vs time and No. of iterations.}
    \label{fig:3}
\end{figure}

\paragraph{Comparison between sampling rules: GK on real data} Now, we compare the performance of the GK method with respect to different sampling rules on real data.  In Figure \ref{fig:3}, we plot comparison graphs for the following Netlib LP test instances: \texttt{lp\_brandy}, \texttt{lp\_bandm}, \texttt{lp\_scorpion} and \texttt{lp\_BNL2}. We consider $10^{-07}$ as the relative positive residual error tolerance for these problems, i.e., ($\|(Ax_k-b)^+\|_2/\|(Ax_0-b)^+\|_2 \leq 10^{-07}$). Similarly as random data, we see the same performance trend here as well. The performance of greedy sampling based GK variants are better compared to the maximum distance and uniform sampling rules.

\begin{figure}[htbp]
\centering
    \includegraphics[scale = 0.76]{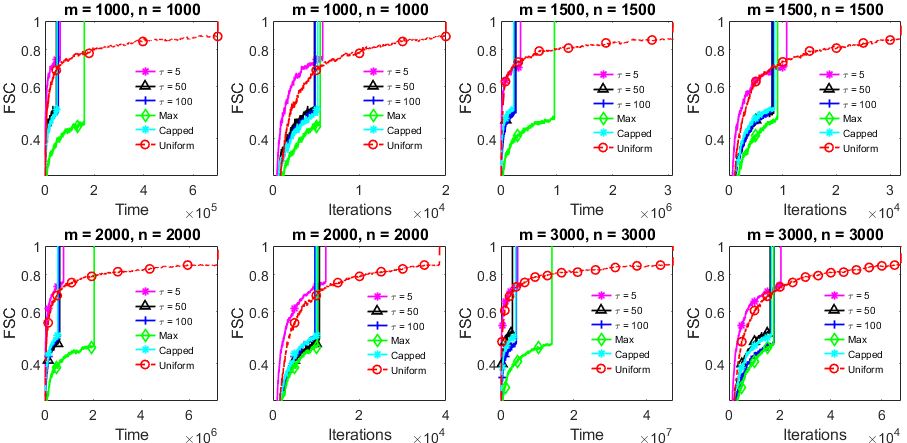}
    \caption{GCD: comparison among sampling rules on Gaussian data, FSC vs time and No. of iterations.}
    \label{fig:5}
\end{figure}

\begin{figure}[htbp]
\centering
    \includegraphics[scale = 0.74]{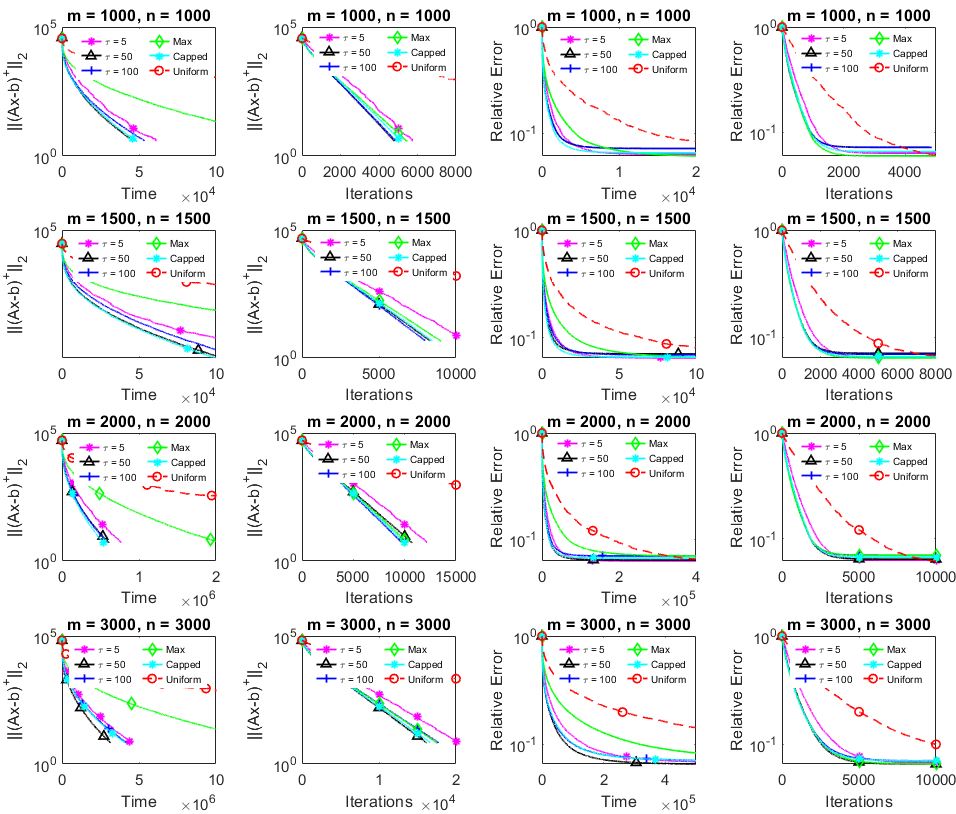}
    \caption{GCD: comparison among sampling rules on Gaussian data, left 2 panels: Positive residual error $\|\left(Ax-b\right)^+\|_2$ vs time and No. of iterations, right 2 panels: relative error $\|x_k-x_{int}\|_B/\|x_0-x_{int}\|_B$ vs time and No. of iterations.}
    \label{fig:4}
\end{figure}

\paragraph{Comparison between sampling rules: \textit{Greedy Co-ordinate Descent} GCD on synthetic data} In this subsection, we compare the performance of different sampling rules on the GCD method. We test GCD variants on four positive definite Gaussian problems of sizes $1000 \times 1000,\ 1500 \times 1500, \ 2000 \times 2000, \ 3000 \times 3000$ respectively. In Figures \ref{fig:4} and \ref{fig:5}, we plot the performance measures of the selected variants of CD. We see a slightly different trend compared to the GK methods. Here, the uniform sampling-based methods have worse performance. However, the GCD methods based on greedy sampling rules heavily outperform both the uniform rule and maximum distance rule-based methods.  From the solution quality performance graphs, we observe a similar trend as well.

\subsection{Greedy Methods with momentum}

In this subsection, we analyze the effect of momentum on the proposed GK and GCD methods equipped with the above-mentioned sampling rules. We carry out three types of experiments. In the first type, we discuss the effect of momentum parameter $\gamma$ on the selection of sketch sample size $\tau = |\phi_k(\tau)|$ as well as on the capped parameter $\theta$. In the second type, we compare the momentum-based GK methods with the basic method without momentum on both random and real-world test problems. In the third type, we compare the momentum-based CD methods with the basic method without momentum on random test instances. Take, $\tau_1 = m$ and $\tau_2 =1$. We perform the test for both GK and GCD methods.

\begin{figure}[htbp]
\centering
    \includegraphics[scale = 0.64]{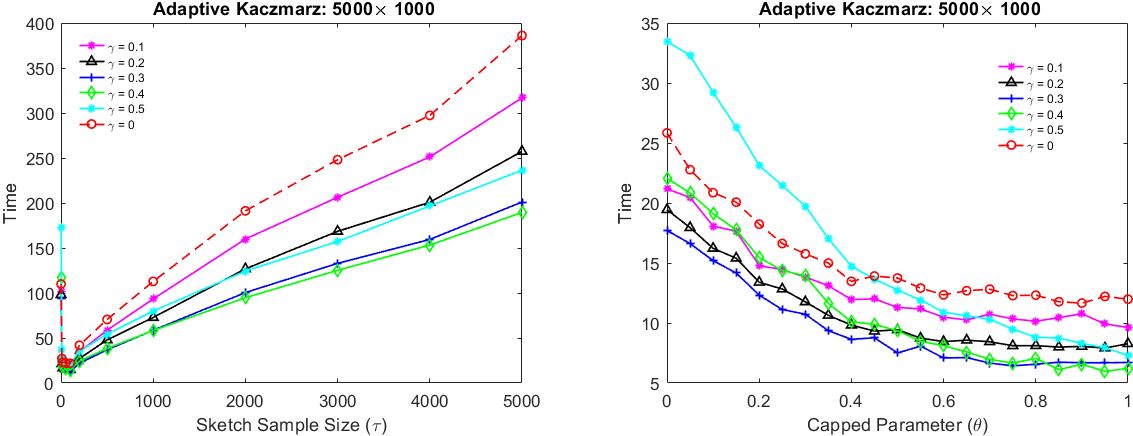}
    \caption{GK with momentum (effect of $\tau$ and $\theta$).}
    \label{fig:6}
\end{figure}

\begin{figure}[htbp]
\centering
    \includegraphics[scale = 0.63]{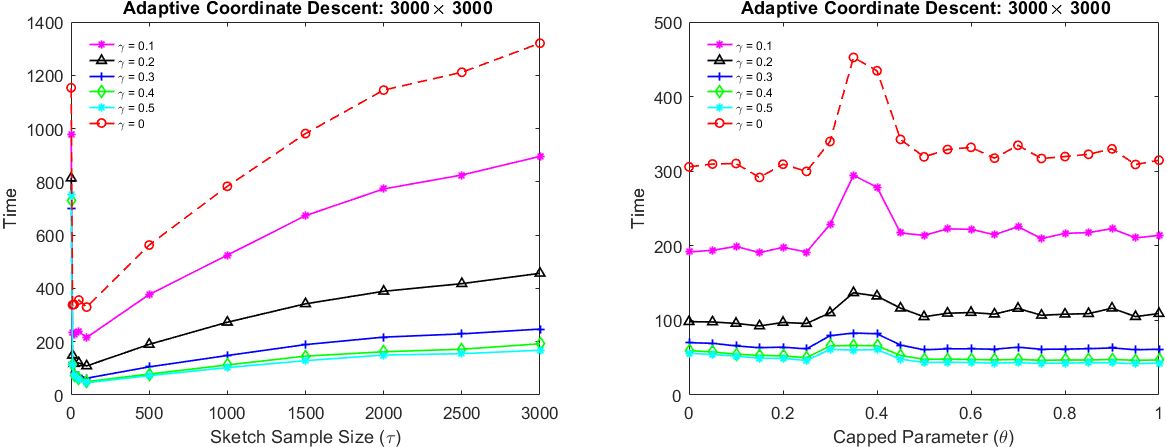}
    \caption{GCD with momentum (effect of $\tau$ and $\theta$).}
    \label{fig:60}
\end{figure}

In Figures \ref{fig:6} and \ref{fig:60}, we plot the momentum variants and the basic algorithm and compare the CPU time with respect to varying sketch sample size $\tau$ and capped parameter $\theta$. From Figure \ref{fig:6}, we see that for the GK method, the optimum sketch sample size $\tau$ occurs at $1 < \tau \ll m$ and capped parameter occurs at $0.5 \leq \theta \leq 1$. And the momentum variants performs well in comparison with the basic method irrespective of $\tau$ and $\theta$. From Figure \ref{fig:60}, we see that for the GCD method, the optimum sketch sample size $\tau$ occurs at $1 < \tau \ll m$. However, performance stays consistent with respect to the capped parameter $\theta$. Note that, similarly as before the momentum variants outperform the basic method. Moreover, it can be noted that for the GCD methods momentum methods heavily outperform the basic method compared to the GK methods.

\subsubsection{GK with momentum} 
To explore the findings of the previous subsection, here we test the GK momentum variants on individual sampling rules in detail. In Figures \ref{fig:7}-\ref{fig:10}, we plot \footnote{Please see Figures \ref{fig:16}-\ref{fig:23} in the Appendix section \ref{appendix:exp} for additional experiments with varying sketch sample size $\tau$ } the comparison measures graph for greedy rule with $\tau = 100$ and capped rule with $\theta = 0.5$. 

\begin{figure}[htbp]
\centering
    \includegraphics[scale = 0.73]{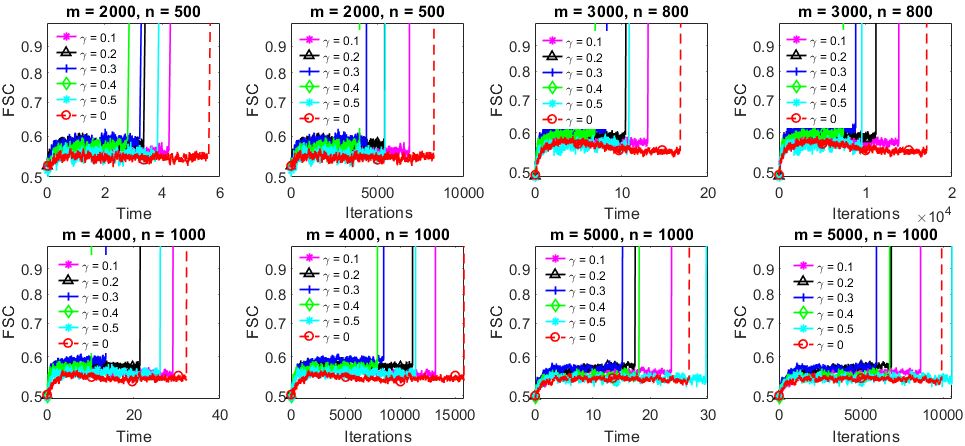}
    \caption{GK with momentum (sketch Sample size, $\tau = 100$): comparison among momentum variants on Gaussian data, FSC vs time and No. of iterations.}
    \label{fig:8}
\end{figure}

\begin{figure}[htbp]
\centering
    \includegraphics[scale = 0.75]{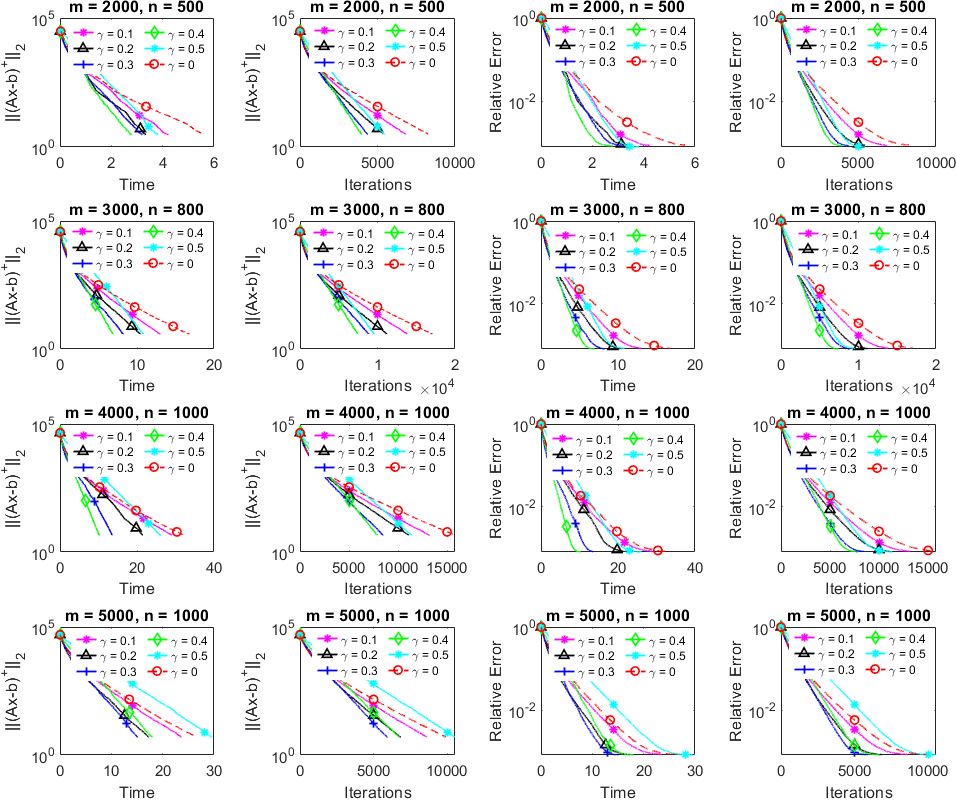}
    \caption{GK with momentum (sketch Sample size, $\tau = 100$): comparison among momentum variants on Gaussian data, left 2 panels: Positive residual error $\|\left(Ax-b\right)^+\|_2$ vs time and No. of iterations, right 2 panels: relative error $\|x_k-x_{int}\|_B/\|x_0-x_{int}\|_B$ vs time and No. of iterations.}
    \label{fig:7}
\end{figure}

\begin{figure}[htbp]
\centering
    \includegraphics[scale = 0.71]{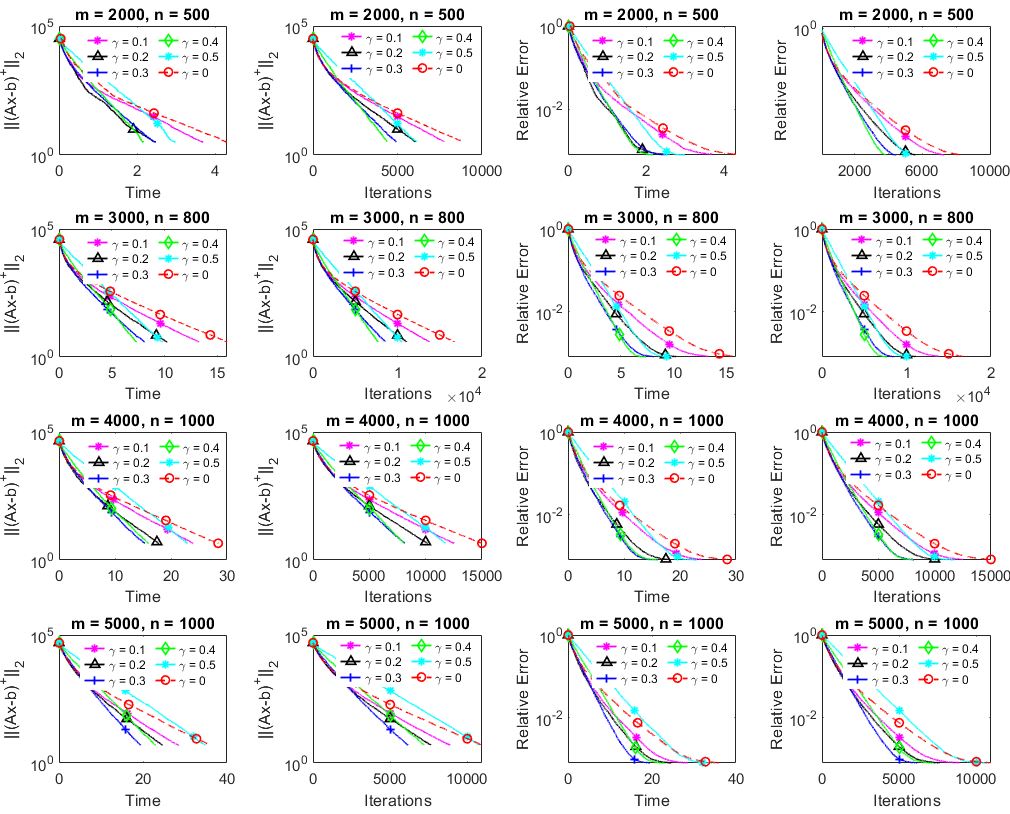}
    \caption{GK with momentum (capped rule, $\tau_1 = 1, \tau_2 = m, \theta = 0.5$): comparison among momentum variants on Gaussian data, left 2 panels: Positive residual error $\|\left(Ax-b\right)^+\|_2$ vs time and No. of iterations, right 2 panels: relative error $\|x_k-x_{int}\|_B/\|x_0-x_{int}\|_B$ vs time and No. of iterations.}
    \label{fig:9}
\end{figure}
\begin{figure}[htbp]
\centering
    \includegraphics[scale = 0.76]{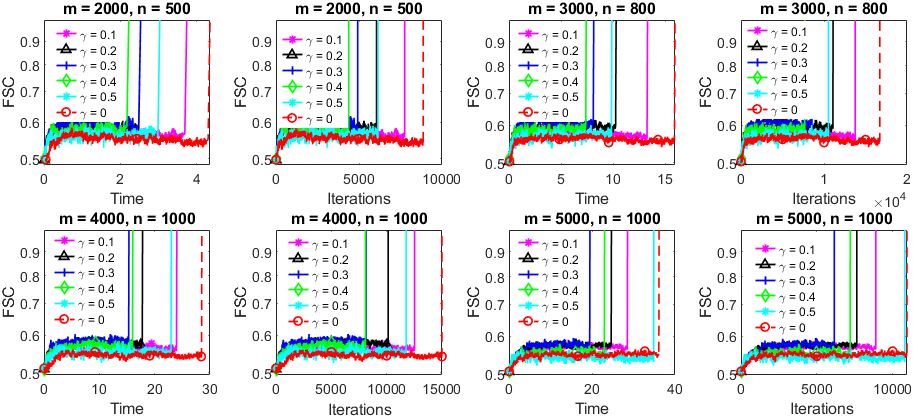}
    \caption{GK with momentum (capped rule, $\tau_1 = 1, \tau_2 = m, \theta = 0.5$): comparison among momentum variants on Gaussian data, FSC vs time and No. of iterations.}
    \label{fig:10}
\end{figure}

\begin{figure}[htbp]
\centering
    \includegraphics[scale = 0.76]{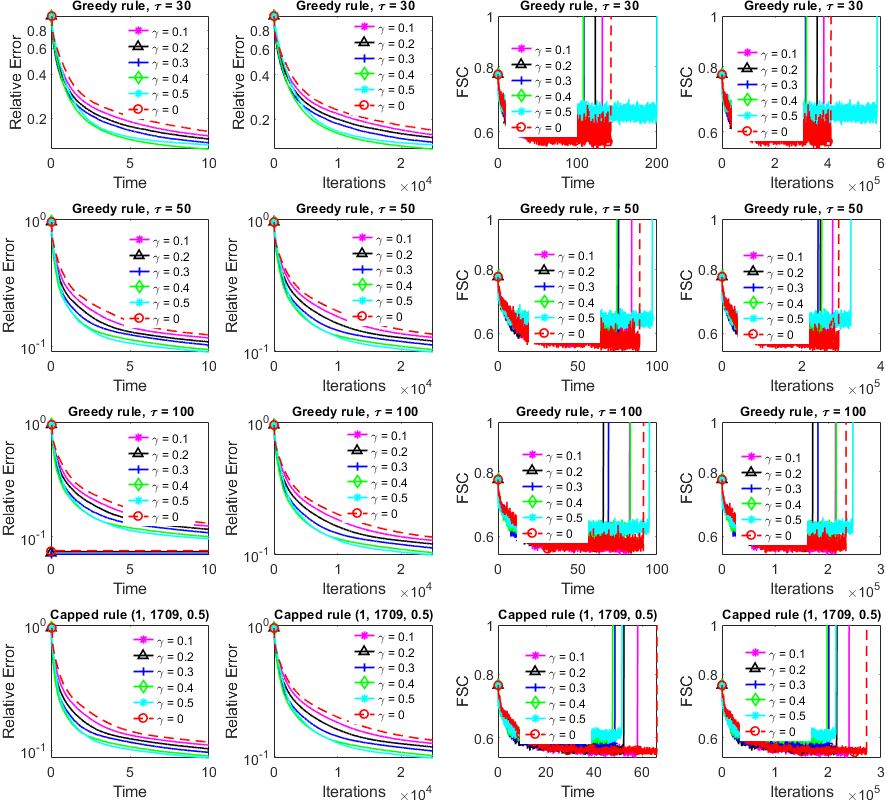}
    \caption{GK with momentum ( Greedy rule: $\tau = 30, 50, 100$ Capped rule, $\tau_1 = 1, \tau_2 = m, \theta = 0.5$) on Netlib \texttt{lp\_scorpion} ($1709 \times 466$) , left 2 panels: relative error $\|x_k-x_{int}\|_B/\|x_0-x_{int}\|_B$ vs time and No. of iterations, right 2 panels: FSC vs time and No. of iterations.}
    \label{fig:11}
\end{figure}

\paragraph{Netlib instance} In Figure \ref{fig:11}, we test the GK method with different sampling rules on sparse \texttt{lp\_scorpion} dataset.

From Figures \ref{fig:7}-\ref{fig:11} and \ref{fig:16}-\ref{fig:23} in the Appendix section \ref{appendix:exp}, we see that the proposed momentum variants outperform the basic GK method for all of the sampling rules considered on both real and random test instances.

\subsubsection{GCD with momentum}
To explore the findings of the previous subsection, here we test the GCD momentum variants on individual sampling rules in detail. In Figures \ref{fig:12}-\ref{fig:15}, we plot \footnote{Please see Figures \ref{fig:24}-\ref{fig:31} in the Appendix section \ref{appendix:exp} for additional experiments with varying sketch sample size $\tau$. } the comparison measures graph for greedy rule with $\tau = 100$ and capped rule with $\theta = 0.5$. From Figures \ref{fig:12}-\ref{fig:15} and \ref{fig:24}-\ref{fig:31} in the Appendix section \ref{appendix:exp}, we see that the proposed momentum variants outperform the basic GCD method for all of the sampling rules considered on both random Gaussian test instances.

\begin{figure}[htbp]
\centering
    \includegraphics[scale = 0.75]{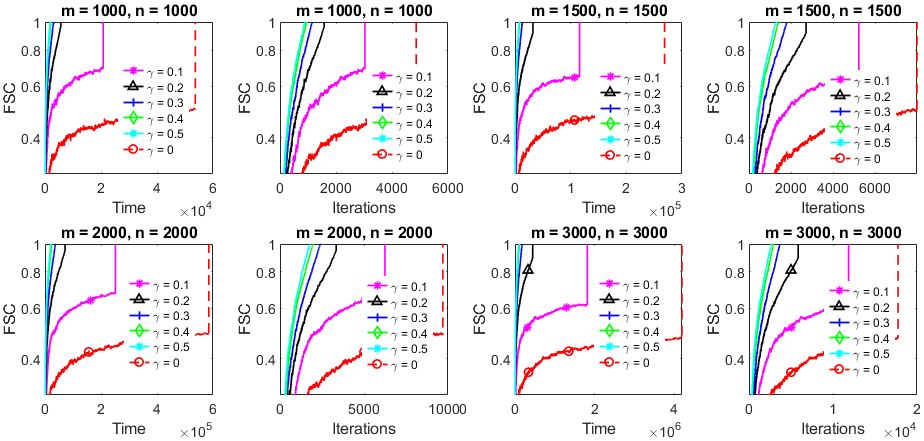}
    \caption{GCD with momentum (sketch Sample size, $\tau = 100$): comparison among momentum variants on Gaussian data, FSC vs time and No. of iterations.}
    \label{fig:13}
\end{figure}
\begin{figure}[htbp]
\centering
    \includegraphics[scale = 0.76]{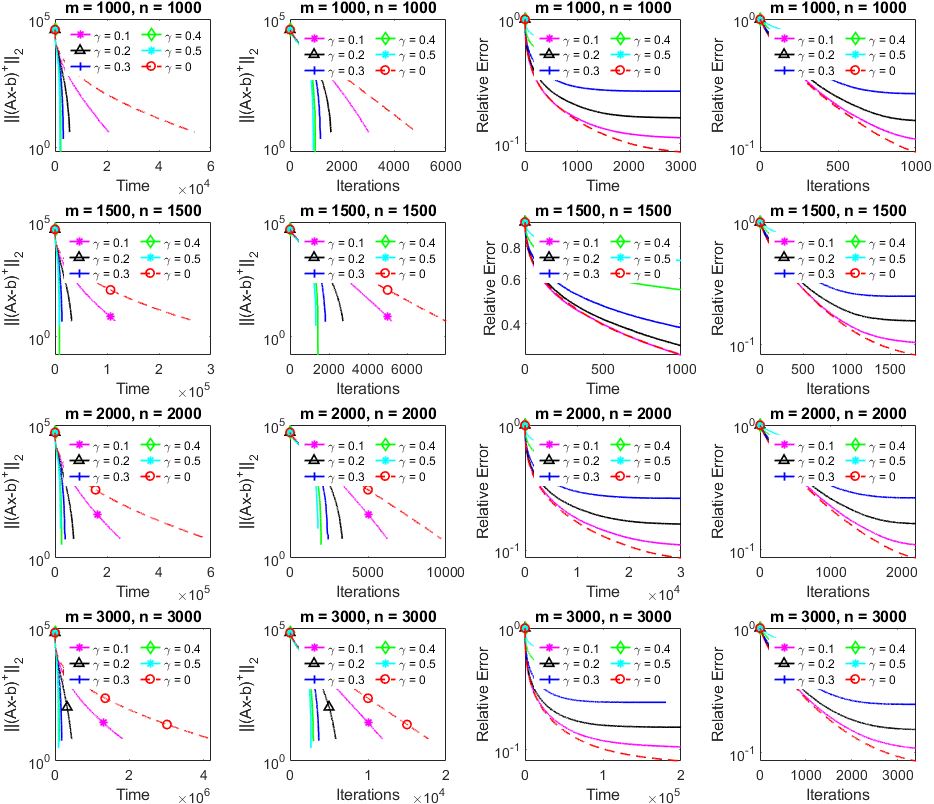}
    \caption{GCD with momentum (sketch Sample size, $\tau = 100$): comparison among momentum variants on Gaussian data, left 2 panels: Positive residual error $\|\left(Ax-b\right)^+\|_2$ vs time and No. of iterations, right 2 panels: relative error $\|x_k-x_{int}\|_B/\|x_0-x_{int}\|_B$ vs time and No. of iterations.}
    \label{fig:12}
\end{figure}

\begin{figure}[htbp]
\centering
    \includegraphics[scale = 0.74]{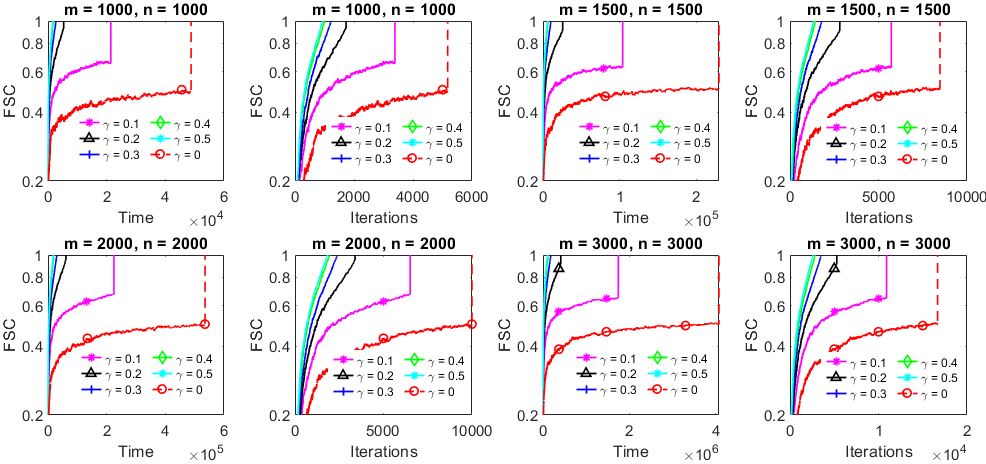}
    \caption{GCD with momentum (capped rule, $\tau_1 = 1, \tau_2 = m, \theta = 0.5$): comparison among momentum variants on Gaussian data, FSC vs time and No. of iterations.}
    \label{fig:15}
\end{figure}

\begin{figure}[htbp]
\centering
    \includegraphics[scale = 0.75]{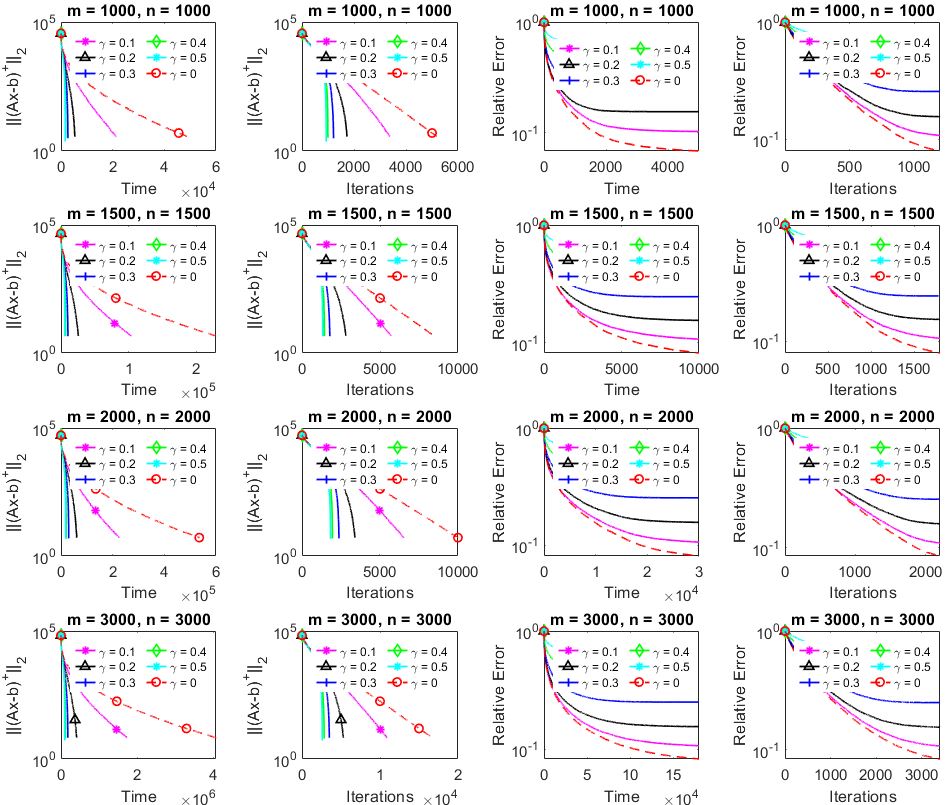}
    \caption{GCD with momentum (capped rule, $\tau_1 = 1, \tau_2 = m, \theta = 0.5$): comparison among momentum variants on Gaussian data, left 2 panels: Positive residual error $\|\left(Ax-b\right)^+\|_2$ vs time and No. of iterations, right 2 panels: relative error $\|x_k-x_{int}\|_B/\|x_0-x_{int}\|_B$ vs time and No. of iterations.}
    \label{fig:14}
\end{figure}

\newpage

\section{Conclusions}
\label{sec:conclusions}

In this work, we propose a \textit{Sketch \& Project} algorithmic framework equipped with greedy sampling strategies for solving linear feasibility problems. The proposed method synthesizes several well-known algorithms and their respective convergence results into one algorithm. Furthermore, we develop efficient algorithmic variants of the proposed method by incorporating the heavy ball momentum technique. We design a comprehensive numerical experimental setup to test the proposed greedy sampling rules as well as the momentum algorithms. For an unbiased conclusion about the performance and applicability, we test the proposed methods on random and real-world feasibility test instances. From our computational experiments, we conclude that the proposed greedy sampling rules heavily outperform the existing sampling rules. Furthermore, the proposed momentum variants accelerate the algorithmic performance of the basic method. We conclude the paper with some noteworthy future research directions. From the computational experiments, we find that the greedy sampling strategies based on sketched loss functions are very efficient. It is natural to think of the extension of the proposed methods in terms of 1) optimal sketched sample size $\tau^*$ selection based on information of matrix $A$, 2) greedy sketch sample size (sample size varies at each iteration). Another important extension of the proposed methods would be the design of efficient sparse variants.

\appendix

\section{Preliminary results}
\label{sec:prel}

\begin{mdframed}[backgroundcolor=gray!20,   topline=false,   bottomline =false,   rightline=false,   leftline=false]  \begin{lemma}
\label{lem0}
(Hoffman \cite{hoffman}, Theorem 4.4 in \cite{lewis:2010}) Let $x \in \R^n$ and $P$ be the feasible region, then there exists a constant $L > 0$ such that the following identity holds:
\begin{align*}
  d(x,P)^2 \leq L^2 \ \|(Ax-b)^+\|^2,
\end{align*}
\end{lemma} \end{mdframed}
where $L$ is the so-called Hoffman constant. Note that, when the system is consistent (i.e., there exists a unique $x^*$ such that $Ax = b$), $L$ can be calculated as follows:
\[L^2 = \frac{1}{\|A^{-1}\|^2} = \frac{1}{\lambda_{min}^{+}(A^TA)}.\]

\begin{mdframed}[backgroundcolor=gray!20,   topline=false,   bottomline =false,   rightline=false,   leftline=false]  \begin{lemma}
\label{lem:skmseq}
(Lemma 2.1 in \cite{haddock:2017})  Let $\{x_k\},\ \{y_k\}$ be real non-negative sequences such that $x_{k+1} > x_k > 0$ and $y_{k+1} \geq y_k \geq 0$, then
\begin{align*}
  \sum\limits_{k=1}^{n} x_k y_k \ \geq \   \sum\limits_{k=1}^{n} \overline{x} y_k, \quad \text{where} \ \ \overline{x} = \frac{1}{n}\sum\limits_{k=1}^{n} x_k.
\end{align*}
\end{lemma} \end{mdframed}

The following two Theorems deal with the decay of non-negative sequences that satisfy certain homogeneous recurrence inequalities.

\allowdisplaybreaks{\begin{mdframed}[backgroundcolor=gray!20,   topline=false,   bottomline =false,   rightline=false,   leftline=false]  \begin{theorem}
\label{th:seq2}
(Theorem 2 in \cite{morshed2020generalization}) Let the real sequences $H_k \geq 0$ and $F_k \geq 0$ satisfy the following recurrence relation:
\begin{align}
\label{t-1}
\begin{bmatrix}
H_{k+1} \\
F_{k+1} 
\end{bmatrix} & \leq   \begin{bmatrix}
\Pi_1 & \Pi_2 \\
\Pi_3  & \ \Pi_4
\end{bmatrix} \begin{bmatrix}
H_{k} \\
F_{k}
\end{bmatrix},
 \end{align}
where, $\Pi_1, \Pi_2, \Pi_3, \Pi_4 \geq 0$ such that the following relation
\begin{align}
    \label{t0}
 \Pi_1+ \Pi_4 < 1+ \min\{1, \Pi_1\Pi_4 - \Pi_2\Pi_3 \}
\end{align}
holds. Then the sequence $\{H_k\}$ and $ \{F_k\}$ converges and the following result holds:
\allowdisplaybreaks{\begin{align*}
  \begin{bmatrix}
H_{k+1}  \\[6pt]
F_{k+1} 
\end{bmatrix} & \leq   \begin{bmatrix}
\Pi_1 & \Pi_2 \\
\Pi_3  & \ \Pi_4
\end{bmatrix}^k \begin{bmatrix}
H_{1} \\
F_{1}
\end{bmatrix} =  \begin{bmatrix}
\Gamma_2 \Gamma_3 (\Gamma_1-1) \ \rho_1^{k}+ \Gamma_1 \Gamma_3 (\Gamma_2+1)\ \rho_2^{k} \\[6pt]
\Gamma_3 (\Gamma_1-1) \ \rho_1^{k}+ \Gamma_3 (\Gamma_2+1)\ \rho_2^{k}
\end{bmatrix} \ \begin{bmatrix}
H_{1} \\
F_1 
\end{bmatrix}.
\end{align*}}
where,
\allowdisplaybreaks{\begin{align}
    \label{t1}
    & \Gamma_1 = \frac{\Pi_1-\Pi_4+\sqrt{(\Pi_1-\Pi_4)^2+4\Pi_2\Pi_3}}{2\Pi_3},  \nonumber \\
    & \Gamma_2 = \frac{\Pi_1-\Pi_4-\sqrt{(\Pi_1-\Pi_4)^2+4\Pi_2\Pi_3}}{2\Pi_3}, \ \Gamma_3 = \frac{\Pi_3}{\sqrt{(\Pi_1-\Pi_4)^2+4\Pi_2\Pi_3}}, \nonumber \\
    & \rho_1 = \frac{1}{2} \left[\Pi_1+\Pi_4 - \sqrt{(\Pi_1-\Pi_4)^2+4\Pi_2\Pi_3}\right], \nonumber \\
    & \rho_2 = \frac{1}{2} \left[\Pi_1+\Pi_4 + \sqrt{(\Pi_1-\Pi_4)^2+4\Pi_2\Pi_3}\right].
\end{align}}
and $  \Gamma_1, \Gamma_3 \geq 0$ and $ 0 \leq |\rho_1| \leq \rho_2 < 1$.

\end{theorem} \end{mdframed}}

\begin{mdframed}[backgroundcolor=gray!20,   topline=false,   bottomline =false,   rightline=false,   leftline=false]  \begin{theorem}
\label{seq}
(Lemma 1 in \cite{ghadimi}) Let $\{H_{k}\}_{k\geq 0}$,  $\ \{F_{k}\}_{k\geq 0}$ and $\{G_{k}\}_{k\geq 0}$ be non-negative sequences of real numbers satisfying 
\begin{align}
\label{seq:1}
    H_{k+1} + \alpha_1 F_{k+1} \ \leq \  \beta_1 H_{k} + \beta_2 H_{k-1}+  \beta_3 F_{k},
\end{align}
 with constants $\beta_1, \beta_2,  \alpha_1 \geq 0$ and $\beta_3 \in \R$. Moreover, assume that
 \begin{align*}
     H_{1} = H_0, \quad \beta_1 + \beta_2 < 1, \quad \beta_3 < \alpha_1,
 \end{align*}
holds. Then the sequence $\{H_{k}\}_{k\geq 0}$ generated by \eqref{seq:1} satisfies
\begin{align}
    \label{seq:2}
    H_{k+1}+ \alpha H_k + \alpha_1 F_{k+1} \leq \ \rho^k \left[(1+\alpha)H_1+ \alpha_1 F_1\right],
\end{align}
where $ \alpha \geq 0$ and $ \rho \in [0,1)$ are given by
\begin{align*}
    \alpha = \max \left\{0, \frac{\beta_3}{\alpha_1}-\beta_1, \frac{-\beta_1+ \sqrt{\beta_1^2+4 \beta_2}}{2} \right\}, \quad \rho = \beta_1 + \alpha.
\end{align*}
\end{theorem} \end{mdframed}

Next, we discuss some well-known results from literature for developing feasibility certification bounds for linear feasibility problems. In our derivation of halting certification, we will use these results frequently. For a detailed discussion and implication of these results, we refer interested readers to the works \cite{haddock:2017,morshed2020generalization,morshed:momentum} and the references therein.

\begin{mdframed}[backgroundcolor=gray!20,   topline=false,   bottomline =false,   rightline=false,   leftline=false]  \begin{lemma}
\label{lem:skm1}
(Lemma 1 in \cite{haddock:2017}, Lemma 10 in \cite{ morshed2020generalization}) Define, $\theta(x) = \max_{i}(a_i^Tx-b_i)^+ $ as the maximum violation of point $x \in \R^n$ and the length of the binary encoding of a linear feasibility problem with rational data-points as
\begin{align*}
  \sigma = \sum \limits_{i}  \sum \limits_{j} \ln{\left(|a_{ij}|+1\right)} + \sum \limits_{i}  \ln{\left(|b_{i}|+1\right)} + \ln{(mn)} +2.
\end{align*}
Then if the rational system $Ax \leq b$ is infeasible, for any $x\in \R^n$, the maximum violation $\theta(x)$ satisfies the following lower bound:
\begin{align*}
    \theta(x) \ \geq \ \frac{2}{2^{\sigma}}.
\end{align*}
\end{lemma} \end{mdframed}

\begin{mdframed}[backgroundcolor=gray!20,   topline=false,   bottomline =false,   rightline=false,   leftline=false]  \begin{lemma}
\label{lem:skm4}
(\cite{KHACHIYAN:1980}) If the rational system $ Ax \leq b$ is feasible, then there is a feasible solution $x^{*}$ whose coordinates satisfy $|x^{*}_j| \leq \frac{2^{\sigma}}{2n}$ for $j = 1, ..., n$.
\end{lemma} \end{mdframed}

\section{Special Case}
\label{extracase}
Assume, $A \geq 0$. Then, we can derive the following method.

\paragraph{Momentum Sampling Co-ordinate Descent-Least Square (MSCD-LS)} Take, $q = n, \ B = A^TA^+, \ S= A^+e_i = A^+_i$. Then the ASPM method with the greedy sketched loss sampling (i.e., $i \sim \mathcal{G}(\tau)$ subsection \ref{ss}) resolves into the following update formula:
\begin{align}
\label{rcd-pd1}
    x_{k+1} = x_k- \delta \frac{[A_{i}^{+T}(Ax_k -b)]^+}{\|A^+_{i}\|^2} e_{i} + \gamma(x_k-x_{k-1}), 
\end{align}
where, $A_{i}$ is the $i^{th}$ column of matrix $A$ and $i = \argmax_{i \in \phi_k(\tau)}\frac{|\left(A_i^{+T}(Ax_k -b)\right)^+|^2}{\|A^+_{i}\|^2}  $ and $\phi_k(\tau)$ denotes the collection of $\tau$ columns chosen uniformly at random out of $n$ columns of the constraint matrix $A$. Using the above parameter choice in subsection \ref{ss}, we get $R^T = A^+$ and $\lambda_{\max}(B^{-\frac{1}{2}}A^TR^TRA B^{-\frac{1}{2}}) = \lambda_{\max}(A^TA^+)$. Take, $A_{\min} = \min_{i \in \{1,2,...,n\}} \|A_i^+\|^2$ and $A_{\max} = \max_{i \in \{1,2,...,n\}} \|A_i^+\|^2$. Then, we can calculate the constants $\mu_1$, $\mu_2$ as follows:
\begin{align}
\label{mu:scd-ls}
  \mu_1 =  \frac{1}{\sigma_1 A_{\max}}\min \left\{\frac{1}{n-\tau+1}, \frac{1}{n-s}\right\} \geq \frac{1}{n \sigma_1 A_{\max}}, \ \   \mu_2 = \min \left\{1, \frac{\tau \lambda_{\max}(A^TA^+) }{n A_{\min}} \right\}.
\end{align}
Here, $s$ is the number of zero entries in the residual $A^{+T}(Ax-b)^+$ $\sigma_1$ is the Hoffman constant. Using these parameter values, in the following Corollaries, we derive the convergence result for both \textit{Sampling Co-ordinate Descent-Least Square} (SCD-LS) and MSCD-LS methods. Note that, with the choice $\tau =1$ in SCD-LS, we can recover the \textit{Randomized Co-ordinate Descent-Least Square} (RCD-LS) method proposed in \cite{lewis:2010} for solving a system of linear equations.

\begin{mdframed}[backgroundcolor=gray!20,   topline=false,   bottomline =false,   rightline=false,   leftline=false]  \begin{corollary}
\label{cor:rcd-ls}
(New Theorem) Let, $x_k$ be the random iterate generated by the SCD-LS method with $0 < \delta < 2$. Then, the following identities
\begin{align*}
 \E[ d_{A^TA^+}(x_{k},\mathcal{X})^2]  \leq  \left(1-\frac{\eta}{n \sigma_1 A_{\max}}\right)^{k}  d_{A^TA^+}(x_0,\mathcal{X})^2,
\end{align*}
and 
\begin{align*}
   \E\left[\frac{\big |[A_{i}^{+T}(Ax_k -b)]^+\big |^2}{\|A^+_{i}\|^2}\right]  \leq  \frac{\tau \lambda_{\max}(A^TA^+) }{n A_{\min}} \left(1-\frac{\eta}{n \sigma_1 A_{\max}}\right)^{k} d_{A^TA^+}(x_0,\mathcal{X})^2,
\end{align*}
hold. Also the average iterate $\Tilde{x}_k = \sum \nolimits_{l=0}^{k-1} x_l$ for all $k \geq 1$ satisfies
\begin{align*}
    \E[d_{A^TA^+}(\Tilde{x}_k,\mathcal{X})^2]  \leq \frac{n \sigma_1 A_{\max}}{2\delta k(2-\delta)} \  d_{A^TA^+}(x_0,\mathcal{X})^2,
\end{align*}
and
\begin{align*}
  \E\left[\frac{\big |[A_{i}^{+T}(Ax_k -b)]^+\big |^2}{\|A^+_{i}\|^2}\right]  \leq \frac{d_{A^TA^+}(x_0,\mathcal{X})^2}{2\delta k(2-\delta)}.
\end{align*}
\end{corollary} \end{mdframed}

\begin{proof}
Consider, $\mu_1$ and $\mu_2$ values from \eqref{mu:scd-ls} in Theorem \ref{th:b1}. Then, with simplification, we get the result of Corollary \ref{cor:rcd-ls}.
\end{proof}

\begin{mdframed}[backgroundcolor=gray!20,   topline=false,   bottomline =false,   rightline=false,   leftline=false]  \begin{corollary} 
\label{cor:mscd-ls1}
(New Theorem) Let $\{x_k\}$ be the sequence of random iterates generated by the MSCD-LS algorithm starting with $x_0 = x_1 \in \R^n$. With $0 < \delta < 2$, the sequence of iterates $\{x_k\}$ converges and the following result holds:
{\allowdisplaybreaks
\begin{align*}
\E \begin{bmatrix}
d_{A^TA^+}(x_{k+1}, \mathcal{X})  \\[6pt]
\|x_{k+1}-x_k\|_{A^TA^+}
\end{bmatrix} & \leq \begin{bmatrix}
-\Gamma_2 \Gamma_3 \ \rho_1^{k}+ \Gamma_1 \Gamma_3 \ \rho_2^{k} \\[6pt]
- \Gamma_3 \ \rho_1^{k}+ \Gamma_3 \ \rho_2^{k} 
\end{bmatrix} \ d_{A^TA^+}(x_{0}, \mathcal{X})  \nonumber \\
& \leq \begin{bmatrix}
 \rho_2^{k} \\[6pt]
2 \Gamma_3 \ \rho_2^{k} 
\end{bmatrix} \ d_{A^TA^+}(x_{0}, \mathcal{X}).
\end{align*}}
\end{corollary} \end{mdframed}

\begin{proof}
Consider, $\mu_1$ and $\mu_2$ values from \eqref{mu:scd-ls} in Theorem \ref{th:1}. Then, with simplification, we get the result of Corollary \ref{cor:mscd-ls1}.
\end{proof}

\begin{mdframed}[backgroundcolor=gray!20,   topline=false,   bottomline =false,   rightline=false,   leftline=false]  \begin{corollary} 
\label{cor:mscd-ls2}
(New Theorem) Let $\{x_k\}$ be the sequence of random iterates generated by the MSCD-LS algorithm with $x_0 = x_1 \in \R^n$. With $0 < \delta < 2$, the sequence of iterates $\{x_k\}$ converges and the following results hold:
\begin{align*}
  \E [d_{A^TA^+}(x_{k+1},\mathcal{X})^2]  \leq  \rho^k (1+\alpha) \  d_{A^TA^+}(x_0,\mathcal{X})^2 ,
\end{align*}
and
\begin{align*}
  \E [f(x_{k+1})]  \leq \frac{\mu_2(1+\alpha)}{2} \rho^k \ d_{A^TA^+}(x_0,\mathcal{X})^2, 
\end{align*}
where, $\alpha \geq 0$,  $0 < \rho < 1$ are provided in \eqref{mom28} \eqref{mom29} respectively.
\end{corollary} \end{mdframed}
\begin{proof}
Consider, $\mu_1$ and $\mu_2$ values from \eqref{mu:scd-ls} in Theorem \ref{th:mom2}. Then, with simplification, we get the result of Corollary \ref{cor:mscd-ls2}.
\end{proof}

\section{Proofs}
\label{sec:prf}

\subsection[Proof of Thm]{Proof of \cref{neco}}
\label{pstoc:reg}
\begin{proof}
We already have $\mathcal{X} \subseteq \mathcal{X}'$. That means we just need to show that the relation $\mathcal{X}' \subseteq \mathcal{X}$ holds. For that, first assume $y \in \mathcal{X}'$. This means $y$ solves the stochastic optimization problem of \eqref{prob:stoc}. Now, it can be easily shown that the optimum value of optimization problem \eqref{prob:stoc} is zero, i.e., $f^* = \min_x f(x) = \min_x \E[f_i(x) \ | \ i \sim \mathcal{R}] = \min_x \E[d_B(x,\mathcal{X}_{S_i}) \ | \ i \sim \mathcal{R}] = 0$. Therefore, we can say for any $y \in \mathcal{X}'$, we have $\E[d_B(y,\mathcal{X}_{S_i}) \ | \ i \sim \mathcal{R}] = 0$. Now, we assume that the property of \eqref{stoc:reg} holds. Then, we have
\begin{align*}
    \mu \ d_B(y,\mathcal{X})^2 \leq \E [d_B(y,\mathcal{X}_{S_i})^2 \ | \ i \sim \mathcal{R}] = 0, \quad \Rightarrow  d_B(y,\mathcal{X})^2 = 0.
\end{align*}
That implies $y \in \mathcal{X}$. This proves the Lemma.
\end{proof}

\subsection[Proof of Thm]{Proof of \cref{1}}
\label{p1}
\begin{proof}
The first and second part of the Lemma follows from the following identity:
\begin{align*}
    d_B(x,\mathcal{X}_{S_i})^2 = \frac{\big |[S_{i}^T(Ax-b)]^+\big |^2}{\|A^TS_{i}\|^2_{B^{-1}}} = 2 f_i(x).
\end{align*}
The third part follows from the fact that for this specific setup the differentiation operator and expectation operator can be interchanged.
\end{proof}

\subsection[Proof of Thm]{Proof of \cref{2}}
\label{p2}
From the definition, we have
  \begin{align*}
      \big \langle  \bar{x}-x, \nabla^{B} f_i(x) \big \rangle_{B} & = \frac{\left[S_{i}^T(Ax-b)\right]^+}{\|A^TS_{i}\|^2_{B^{-1}}}  S_{i}^T A B^{-1} B (\bar{x}-x)  \leq \frac{\left[S_{i}^T(Ax-b)\right]^+}{\|A^TS_{i}\|^2_{B^{-1}}}  S_{i}^T (b-Ax) \\
      &  \leq -\frac{\big | \left[S_{i}^T(Ax-b)\right]^+ \big |^2}{\|A^TS_{i}\|^2_{B^{-1}}} = - 2 f_i(x).
  \end{align*}
This proves the first part of the Lemma. For proving the second part, let's take $\bar{x} = \mathcal{P}^B_{\mathcal{X}}(x)$ then we get the following:
\begin{align}
 \big \langle  x- \mathcal{P}^B_{\mathcal{X}}(x), \E[\nabla^{B} f_i(x)] \big  \rangle_{B} & = \big \langle  x- \mathcal{P}^B_{\mathcal{X}}(x), \E\left[x- \mathcal{P}^B_{\mathcal{X}_{S_i}}(x)\right] \big \rangle_{B}  \nonumber \\
  \geq  2 \E[f_{i}(x)] & = 2 f(x) = \E \left[\big \|x- \mathcal{P}^B_{\mathcal{X}_{S_i}}(x)\big\|^2_B\right] = \E \left[\big \|\nabla^{B} f_i(x)\big\|^2_B\right].
\end{align}
Using Cauchy–Schwarz inequality, we can get the following
\begin{align}
    2 f(x) \leq \big \langle  x- \mathcal{P}^B_{\mathcal{X}}(x), \E[\nabla^{B} f_i(x)] \big  \rangle_{B} \leq \|x- \mathcal{P}^B_{\mathcal{X}}(x)\|_B \|\E[\nabla^{B} f_i(x)]\|_B.
\end{align}
This proves the second part of Lemma.

\subsection[Proof of Thm]{Proof of \cref{th:upper}}
\label{pupper}
First, note that if the following identity holds 
\begin{align*}
  \big \| \E [\nabla^B f_i(x) \ | \ i \sim \mathcal{R}]\big \|_B^2 & \leq 2 \mu_2 \ \E[f_i(x) \ | \ i \sim \mathcal{R}],
\end{align*}
for some $\mu_2 \geq 0$, then we have
\begin{align*}
   f(x)  \leq  \frac{1}{2} d_B(x,\mathcal{X}) \ \|\E[\nabla^{B} f_i(x)\ | \ i \sim \mathcal{R} ]\|_B \leq \frac{\sqrt{2\mu_2}}{2} d_B(x,\mathcal{X}) \sqrt{f(x)},
\end{align*}
here, we used part 2 of Lemma \ref{2}. Simplifying further, we get the required bound for $f(x)$. This implies, we need to show there exists $\mu_2 \geq 0$ such that the following holds:
\begin{align}
\label{p}
  \big \| \E [\nabla^B f_i(x) \ | \ i \sim \mathcal{R}]\big \|_B^2 & \leq 2 \mu_2 \ \E[f_i(x) \ | \ i \sim \mathcal{R}].
\end{align}
It can be noted that, for any $x \in \mathcal{X} = \{x: Ax \leq b\}$, equation \eqref{p} is trivially true for any $\mu_2 \geq 0$. This implies we just need to consider the case $x$ with $Ax \nleq b$. Let's define, $\mathcal{T}(x) = \{i \in \{1,2,...,q\} \ | \ S_i^T(Ax-b) > 0\}$ and denote the following matrices:
\begin{align}
   W(x) = \sum \limits_{i \in \mathcal{T}(x)} p_i \frac{S_iS_i^T}{\|A^TS_i\|^2_{B^{-1}}}  \quad \text{and} \quad W = \sum \limits_{i =1}^{q} p_i \frac{S_iS_i^T}{\|A^TS_i\|^2_{B^{-1}}}.
\end{align}
It can be easily checked that matrices $W$ and $W(x)$ are positive semi-definite and $W(x) \preceq W$ for all $x$ with $Ax \nleq b$. Now, we have
\begin{align}
    \big \| \E [\nabla^B f_i(x) \ | \ i \sim \mathcal{R}]\big \|_B^2 & =  \Big \| B^{-1} A^T   \E \left[ \frac{S_i\left[S_i^T(Ax-b)\right]^+}{\|A^TS_i\|^2_{B^{-1}}}   \ | \ i \sim \mathcal{R} \right]\Big \|_B^2  \nonumber \\
    &  = \Big \| B^{-1} A^T   \left[ \sum \limits_{i \in \mathcal{T}(x)} p_i \frac{S_iS_i^T}{\|A^TS_i\|^2_{B^{-1}}} \right] (Ax-b)\Big \|_B^2  = \big \| B^{-1} A^T W(x)(Ax-b) \big \|^2_B.
\end{align}
Similarly, we have
\begin{align}
    2 \E[f_i(x) \ | \ i \sim \mathcal{R}] & = \E \left[\big \| \nabla^{B} f_i(x) \big \|^2_B \ | \ i \sim \mathcal{R} \right]  = \E \left[ \frac{\Big | \left[S_i^T(Ax-b)\right]^+ \Big |^2}{\|A^TS_i\|^2_{B^{-1}}} \ | \ i \sim \mathcal{R} \right] \nonumber \\
    & = (Ax-b)^T \left[ \sum \limits_{i \in \mathcal{T}(x)} p_i \frac{S_iS_i^T}{\|A^TS_i\|^2_{B^{-1}}} \right] (Ax-b)  = (Ax-b)^T W(x) (Ax-b).
\end{align}
Take, $y = \sqrt{W(x)} (Ax-b) \neq 0$ \footnote{When $y = \sqrt{W(x)} (Ax-b) = 0$ the required inequality of \eqref{p} holds trivially.}. Then, for all $x$ with $Ax \nleq b$ the following holds:
\begin{align}
    \frac{\big \| \E [\nabla^B f_i(x) \ | \ i \sim \mathcal{R}]\big \|_B^2}{2 \E[f_i(x) \ | \ i \sim \mathcal{R}]} & = \frac{\big \| B^{-1} A^T W(x)(Ax-b) \big \|^2_B}{(Ax-b)^T W(x) (Ax-b)} = \frac{\|B^{-1} A^T \sqrt{W(x)}y\|^2_B}{\|y\|^2} = \frac{\|B^{-\frac{1}{2}} A^T \sqrt{W(x)}y\|^2}{\|y\|^2} \nonumber \\
    & \leq \|B^{-\frac{1}{2}} A^T \sqrt{W(x)}\|^2   =  \lambda_{\max} \left(B^{-\frac{1}{2}} A^T W(x)A B^{-\frac{1}{2}}\right).
\end{align}
As $W(x) \preceq W$, then we must have $A^T W(x)A \preceq A^T W A$. That implies, we have the following
\begin{align}
   \lambda_{\max} \left(B^{-\frac{1}{2}} A^T W(x)A B^{-\frac{1}{2}}\right) \leq \lambda_{\max} \left(B^{-\frac{1}{2}} A^T W A B^{-\frac{1}{2}}\right).
\end{align}
Furthermore, it can be noted that, there exists $x$ such that $\mathcal{T}(x) = \{1,2,...,q\} $. Therefore, we have
\begin{align}
   \lambda_{\max} \left(B^{-\frac{1}{2}} A^T W A B^{-\frac{1}{2}}\right) = \lambda_{\max} \left(B^{-\frac{1}{2}} A^T \E\left[\frac{S_iS_i^T}{\|A^TS_i\|^2_{B^{-1}}} \ | \ i \sim \mathcal{R}] \right] A B^{-\frac{1}{2}}\right).
\end{align}
Moreover, as the function $\lambda_{\max}$ is convex over the space of positive semi-definite matrices, using Jensen's inequality we have the following:
\begin{align}
    \lambda_{\max} & \left(B^{-\frac{1}{2}} A^T \E\left[\frac{S_iS_i^T}{\|A^TS_i\|^2_{B^{-1}}}\ | \ i \sim \mathcal{R}\right] A B^{-\frac{1}{2}}\right)  \leq \E \left[ \lambda_{\max} \left(B^{-\frac{1}{2}} A^T \frac{S_iS_i^T}{\|A^TS_i\|^2_{B^{-1}}} A B^{-\frac{1}{2}}\right) \ | \ i \sim \mathcal{R} \right].
\end{align}
Now, denote $T_i = B^{-\frac{1}{2}} A^T \frac{S_iS_i^T}{\|A^TS_i\|^2_{B^{-1}}} A B^{-\frac{1}{2}}$, then the following identity holds: 
\begin{align}
 T_i^2  = \frac{B^{-\frac{1}{2}} A^T S_i \left(S_i^T A B^{-1} A^T S_i \right)S_i^T A B^{-\frac{1}{2}}}{\|A^TS_i\|^4_{B^{-1}}}  = \frac{B^{-\frac{1}{2}} A^T S_iS_i^T A B^{-\frac{1}{2}}}{\|A^TS_i\|^2_{B^{-1}}} = T_i.
\end{align}
Therefore, $\lambda_{\max}(T_i) = 1$. Now, denote, $\mu_2 =  \lambda_{\max} \left(B^{-\frac{1}{2}} A^T W A B^{-\frac{1}{2}}\right)$. Then, using we have
\begin{align}
    \mu_2 =  \lambda_{\max} \left(B^{-\frac{1}{2}} A^T Z A B^{-\frac{1}{2}}\right) \leq \E \left[ \lambda_{\max} \left(T_i\right) \right] \leq 1.
\end{align}
Which proves the Lemma.

\subsection[Proof of Thm]{Proof of \cref{lem:mu}}
\label{pmu}
From Lemma \ref{th:upper}, we have
\begin{align}
    \mu_2 =  \lambda_{\max} \left(B^{-\frac{1}{2}} A^T \E\left[\frac{S_iS_i^T}{\|A^TS_i\|^2_{B^{-1}}} \ | \ i \sim \mathcal{G}(\tau)] \right] A B^{-\frac{1}{2}}\right).
\end{align}
Using the construction, we have
\begin{align*}
   B^{-\frac{1}{2}} A^T \E\left[\frac{S_iS_i^T}{\|A^TS_i\|^2_{B^{-1}}} \ | \ i \sim \mathcal{G}(\tau)] \right] & A B^{-\frac{1}{2}}  \overset{\eqref{def:exp}}{=} \frac{1}{\binom{q}{\tau}} \sum\limits_{j = 0}^{q-\tau} \binom{\tau-1+j}{\tau-1} \frac{B^{-\frac{1}{2}} A^T S_{\underline{\mathbf{i_j}}}S_{\underline{\mathbf{i_j}}}^T A B^{-\frac{1}{2}}}{\|A^TS_{\underline{\mathbf{i_j}}}\|^2_{B^{-1}}} \\
    & \preceq \ \frac{1}{\omega_1 \binom{q}{\tau}} \sum\limits_{j = 0}^{q-\tau} \binom{\tau-1+j}{\tau-1} B^{-\frac{1}{2}} A^T S_{\underline{\mathbf{i_j}}}S_{\underline{\mathbf{i_j}}}^T A B^{-\frac{1}{2}} \preceq \ \frac{\binom{q-1}{\tau-1}}{\omega_1 \binom{q}{\tau}} \sum\limits_{j = 0}^{q-\tau}  B^{-\frac{1}{2}} A^T S_{\underline{\mathbf{i_j}}}S_{\underline{\mathbf{i_j}}}^T A B^{-\frac{1}{2}} \\
    & \preceq \ \frac{\tau}{\omega_1 q}  B^{-\frac{1}{2}} A^T  \sum\limits_{j = 1}^{q}  S_{j}S_{j}^T  A B^{-\frac{1}{2}} =  \ \frac{\tau}{\omega_1 q}  B^{-\frac{1}{2}} A^T  R^T RA B^{-\frac{1}{2}}.
\end{align*}
Now, considering Lemma \ref{th:upper}, we get the required result. Similarly, we have
\begin{align*}
   \E\left[d_B(x,\mathcal{X}_{i})^2 \ | \ i \sim \mathcal{G}(\tau)] \right] & = \E\left[\frac{\big | \left[S_i^T(Ax-b)\right]^+ \big |^2}{\|A^TS_i\|^2_{B^{-1}}} \ | \ i \sim \mathcal{G}(\tau)] \right]  \overset{\eqref{def:exp}}{=} \frac{1}{\binom{q}{\tau}} \sum\limits_{j = 0}^{q-\tau} \binom{\tau-1+j}{\tau-1} \frac{\big | \left[S_{\underline{\mathbf{i_j}}}^T(Ax-b)\right]^+ \big |^2}{\|A^TS_{\underline{\mathbf{i_j}}}\|^2_{B^{-1}}} \\
   & \overset{\text{Lemma} \ \ref{lem:skmseq}}{\geq} \frac{1}{\binom{q}{\tau}} \sum\limits_{j = 0}^{q-\tau}  \frac{\sum\limits_{l = 0}^{q-\tau} \binom{\tau-1+l}{\tau-1}}{q-\tau +1} \frac{\big | \left[S_{\underline{\mathbf{i_j}}}^T(Ax-b)\right]^+ \big |^2}{\|A^TS_{\underline{\mathbf{i_j}}}\|^2_{B^{-1}}} \\
    & \geq \frac{1}{\omega_2 (q-\tau+1)} \sum\limits_{j = 0}^{q-\tau} \big | \left[S_{\underline{\mathbf{i_j}}}^T(Ax-b)\right]^+ \big |^2 \\
     \geq & \frac{1}{\omega_2 (q-\tau+1)} \min \left\{1, \frac{q-\tau+1}{q-s}\right\} \sum\limits_{j = 1}^{q} \big | \left[S_j^T(Ax-b)\right]^+ \big |^2 \\
     = & \min \left\{\frac{1}{\omega_2(q-\tau+1)}, \frac{1}{\omega_2(q-s)}\right\} \|\left[R(Ax-b)\right]^+\|^2.
\end{align*}
The quantity $s$ denotes the number of zero entries in the vector $\left[R(Ax-b)\right]^+$ (i.e., $s = q- \|\left[R(Ax-b)\right]^+\|_0$, where $\|\cdot\|_0$ denotes the zero norm). Now, there exists some positive constant $\sigma > 0$ such that the following identity holds
\begin{align}
   d_B(x,\mathcal{X})^2 \leq \sigma  \big \| \left[R(Ax-b)\right]^+ \big \|^2, 
\end{align}
for all $x \in \R^n$. The constant $\sigma$ is the so-called Hoffman constant. Using the Hoffman bound in we have the following
\begin{align}
     \E\left[d_B(x,\mathcal{X}_{S})^2\right]  \geq   \frac{1}{\sigma} \min \left\{\frac{1}{\omega_2(q-\tau+1)}, \frac{1}{\omega_2(q-s)}\right\} d_B(x,\mathcal{X})^2 \geq \frac{1}{q \sigma \omega_2} \ d_B(x,\mathcal{X})^2.
\end{align}
This proves the second part of the Lemma.

\subsection[Proof of Thm]{Proof of \cref{lem:c}}
\label{pc}
First, note that from the expectation calculation, we have
 \begin{align}
 \label{cap1}
    \E[f_i(x) \ | \ i \sim \mathcal{C}(\theta, \tau_1,\tau_2)] & = \sum \limits_{j \in \mathcal{W}} p_j f_j(x) \nonumber \\
    & \geq \theta  \E[f_j(x) \ | \ j \sim \mathcal{G}(\tau_1)] + (1-\theta) \E[f_j(x) \ | \ j \sim \mathcal{G}(\tau_2)] \nonumber \\
    & \geq \frac{\theta \mu_1(\tau_1)}{2}d_B(x,\mathcal{X})^2 + \frac{(1-\theta) \mu_1 (\tau_2)}{2} d_B(x,\mathcal{X})^2 \nonumber \\
    & = \frac{ \theta \mu_1(\tau_1)+ (1-\theta) \mu_1 (\tau_2)}{2} d_B(x,\mathcal{X})^2.
\end{align}
Similarly, we have
 \begin{align}
 \label{cap2}
    \E[f_i(x) \ | \ i \sim \mathcal{C}(\theta, \tau_1,\tau_2)]  = \sum \limits_{j \in \mathcal{W}} p_j f_j(x) \leq \max_{i \in \{1,2,...,q\}} f_i(x) \leq \frac{\mu_2(q)}{2} d_B(x,\mathcal{X})^2.
\end{align}
Combining \eqref{cap1} and \eqref{cap2}, we get the required Lemma.

\subsection[Proof of Thm]{Proof of \cref{th:b1}}
\label{pb1}
Since, $\mathcal{P}_{\mathcal{X}}^{B}(x_k)  \in \mathcal{X}$, From \eqref{eq:sk1} we have
{\allowdisplaybreaks
\begin{align}
\label{eq:b1}
d_B(x_{k+1},\mathcal{X})^2 & = \|x_{k+1}-  \mathcal{P}_{\mathcal{X}}^{B}(x_{k+1})\|^2_B  \overset{\text{Lemma} \ \ref{lem:distance}}{ \leq} \|x_{k+1}-  \mathcal{P}_{\mathcal{X}}^{B}(x_{k})\|^2_B \nonumber \\
& = \| x_k-  \mathcal{P}_{\mathcal{X}}^{B}(x_k) - \delta \ \nabla^{B} f_i(x_k) \|^2_B  \nonumber \\
    & \overset{\text{Lemma} \ \ref{1}}{=} \ \|x_k - \mathcal{P}_{\mathcal{X}}^{B}(x_k) \|_B^2 + 2\delta^2 f_i(x_k)  + 2 \delta \ \big  \langle  \mathcal{P}_{\mathcal{X}}^{B}(x_k) -x_k, \nabla^B f_i(x_k) \big \rangle_B \nonumber \\
    & \overset{\text{Lemma} \ \ref{2} }{\leq} \ \|x_k - \mathcal{P}_{\mathcal{X}}^{B}(x_k) \|_B^2  -2(2\delta-\delta^2) f_i(x_k).
\end{align}}
Now, taking expectation with respect to index $i$, we get the following:
\begin{align}
\label{eq:b2}
\E[d_B(x_{k+1},\mathcal{X})^2 \ | \ i \sim \mathcal{R}]  & \leq  \E [\| x_k-  \mathcal{P}_{\mathcal{X}}^{B}(x_k) - \delta \ \nabla^{B} f_i(x_k) \|^2_B \ | \ i \sim \mathcal{R}] \nonumber \\
& \leq \|x_k - \mathcal{P}_{\mathcal{X}}^{B}(x_k) \|_B^2  -2(2\delta-\delta^2) \E[f_i(x_k) \ | \ i \sim \mathcal{R}] \nonumber \\
& \overset{\eqref{r:mu1} }{\leq}  \ d_B(x_k,\mathcal{X})^2  -  \mu_1 (2\delta-\delta^2) \ d_B(x_k,\mathcal{X})^2  \nonumber \\
& = h_{\mathcal{R}}(\delta) \ d_B(x_k,\mathcal{X})^2.
\end{align}
Here, we used the lower bound of the function $f(x)$. Now, taking expectation again and using the tower property along with induction we get the first part of Theorem \ref{th:b1}. Similarly, considering \eqref{eq:b2} along with the bound of Theorem \ref{th:upper} we get the following:
\begin{align*}
     \E[f(x_{k+1})] \leq \frac{\mu_2}{2} \E[d_B(x_{k+1},\mathcal{X})^2] \leq    \frac{\mu_2}{2} [h_{\mathcal{R}}(\delta)]^{k+1} d_B(x_0,\mathcal{X})^2. 
\end{align*}
This proves the first part results of Theorem \ref{th:b1}. Since, $\frac{1}{k} \sum \limits_{l=0}^{k-1} \mathcal{P}_{\mathcal{X}}^{B}(x_l) \in \mathcal{X}$, using Lemma \ref{lem:distance} we have
\begin{align}
\label{eq:b3}
  \E[d_B(\Tilde{x}_k,\mathcal{X})^2]  & = \E[\|\Tilde{x}_k-  \mathcal{P}_{\mathcal{X}}^{B}(\Tilde{x}_k)\|^2_B]  \overset{\text{Lemma} \ \ref{lem:distance}}{ \leq}  \E \left[\Big \| \frac{1}{k} \sum \limits_{l=0}^{k-1} \left(x_l-\mathcal{P}_{\mathcal{X}}^{B}(x_l)\right)\Big \|^2_B\right] \nonumber \\
  & \leq  \E \left[\frac{1}{k} \sum \limits_{l=0}^{k-1} \big \| x_l-\mathcal{P}_{\mathcal{X}}^{B}(x_l)\big \|^2_B\right] = \frac{1}{k} \sum \limits_{l=0}^{k-1} \E[d_B(x_l,\mathcal{X})^2] \nonumber \\
  & \leq \frac{d_B(x_0,\mathcal{X})^2}{k} \sum \limits_{l=0}^{k-1} \left[h_{\mathcal{R}}(\delta)\right]^{l} \leq \frac{d_B(x_0,\mathcal{X})^2}{2 \delta k(2-\delta) \mu_1}.
\end{align}
Furthermore, denote $r_{k+1} = \E[d_B(x_{k+1},\mathcal{X})^2]$. Now, using \eqref{eq:b2} we have the following
\begin{align}
\label{eq:b4}
    2(2\delta-\delta^2) \sum \limits_{l=0}^{k-1} \E[f(x_l)] \ \leq \ \sum \limits_{l=0}^{k-1} (r_l-r_{l+1}) = r_0-r_k \leq r_0 = d_B(x_{0},\mathcal{X})^2.
\end{align}
Then, we get
\begin{align}
  \E[f(\Tilde{x}_k)]  & \leq  \E \left[ \frac{1}{k} \sum \limits_{l=0}^{k-1} f(x_l)\right] = \frac{1}{k} \sum \limits_{l=0}^{k-1} \E[f(x_l)] \ \leq \  \frac{d_B(x_0,\mathcal{X})^2}{2\delta k(2-\delta)}.
\end{align}
This proves the second part of Theorem \ref{th:b1}.

\subsection[Proof of Thm]{Proof of \cref{th:1}}
\label{pt1}
From the update formula of the ASPM algorithm, we get,
\begin{align}
\label{m1}
  \E[d_B(x_{k+1}, & \mathcal{X}) \ | \ i \sim \mathcal{R}]  =  \E [ \|  x_{k+1}-  \mathcal{P}_{\mathcal{X}}^{B}(x_{k+1})\|_B\ | \ i \sim \mathcal{R}]  \nonumber \\
  & \overset{\text{Lemma} \ \ref{lem:distance}}{ \leq}   \  \E [\| x_{k+1}- \mathcal{P}_{\mathcal{X}}^{B}(x_{k}) \|_B \ | \ i \sim \mathcal{R}] \nonumber \\
    & = \E [\|x_k-\mathcal{P}_{\mathcal{X}}^{B}(x_{k}) - \delta \ \nabla^{B} f_i(x_k) - \gamma (x_k-x_{k-1}) \|_B\ | \ i \sim \mathcal{R}] \nonumber \\
    & \leq  \E [\|x_k-\mathcal{P}_{\mathcal{X}}^{B}(x_{k})- \delta \ \nabla^{B} f_i(x_k)\|_B\ | \ i \sim \mathcal{R}] + \gamma \|x_k-x_{k-1}\|_B \nonumber \\
    & \leq \left\{\E [ \|x_k-\mathcal{P}_{\mathcal{X}}^{B}(x_{k})- \delta \ \nabla^{B} f_i(x_k)\|_B^2\ | \ i \sim \mathcal{R}]\right\}^{\frac{1}{2}} +  \gamma  \|x_k-x_{k-1}\|_B \nonumber \\
     & \overset{\text{Theorem} \ \ref{th:b1} }{\leq}  \sqrt{h_{\mathcal{R}}(\delta)} \  \|x_k-\mathcal{P}_{\mathcal{X}}^{B}(x_{k})\|_B + \gamma  \|x_k-x_{k-1}\|_B.
\end{align}
Now, taking expectation again in \eqref{m1} and using the tower property, we have,
\begin{align}
\label{m2}
\E[d_B(x_{k+1}, \mathcal{X})]  & \leq  \sqrt{h_{\mathcal{R}}(\delta)} \  \E[d_B(x_{k}, \mathcal{X})] + \gamma \ \E [\|x_k-x_{k-1}\|_B].
\end{align}
Similarly, using the update formula for $x_{k+1}$, we have
\begin{align}
    \label{m3}
   \E [\| x_{k+1}-x_k \|_B \ | \ i \sim \mathcal{R}] &  = \E [\|\gamma (x_k-x_{k-1}) - \delta \ \nabla^{B} f_i(x_k)\|_B \ | \ i \sim \mathcal{R}] \nonumber \\
    & \leq  \gamma \ [\|x_k-x_{k-1}\|_B + \delta  \E[\| \nabla^{B} f_i(x_k)\|_B \ | \ i \sim \mathcal{R}] \nonumber \\
    & \leq \gamma \ \|x_k-x_{k-1}\|_B + \delta  \left\{\E[\| \nabla^{B} f_i(x_k)\|^2_B \ | \ i \sim \mathcal{R}] \right\}^{\frac{1}{2}} \nonumber \\
    & = \gamma \ \|x_k-x_{k-1}\|_B + \sqrt{2} \delta  \sqrt{f(x)} \nonumber \\
    &  \overset{\text{Lemma} \ \ref{th:upper} }{\leq}  \gamma \ \|x_k-x_{k-1}\|_B + \delta \sqrt{\mu_2} \   d_B(x_{k}, \mathcal{X}).
\end{align}
Taking expectation in \eqref{m3} and using the tower property, we have,
\begin{align}
    \label{m4}
    \E [\|& x_{k+1}-  x_{k}\|_B]  \   \leq \   \gamma \ \E [\|x_k-x_{k-1}\|_B] + \delta \sqrt{\mu_2 }  \E [d_B(x_{k}, \mathcal{X})].
\end{align}
Combining both \eqref{m2} and \eqref{m4},
we can deduce the following matrix inequality:
{\allowdisplaybreaks
\begin{align}
\label{m5}
\E \begin{bmatrix}
d_B(x_{k+1}, \mathcal{X})  \\[6pt]
\|x_{k+1}-x_k\|_B 
\end{bmatrix}  & \leq  \begin{bmatrix}
\sqrt{h_{\mathcal{R}}(\delta)} &  \gamma  \\
\delta \sqrt{\mu_2 }   &  \ \gamma
\end{bmatrix} \E \begin{bmatrix}
d_B(x_{k}, \mathcal{X}) \\
\|x_k-x_{k-1}\|_B
\end{bmatrix}.
\end{align}}
Since, $(\delta, \gamma) \in Q  = \{(\delta, \gamma) \ | \ 0 < \delta < 2, \ 0 \leq \gamma <  \frac{1-\sqrt{h_{\mathcal{R}}(\delta)}}{1-\sqrt{h_{\mathcal{R}}(\delta)} + \delta \sqrt{\mu_2}}\}$, we have
\begin{align}
\label{m6}
\Pi_1  + \Pi_4-  \Pi_1 \Pi_4+ & \Pi_2 \Pi_3  = \gamma + \sqrt{h_{\mathcal{R}}(\delta)} +\gamma \delta \sqrt{\mu_2} - \gamma \sqrt{h_{\mathcal{R}}(\delta)} < 1.
\end{align}
Also, from the definition, it can be easily checked that $\Pi_1, \Pi_2, \Pi_3, \Pi_4 \geq 0$. Considering \eqref{m6}, we can check that $\Pi_1  + \Pi_4 < 1+\gamma \sqrt{h_{\mathcal{R}}(\delta)}- \gamma \delta \sqrt{\mu_2} = 1+ \min\{1, \gamma \sqrt{h_{\mathcal{R}}(\delta)}-\gamma \delta \sqrt{\mu_2}\}$. Let's define the sequences $F_k = \E [\|x_k-x_{k-1}\|_B]$ and $H_k = \E [d_B(x_{k}, \mathcal{X})]$. Now, using Theorem \ref{th:seq2}, we have
\allowdisplaybreaks{\begin{align} \label{m7}
\begin{bmatrix}
H_{k+1}  \\[6pt]
F_{k+1} 
\end{bmatrix} & \leq   \begin{bmatrix}
\Gamma_2 \Gamma_3 (\Gamma_1-1) \ \rho_1^{k}+ \Gamma_1 \Gamma_3 (\Gamma_2+1)\ \rho_2^{k} \\[6pt]
\Gamma_3 (\Gamma_1-1) \ \rho_1^{k}+ \Gamma_3 (\Gamma_2+1)\ \rho_2^{k}
\end{bmatrix} \ \begin{bmatrix}
H_{1} \\
F_1 
\end{bmatrix}.
\end{align}}
where, $\Gamma_1 , \Gamma_2, \Gamma_3, \rho_1, \rho_2$ can be derived from \eqref{t1} using the above parameter choice. Note that, from the ASPM algorithm we have, $x_1 = x_0$. Therefore we can easily check that, $F_1 = \E [\|x_1-x_{0}\|_B] = 0$ and $H_1 = \E [d_B(x_{1}, \mathcal{X})] = \E [d_B(x_{0}, \mathcal{X})] = d_B(x_{0}, \mathcal{X}) = H_0 $. Now, substituting the values of $H_1$ and $F_1$ in \eqref{m7}, we have
{\allowdisplaybreaks
\begin{align}
\label{m8}
\E \begin{bmatrix}
d_B(x_{k+1}, \mathcal{X})  \\[6pt]
\|x_{k+1}-x_k\|_B 
\end{bmatrix} & \leq \begin{bmatrix}
-\Gamma_2 \Gamma_3 \ \rho_1^{k}+ \Gamma_1 \Gamma_3 \ \rho_2^{k} \\[6pt]
- \Gamma_3 \ \rho_1^{k}+ \Gamma_3 \ \rho_2^{k} 
\end{bmatrix} \ d_B(x_{0}, \mathcal{X})  \leq \begin{bmatrix}
 \rho_2^{k} \\[6pt]
2 \Gamma_3 \ \rho_2^{k} 
\end{bmatrix} \ d_B(x_{0}, \mathcal{X}).
\end{align}}
Also from Theorem \ref{th:seq2} we have, $\Gamma_1,  \Gamma_3 \geq 0$ and $ 0 \leq |\rho_1| \leq \rho_2 < 1$. Which proves the Theorem.

\subsection[Proof of Thm]{Proof of \cref{lem:skm2}}
\label{pskm2}
\begin{proof}
Let's assume $y_{k+1} = x_k-\delta \nabla^B f_i(x_k)$. Which implies $\|y_{k+1}-\bar{x}\|_B \ \leq \  \|x_{k}-\bar{x}\|_B$ for all $k$ since $\bar{x} \in \mathcal{X} \subset \mathcal{X}_{t_k}= \{x: S_{t_k}^T(Ax-b) \leq 0\}$ and $y_{k+1}$ is the projection of $x_k$ towards or into the half-space $\mathcal{X}_{t_k}$ with respect to the $B-$norm (we assumed $x_k \notin \mathcal{X}_{t_k}$, if $x_k \in \mathcal{X}_{t_k}$ the inequality is true with equality). Furthermore, from our update formula we have $x_{k+1} = y_{k+1}+ \gamma(x_k-x_{k-1})$, that implies that the momentum term $\gamma(x_k-x_{k-1})$ forces the iterate $x_{k+1}$ to be closer to the feasible region $\mathcal{X}$ than the corresponding $B-$projection $y_{k+1}$. Therefore, we have $ \|x_{k+1}-\bar{x}\|_B \ \leq \  \|x_{k}-\bar{x}\|_B$ for all $\bar{x} \in \mathcal{X}$.
\end{proof}

\subsection[Proof of Thm]{Proof of \cref{lem:mom1}}
\label{pmom1}
From the general update formula, we have
  \begin{align}
  \label{eq:mom1}
   d_B(x_{k+1},\mathcal{X})^2 & =  \big \| x_{k+1}- \mathcal{P}_{\mathcal{X}}^{B}( x_{k+1}) \big \|^2_B   \overset{\text{Lemma} \ \ref{lem:distance}}{ \leq} \big \| x_{k+1}- \mathcal{P}_{\mathcal{X}}^{B}( x_{k}) \big \|^2_B \nonumber \\
   & = \|x_k-\mathcal{P}_{\mathcal{X}}^{B}( x_{k})-\delta \ \nabla^B f_i(x_k) + \gamma (x_k-x_{k-1}) \|^2_B \nonumber \\
   & = \|x_k-\mathcal{P}_{\mathcal{X}}^{B}( x_{k})-\delta \ \nabla^B f_i(x_k) \|^2_B + \gamma^2 \|(x_k-x_{k-1})\|^2_B  \nonumber \\
   & + 2 \gamma \delta  \langle x_{k-1}-x_{k}, \nabla^B f_i(x_k) \rangle_B  - 2 \gamma \langle x_{k-1}-x_{k}, x_k-\mathcal{P}_{\mathcal{X}}^{B}( x_{k}) \rangle_B \nonumber \\
   & = \|x_k-\mathcal{P}_{\mathcal{X}}^{B}( x_{k})\|^2_B- 2(2\delta-\delta^2)f_i(x_k)  + 2 \gamma \delta  \langle x_{k-1}-x_{k}, \nabla^B f_i(x_k) \rangle_B \nonumber \\
    &  + (\gamma^2+ \gamma) \|x_k-x_{k-1}\|^2_B + \gamma \|x_k-\mathcal{P}_{\mathcal{X}}^{B}( x_{k})\|^2_B- \gamma \|x_{k-1}-\mathcal{P}_{\mathcal{X}}^{B}( x_{k})\|^2_B.
\end{align}
Here, we used the identity $2 \langle x_{k-1}-x_{k}, x_k-\mathcal{P}_{\mathcal{X}}^{B}( x_{k})  \rangle_B = - \|x_{k-1}-\mathcal{P}_{\mathcal{X}}^{B}( x_{k}) \|^2_B + \|x_k-x_{k-1}\|^2_B + \|x_k- \mathcal{P}_{\mathcal{X}}^{B}( x_{k}) \|^2_B$. Taking expectation in \eqref{eq:mom1} with respect to index $i$ and simplifying we get 
  \begin{align}
  \label{eq:mom2}
   \E & [d_B(x_{k+1}  ,\mathcal{X})^2  \ | \ i \sim \mathcal{R}] = \|x_k-\mathcal{P}_{\mathcal{X}}^{B}( x_{k})\|^2_B- 2(2\delta-\delta^2)f(x_k) + \gamma \|x_k-\mathcal{P}_{\mathcal{X}}^{B}( x_{k})\|^2_B \nonumber \\
   & + 2 \gamma \delta  \langle x_{k-1}-x_{k}, \nabla f(x_k) \rangle  + (\gamma^2+ \gamma) \|x_k-x_{k-1}\|^2_B - \gamma \|x_{k-1}-\mathcal{P}_{\mathcal{X}}^{B}( x_{k})\|^2_B \nonumber \\
     & \overset{\text{Lemma}  \ \ref{lem:GRsketching2} }{\leq}  (1+\gamma) \| x_{k}- \mathcal{P}_{\mathcal{X}}^{B}( x_{k}) \big \|^2_B - \gamma \| x_{k-1}- \mathcal{P}_{\mathcal{X}}^{B}( x_{k-1}) \big \|^2_B - 2(2\delta-\delta^2)f(x_k) \nonumber \\
    &  \quad \quad + 2 \gamma \delta [f(x_{k-1})-f(x_k)] + (\gamma^2+ \gamma) \|x_k-x_{k-1}\|^2_B   \nonumber \\
    & = (1+\gamma) \ d_B(x_{k},\mathcal{X})^2 - \gamma \ d_B(x_{k-1},\mathcal{X})^2 + (\gamma^2+\gamma) \ \|x_{k}-x_{k-1}\|_B^2 \nonumber \\
    & \quad \quad  + 2 \gamma \delta f(x_{k-1}) - 2 \delta (2-\delta + \gamma) f(x_{k}),
\end{align}
here, we use the identity $\|x_{k-1}-\mathcal{P}_{\mathcal{X}}^{B}( x_{k-1}) \|^2_B \leq \|x_{k-1}-\mathcal{P}_{\mathcal{X}}^{B}( x_{k}) \|^2_B$. Similarly, we have,
\begin{align}
\label{eq:mom3}
   \E & [\|x_{k+1}-x_{k}\|_B^2\ | \ i \sim \mathcal{R}]   =  \E[\|\gamma (x_k-x_{k-1}) - \delta \ \nabla^B f_i(x_k)\|^2_B\ | \ i \sim \mathcal{R}] \nonumber \\
   & = \gamma^2 \ \|x_{k}-x_{k-1}\|_B^2 + \delta^2 \E[\|\nabla^B f_i(x_k)\|_B^2\ | \ i \sim \mathcal{R}] - 2 \gamma \delta \langle x_{k}-x_{k-1}, \nabla f(x_k) \rangle \nonumber \\
   & \overset{\text{Lemma}  \ \ref{lem:GRsketching2} }{\leq} \gamma^2 \ \|x_{k}-x_{k-1}\|_B^2 + 2 \delta^2 f(x_k)  + 2 \gamma \delta [f(x_{k-1})-f(x_k)] \nonumber \\
   &  = \gamma^2 \ \|x_{k}-x_{k-1}\|_B^2 + 2 \gamma \delta f(x_{k-1}) + 2 \delta (\delta - \gamma) f(x_{k}).
\end{align}
Take, $\zeta \geq 0$. Combining \eqref{eq:mom2} and \eqref{eq:mom3}, we get the required result.

\subsection[Proof of Thm]{Proof of \cref{th:mom2}}
\label{pmom2}
From Lemma \ref{lem:mom1}, we get the following
\begin{align}
\label{eq:mom4}
   \E[d_B(x_{k+1} & ,\mathcal{X})^2  \ | \ i \sim \mathcal{R}] + \zeta \E[\|x_{k+1}-x_{k}\|_B^2\ | \ i \sim \mathcal{R}]  \nonumber \\
   & \leq (1+\gamma) \ d_B(x_{k},\mathcal{X})^2 - \gamma \ d_B(x_{k-1},\mathcal{X})^2  + (\gamma^2+\gamma+ \zeta \gamma^2) \ \|x_{k}-x_{k-1}\|_B^2 \nonumber \\
    & + 2 \gamma \delta (1+\zeta) f(x_{k-1}) - 2 \delta [2-(\delta -\gamma)(1+\zeta)] f(x_{k}) \nonumber \\
    \leq &  \left\{1+\gamma + \delta \mu_1 [(1+\zeta)(\delta-\gamma)-2] \right\} \ d_B(x_{k},\mathcal{X})^2 \nonumber \\
    & + \gamma \left[\delta(1+\zeta)\mu_2-1 \right]  \ d_B(x_{k-1},\mathcal{X})^2 + (\zeta \gamma^2+ \gamma^2+ \gamma) \ \|x_{k}-x_{k-1}\|_B^2.
\end{align}
Now, let's define the sequences $H_k = \E[d_B(x_{k},\mathcal{X})^2]$ and $F_k = \E [\|x_k-x_{k-1}\|^2_B]$. Since, $x_1 = x_0$ one can easily check that, $F_1 = \E [\|x_1-x_{0}\|^2_B] = 0$ and $H_1 =  \E[d_B(x_{1},\mathcal{X})^2] = \E[d_B(x_{0},\mathcal{X})^2] = H_0 $. Taking expectation in \eqref{eq:mom4}, and using the tower property of expectation, we get
\begin{align}
    \label{mom21}
    H_{k+1} + \zeta F_{k+1}  &  \leq \left\{1+\gamma + \delta \mu_1 [(1+\zeta)(\delta-\gamma)-2] \right\} H_k \nonumber \\
    & + \gamma \left[\delta(1+\zeta)\mu_2-1 \right]  H_{k-1} + (\zeta \gamma^2+ \gamma^2+ \gamma) F_k.
\end{align}
Since, $(\delta, \gamma, \zeta) \in R \cap S$, we can easily check that the following conditions hold:
\begin{align}
\label{mom22}
   & (1+\zeta)(\delta-\gamma) \leq 2  \quad \text{and} \quad 1+\gamma + \delta \mu_1 [(1+\zeta)(\delta-\gamma)-2] \geq 0 \\
    & 0 \leq \gamma < \frac{\zeta}{1+\zeta} \quad \text{and} \quad  \gamma (1+\zeta)(\mu_2-\mu_1)+ \delta \mu_1(1+\zeta) < 2 \mu_1.
\end{align}
Next, we will analyze the recurrence relation of \eqref{mom21} with respect to the following cases: 1) $0 < \delta \mu_2(1+\zeta) \leq 1$, and 2) $ 1 < \delta \mu_2(1+\zeta) < 2 \mu_2(1+\zeta)$. In other words, we will divide the interval $(0,2]$ as  $(0,2) = (0, \frac{1}{\mu_2(1+\zeta)}] \cup (\frac{1}{\mu_2(1+\zeta)},2)$.
\paragraph{Case 1:} Assume, $0 < \delta \leq \frac{1}{\mu_2(1+\zeta)}$, then from \eqref{mom21} we have,
\begin{align}
    \label{mom23}
     H_{k+1}  + \zeta F_{k+1}    &  \leq \left\{1+\gamma \delta \mu_2(1+\zeta) + \delta \mu_1 [(1+\zeta)(\delta-\gamma)-2] \right\} H_k  + (\zeta \gamma^2+ \gamma^2+ \gamma) F_k. 
\end{align}
We used the identity $H_k \leq H_{k-1}$ from Lemma \ref{lem:skm2}. Let's take $ \alpha_1 = \zeta, \ \beta_2 = \gamma \left[\delta(1+\zeta)\mu_2-1 \right], \ \beta_3 = \zeta \gamma^2 + \gamma^2 + \gamma $ and $\beta_1 = 1+\gamma  + \delta \mu_1 [(1+\zeta)(\delta-\gamma)-2]$. One can easily check that for any $0 \leq \gamma < \frac{\zeta}{1+\zeta}$ we have $\beta_3-\alpha_1 < 0$. Moreover, from the assumed condition we have,
\begin{align*}
   0 \leq  \beta_1 + \beta_2 = 1+\gamma \delta \mu_2(1+\zeta) + \delta \mu_1 [(1+\zeta)(\delta-\gamma)-2] < 1.
\end{align*}
Now, from \eqref{mom23}, we have
\begin{align*}
    H_{k+1}+\zeta F_{k+1} \leq (\beta_1+\beta_2) H_k+ \beta_3 F_k.
\end{align*}
Which means that the sequences $H_k$ and $F_k$ satisfy the conditions of Theorem \ref{seq}. Now, using Theorem \ref{seq} we have
\begin{align}
\label{mom24}
    H_{k+1} \leq H_{k+1} + \alpha H_k +  \zeta F_{k+1} &  \leq  \rho^k \left[(1+\alpha)H_1+\alpha_1 F_{1}\right]  =  \rho^k (1+\alpha)H_0,
\end{align}
where, $ \alpha \geq 0$ and $ \rho \in [0,1)$ are given by
\begin{align}
\label{mom25}
    \alpha = \max \left\{0, \frac{\zeta \gamma^2+ \gamma^2+ \gamma}{\zeta}-\beta_1-\beta_2 \right\}, \ \rho = \max \left\{\beta_1+\beta_2, \frac{\zeta \gamma^2+ \gamma^2+ \gamma}{\zeta}\right\}.
\end{align}
Therefore, if $(\delta, \gamma, \zeta) \in R \cap S $ and $0 < \delta \leq \frac{1}{\mu_2(1+\zeta)}$, then the sequence $x_k$ generated by the ASPM algorithm converges and \eqref{mom24} holds. 

\paragraph{Case 2:} Assume, $ \frac{1}{\mu_2(1+\zeta)} < \delta < 2$, then from \eqref{mom21} we have,
\begin{align}
    \label{mom26}
    H_{k+1}  + \zeta F_{k+1}    \leq \underbrace{\left\{1+\gamma + \delta \mu_1 [(1+\zeta)(\delta-\gamma)-2] \right\}}_{\geq 0} & H_k  + \underbrace{\gamma \left[\delta(1+\zeta)\mu_2-1 \right]}_{\geq 0} H_{k-1} \nonumber \\
    & + (\zeta \gamma^2+ \gamma^2+ \gamma) F_k.
\end{align}
Now, we already show that $\beta_3 < \alpha_1$. Furthermore, we have
\begin{align*}
  0 \leq   \beta_1 + \beta_2 = 1+\gamma \delta \mu_2(1+\zeta) + \delta \mu_1 [(1+\zeta)(\delta-\gamma)-2] < 1,
\end{align*}
which are precisely the conditions of Theorem \ref{seq}. Using Theorem \ref{seq} we have
\begin{align}
\label{mom27}
   H_{k+1} \leq H_{k+1} + \alpha H_k +  \zeta F_{k+1} &  \leq  \rho^k \left[(1+\alpha)H_1+\alpha_1 F_{1}\right]  =  \rho^k (1+\alpha)H_0.
\end{align}
where, $ \alpha \geq 0$ and $ \rho \in [0,1)$ are given by
\begin{align}
    & \alpha = \max \left\{0, \frac{\zeta \gamma^2+ \gamma^2+ \gamma}{\zeta}-\beta_1, \frac{-\beta_1+ \sqrt{\beta_1^2+4 \beta_2}}{2} \right\}, \label{mom28} \\
    & \rho = \max \left\{ \frac{\zeta \gamma^2+ \gamma^2+ \gamma}{\zeta}, \frac{\beta_1+ \sqrt{\beta_1^2+4 \beta_2}}{2} \right\}. \label{mom29}
\end{align}
Therefore, if $(\delta, \gamma, \zeta) \in R \cap S $ and $ \frac{1}{\mu_2(1+\zeta)} < \delta < 2$, then the sequence $x_k$ generated by the ASPM algorithm converges and \eqref{mom27} holds. Now, we will combine the previous two cases. Since, $\beta_1+\beta_2 < 1$, we must have $\frac{\beta_1+ \sqrt{\beta_1^2+4 \beta_2}}{2} > \beta_1 + \beta_2$. Combining the above-mentioned cases, we can deduce that for any $0 < \delta < 2$, if the parameters $\gamma , \delta$ and $\zeta $ satisfies $(\delta, \gamma, \zeta) \in R \cap S$, then the sequence $x_k$ generated by the ASPM algorithm converges and the following relation holds.
\begin{align}
\label{mom30}
  \E [d_B(x_{k+1},\mathcal{X})^2]  \leq  \rho^k (1+\alpha) \  d_B(x_0,\mathcal{X})^2,
\end{align}
where, $\alpha \geq 0$ and $\rho$ are as in \eqref{mom28} and \eqref{mom29}. Furthermore, using \eqref{mom30} along with Theorem \ref{th:upper} we get the following:
\begin{align*}
     \E[f(x_{k+1})] \leq \frac{\mu_2}{2} \E[d_B(x_{k+1},\mathcal{X})^2] \leq \frac{\mu_2(1+\alpha)}{2} \rho^k  d_B(x_0,\mathcal{X})^2.
\end{align*}
This proves the first part results of Theorem \ref{th:cesaro}. Since, $\frac{1}{k} \sum \limits_{l=1}^{k} \mathcal{P}_{\mathcal{X}}^{B}(x_l) \in \mathcal{X}$, using Lemma \ref{lem:distance} we have
\begin{align}
\label{mom31}
  \E[d_B(\Tilde{x}_k,\mathcal{X})^2]  & = \E[\|\Tilde{x}_k-  \mathcal{P}_{\mathcal{X}}^{B}(\Tilde{x}_k)\|^2_B]  \overset{\text{Lemma} \ \ref{lem:distance}}{ \leq}  \E \left[\Big \| \frac{1}{k} \sum \limits_{l=1}^{k} \left(x_l-\mathcal{P}_{\mathcal{X}}^{B}(x_l)\right)\Big \|^2_B\right] \nonumber \\
  & \leq  \E \left[\frac{1}{k} \sum \limits_{l=1}^{k} \big \| x_l-\mathcal{P}_{\mathcal{X}}^{B}(x_l)\big \|^2_B\right] = \frac{1}{k} \sum \limits_{l=1}^{k} \E[d_B(x_l,\mathcal{X})^2] \nonumber \\
  & \leq \frac{d_B(x_0,\mathcal{X})^2}{k} \sum \limits_{l=1}^{k} (1+\alpha) \rho^{l-1} \leq \frac{(1+\alpha) \ d_B(x_0,\mathcal{X})^2}{ k(1-\rho)}.
\end{align}
Moreover, using \eqref{mom31} along with Theorem \ref{th:upper} we get the following
\begin{align}
  \E[f(\Tilde{x}_k)]  & \leq  \frac{\mu_2}{2} \E[d_B(\Tilde{x}_k,\mathcal{X})^2]  \ \leq \  \frac{\mu_2(1+\alpha) }{ 2k(1-\rho)} \ d_B(x_0,\mathcal{X})^2.
\end{align}
This proves the second part of Theorem \ref{th:mom2}.

\subsection[Proof of Thm]{Proof of \cref{th:4}}
\label{p4}
Since, the system $Ax \leq b$ is feasible then there exists a feasible solution $x^*$ such that the condition $|x^{*}_j| \leq \frac{2^{\sigma}}{2n}$ holds for all $j = 1, ..., n$ (see Lemma \ref{lem:skm4}). As $x_0 = 0$, we have the following bound:
\begin{align}
\label{eq:th40}
  d_B(x_0,\mathcal{X}) = \| x_0-\mathcal{P}_{\mathcal{X}}^{B}(x_0)\big \|_B  \ \leq \ \|x^*\|_B \ \leq \ \sqrt{\lambda_2} \ \frac{2^{\sigma -1}}{\sqrt{n}},
\end{align}
Now, considering Lemma \ref{lem:skm1}, we can argue that whenever the momentum algorithm runs on the system $Ax \leq b$, the system is feasible if the condition $\theta (x) < 2^{1-\sigma}$ holds. Furthermore, as $\mathcal{X} = \cap_{i \in } \{1,2,...,m\} \{x \ | \ a_i^T x \leq b_i\} $, we have the following: 
\begin{align}
\label{eq:th41}
  \theta(x) \ = \   \max_{i}(a_i^Tx-b_i)^{+} \ \leq \ \|a_i^T(x-\mathcal{P}_{\mathcal{X}}(x))\|  \leq \psi \|x-\mathcal{P}_{\mathcal{X}}(x)\| \ \leq \  \frac{\psi \  d_B(x,\mathcal{X})}{\sqrt{\lambda_1}}.
\end{align}
The, for any $(\delta, \gamma) \in Q_1$, the following bound
\begin{align}
\label{eq:th420}
  \E \left[\theta(x_{k})\right]  \overset{\eqref{eq:th41}}{\leq}   \frac{\psi \ \E[d_B(x_{k},\mathcal{X})]}{\sqrt{\lambda_1}} &  \overset{\text{Theorem} \ \ref{th:1}}{\leq} \frac{\psi \ \rho_2^{k-1}}{\sqrt{\lambda_1}}   \ d_B(x_{0},\mathcal{X}) \nonumber \\
  & \leq \sqrt{\frac{(1+\alpha)}{\lambda_1}} \psi  \rho_2^{k-1} \ d_B(x_{0},\mathcal{X}),
\end{align}
holds whenever the system $Ax \leq b$ is feasible. In the last inequality we used $\alpha \geq 0$. Similarly, if $(\delta, \gamma, t) \in R_1 \cap S_1$ for some $t \geq 0$, then the following identity holds
\begin{align}
\label{eq:th421}
  \E \left[\theta(x_{k})\right]  \overset{\eqref{eq:th41}}{\leq} \frac{\psi \ \E[d_B(x_{k},\mathcal{X})]}{\sqrt{\lambda_1}} & \leq \frac{\psi \sqrt{\E[d_B(x_{k},\mathcal{X})^2]}}{\sqrt{\lambda_1}} \nonumber \\
  &  \overset{\text{Theorem} \ \ref{th:mom2}}{\leq} \sqrt{\frac{(1+\alpha)}{\lambda_1}} \psi \rho^\frac{k-1}{2} \ d_B(x_{0},\mathcal{X}),
\end{align}
whenever the system $Ax \leq b$ is feasible. Let's denote $\bar{\rho} = \max\{\rho_2^2, \rho\}$. Combining \eqref{eq:th420} and \eqref{eq:th421}, we can deduce that the following identity
\begin{align}
\label{eq:th42}
  \E \left[\theta(x_{k})\right] \overset{\eqref{eq:th420} \ \& \ \eqref{eq:th421}}{\leq} \sqrt{\frac{(1+\alpha)}{\lambda_1}} \psi \ \bar{\rho}^{\frac{k-1}{2}} \ d_B(x_{0},\mathcal{X}) \overset{\eqref{eq:th40}}{\leq} \ \sqrt{\frac{\lambda_2(1+\alpha)}{\lambda_1}}  \psi \bar{\rho}^{\frac{k-1}{2}} \ \frac{2^{\sigma -1}}{\sqrt{n}}.
\end{align}
holds for any $(\delta, \gamma, t) \in Q_1 \cup \left(R_1 \cap S_1\right)$ whenever the system $Ax \leq b$ is feasible. However, for detecting feasibility of the system $Ax \leq b$, we must have $\E [\theta(x_k)] < 2^{1-\sigma}$. Now, from \eqref{eq:th42} we have
\begin{align*}
   \sqrt{\frac{\lambda_2(1+\alpha)}{\lambda_1}}  \psi \bar{\rho}^{\frac{k-1}{2}} \ \frac{2^{\sigma -1}}{\sqrt{n}} = \sqrt{\frac{\xi(1+\alpha)}{n}}  \psi \bar{\rho}^{\frac{k-1}{2}} 2^{\sigma -1}  < 2^{1-\sigma}.
\end{align*}
Simplifying further, we get the following bound:
\begin{align*}
   k-1 \ > \ \frac{4 \sigma - 4 -\log n + \log (1+\alpha)+ \log \xi + 2 \log \psi}{\log \left(\frac{1}{\bar{\rho}}\right)}. 
\end{align*}
Moreover, if the system $Ax \leq b$ is feasible, the probability of not having a certificate of feasibility can be calculated as follows:
\begin{align*}
  p = \mathbb{P} \left(\theta(x_k) \geq 2^{1-\sigma}\right) \ \leq \ \frac{\E \left[\theta(x_k)\right]}{2^{1-\sigma}} \ < \ \sqrt{\frac{\xi(1+\alpha)}{n}}  \psi
  \ 2^{2\sigma -2} \ \bar{\rho}^{\frac{k-1}{2}},
\end{align*}
where, we used the Markov's inequality, $\mathbb{P} (x \geq t) \leq \frac{\E[x]}{t}$. This proves the Theorem.

\subsection[Proof of Thm]{Proof of \cref{th:cesaro}}
\label{pcesaro}
First, let us  define the following sequences:
\begin{align}
\label{def:seq}
   \vartheta_l = \frac{\gamma}{1-\gamma}[x_{l}-x_{l-1}], \quad  \Delta_l = x_l + \vartheta_l, \quad \chi_l = \|x_l+\vartheta_l-\mathcal{P}_{\mathcal{X}}^{B}(\Delta_{l})\|_B^2,
\end{align}
for any natural number $l \geq 1$. Using the update formula of the momentum algorithm, we get the following identity:
\begin{align*}
    x_{l+1} + \vartheta_{l+1} \overset{\eqref{eq:momupdate}}{=} x_l + \vartheta_l - \frac{\delta}{1-\gamma} \nabla^{B} f_i(x_l),
\end{align*}
here, at iteration $l$ the index $i$ is selected based of the sampling processes described in the previous section. Now, we have
\begin{align}
\label{ces:1}
    \chi_{l+1} & = \|x_{l+1}+\vartheta_{l+1}-\mathcal{P}_{\mathcal{X}}^{B}(\Delta_{l+1})\|_B^2  \overset{\text{Lemma} \ \ref{lem:distance}}{ \leq}  \|x_{l+1}+\vartheta_{l+1}-\mathcal{P}_{\mathcal{X}}^{B}(\Delta_{l})\|_B^2 \nonumber \\
    & = \big \| x_l + \vartheta_l - \frac{\delta}{1-\gamma} \nabla^{B} f_i(x_l)-\mathcal{P}_{\mathcal{X}}^{B}(\Delta_{l}) \big \|_B^2 \nonumber \\
    & =   \underbrace{\|x_l+\vartheta_l-\mathcal{P}_{\mathcal{X}}^{B}(\Delta_{l})\|_B^2}_{= \chi_l} + \frac{\delta^2}{(1-\gamma)^2} \underbrace{\|\nabla^{B} f_i(x_l)\|_B^2}_{J_1} \nonumber \\
    & - \frac{2 \delta}{1-\gamma}   \underbrace{\big \langle x_l+\vartheta_l-\mathcal{P}_{\mathcal{X}}^{B}(\Delta_{l}) \ ,\  \nabla^{B} f_i(x_l) \big \rangle_B }_{J_2} \nonumber \\
    & = \chi_l + \frac{\delta^2}{(1-\gamma)^2} J_1- \frac{2 \delta}{1-\gamma} J_2.
\end{align}
Taking expectation with respect to index $i$, we have,
\begin{align}
  \label{ces:2}
  \frac{\delta^2}{(1-\gamma)^2} \E [J_1 \ | \ i \in \mathcal{R}] \overset{\text{Lemma} \ \ref{1} }{=}  \frac{2\delta^2}{(1-\gamma)^2} f(x_l).
\end{align}
The third term of \eqref{ces:1} can be simplified as
\begin{align}
  \label{ces:3}
   - \frac{2\delta}{1-\gamma}  & \E [J_2 \ | \ i \in \mathcal{R}]  \nonumber \\ 
   & \overset{\eqref{def:seq} }{=}  -\frac{2\delta}{1-\gamma} \big \langle x_l-\mathcal{P}_{\mathcal{X}}^{B}(\Delta_{l})  , \nabla f(x_l) \big \rangle + \frac{2\delta \gamma}{(1-\gamma)^2} \big \langle x_{l-1}-x_l ,  \nabla f(x_l) \big \rangle \nonumber \\
  & \overset{\text{Lemma} \ \ref{2} \ \& \ \ref{lem:GRsketching2}}{\leq}  - \frac{4\delta}{1-\gamma} f(x_l) +  \frac{2\delta \gamma}{(1-\gamma)^2} \left[f(x_{l-1}) -  f(x_{l}) \right].
\end{align} 
Substituting the identities of \eqref{ces:2} and \eqref{ces:3} in \eqref{ces:1} and simplifying further, we have
\begin{align}
    \label{ces:4}
    \E[\chi_{l+1} \ | \ i \in \mathcal{R}]  + \frac{2\delta \gamma (1+\delta)}{(1-\gamma)^2} f(x_l) + \varpi f(x_l) \ \leq \  \E[\chi_{l}] + \frac{2\delta \gamma (1+\delta)}{(1-\gamma)^2} f(x_{l-1}),
\end{align}
where, the term $\varpi$ is defined as
\begin{align}
 \label{ces:10}
   \varpi =   \frac{4\delta}{1-\gamma}  -  \frac{2\delta^2}{(1-\gamma)^2} =  \frac{2 \delta  (2-2 \gamma  -\delta)}{(1-\gamma)^2} \ > \ 0.
\end{align}
Taking expectation again in \eqref{ces:4} and using the tower property of expectation, we get,
\begin{align}
\label{ces:5}
    q_{l+1} + \varpi \E[f(x_l)] \leq q_l, \quad l = 1,2,3...,
\end{align}
with the sequence $q_l$ defined as $q_l =\E[\chi_{l}] + \frac{2\delta \gamma (1+\delta)}{(1-\gamma)^2} \E [f(x_{l-1})] $. Summing up identity \eqref{ces:5} for $l=1,2,...,k$, we get the following
\begin{align}
    \label{ces:6}
    \sum \limits_{l=1}^{k} \E [f(x_l)] \ \leq \ \frac{q_1-q_{k+1}}{\varpi} \ \leq \ \frac{q_1}{\varpi}.
\end{align}
Moreover, considering the Jensen inequality, we have
\begin{align*}
    \E \left[f(\bar{x_k})\right] = \E \left[f\left(\sum \limits_{l=1}^{k} \frac{x_l}{k}\right)\right] \ \leq \ \E \left[\frac{1}{k} \sum \limits_{l=1}^{k}f(x_l)\right] \ = \ \frac{1}{k}  \sum \limits_{l=1}^{k} \E [f(x_l)] \ \overset{\eqref{ces:6}}{\leq}  \frac{q_1}{\varpi k}.
\end{align*}
In the ASPM scheme, we assumed $x_0=x_1$. That implies $\vartheta_1 = \frac{ \gamma}{1-\gamma} [x_1-x_0] = 0 $. Considering these special values we have
\begin{align}
    \label{ces:7}
    \E[\chi_1] & = \E \left[\|x_1+\vartheta_1 -\mathcal{P}_{\mathcal{X}}^{B}(\Delta_{1})\|^2_B \right]  \overset{\text{Lemma} \ \ref{lem:distance}}{ \leq} \E \left[\|x_1+\vartheta_1 -\mathcal{P}_{\mathcal{X}}^{B}(x_0)\|^2_B \right] \nonumber \\
    & = \E \left[\|x_0 -\mathcal{P}_{\mathcal{X}}^{B}(x_0) \|^2_B \right] = d_B(x_{0}, \mathcal{X})^2.
\end{align}
Finally, using the definition, we get
\begin{align*}
   q_1  =\E[\chi_{1}] + \frac{2\delta \gamma }{(1-\gamma)^2} \E [f(x_{0})]  \leq \ d_B(x_{0}, \mathcal{X})^2 + \frac{2\delta \gamma }{(1-\gamma)^2} f(x_0).
\end{align*}
Now, substituting the values of $\varpi$ and $q_1$ in the expression of $\E \left[f(\bar{x_k})\right] $, we have the following
\begin{align*}
    \E \left[f(\bar{x}_k)\right] \leq \frac{ (1-\gamma)^2 \ d_B(x_{0}, \mathcal{X})^2+ 2 \gamma \delta f(x_0)}{2 \delta k \left(2-2 \gamma  -\delta\right)}.
\end{align*}
which proves the Theorem.

\section{Additional experimental results: GK \& GCD with momentum for $\tau = 1, 5, m$}
\label{appendix:exp}

\begin{figure}[htbp]
\centering
    \includegraphics[scale = 0.82]{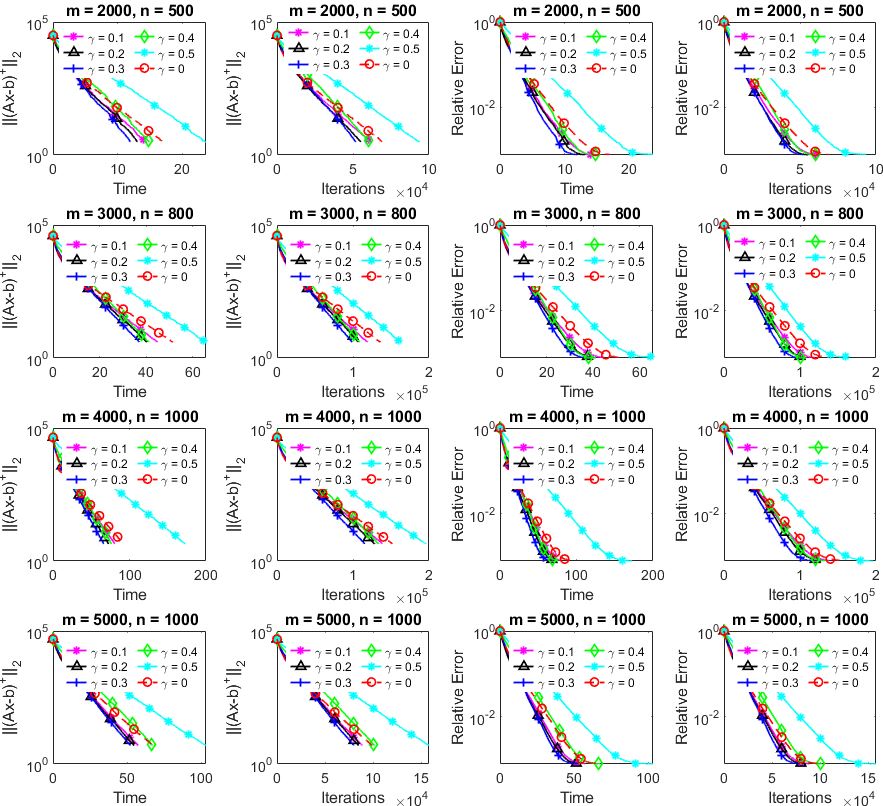}
    \caption{GK with momentum (Uniform, sketch Sample size, $\tau = 1$): comparison among momentum variants on Gaussian data, left 2 panels: Positive residual error $\|\left(Ax-b\right)^+\|_2$ vs time and No. of iterations, right 2 panels: relative error $\|x_k-x_{int}\|_B/\|x_0-x_{int}\|_B$ vs time and No. of iterations.}
    \label{fig:16}
\end{figure}
\begin{figure}[htbp]
\centering
    \includegraphics[scale = 0.79]{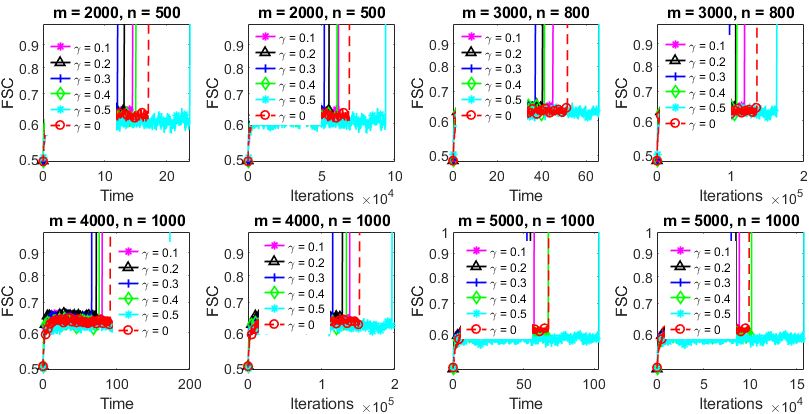}
    \caption{GK with momentum (Uniform, sketch Sample size, $\tau = 1$): comparison among momentum variants on Gaussian data, FSC vs time and No. of iterations.}
    \label{fig:17}
\end{figure}

\begin{figure}[htbp]
\centering
    \includegraphics[scale = 0.78]{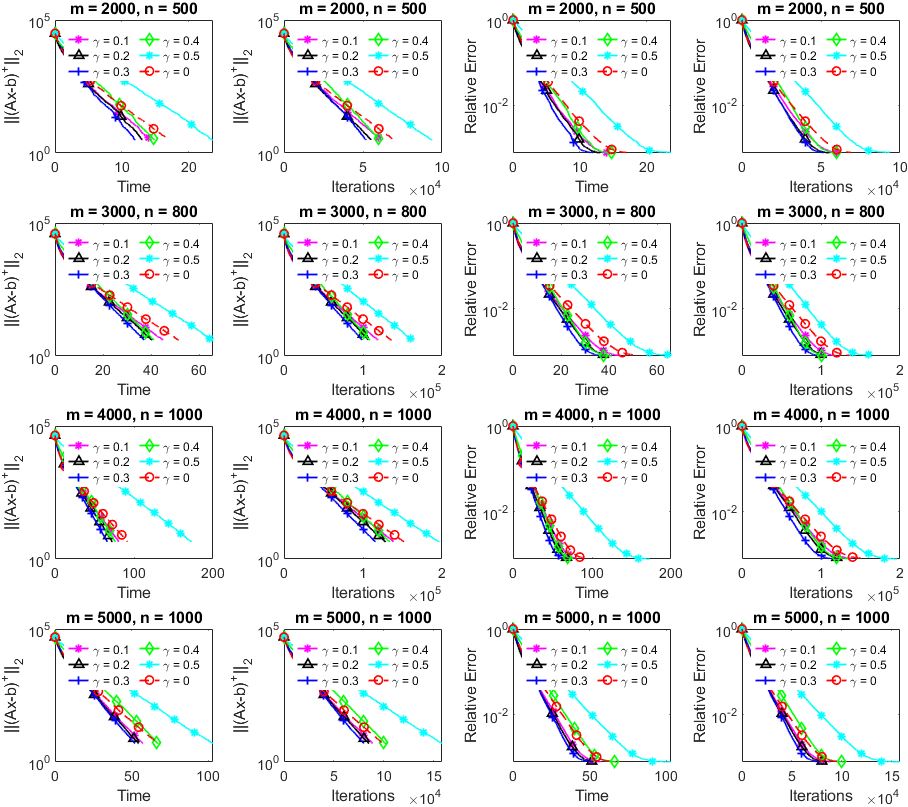}
    \caption{GK with momentum (sketch Sample size, $\tau = 5$): comparison among momentum variants on Gaussian data, left 2 panels: Positive residual error $\|\left(Ax-b\right)^+\|_2$ vs time and No. of iterations, right 2 panels: relative error $\|x_k-x_{int}\|_B/\|x_0-x_{int}\|_B$ vs time and No. of iterations.}
    \label{fig:18}
\end{figure}
\begin{figure}[htbp]
\centering
    \includegraphics[scale = 0.87]{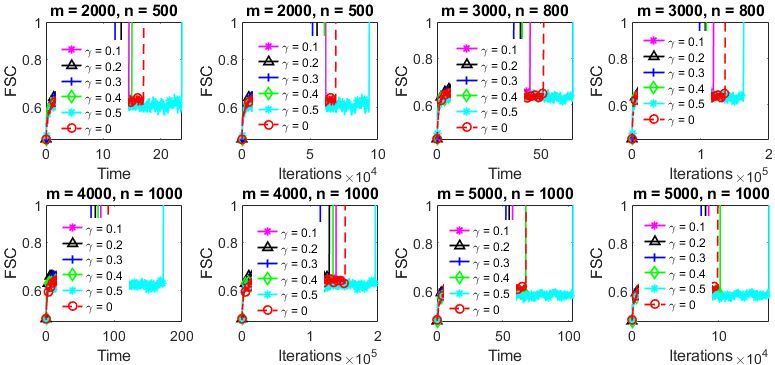}
    \caption{GK with momentum (sketch Sample size, $\tau = 5$): comparison among momentum variants on Gaussian data, FSC vs time and No. of iterations.}
    \label{fig:19}
\end{figure}

\begin{figure}[htbp]
\centering
    \includegraphics[scale = 0.76]{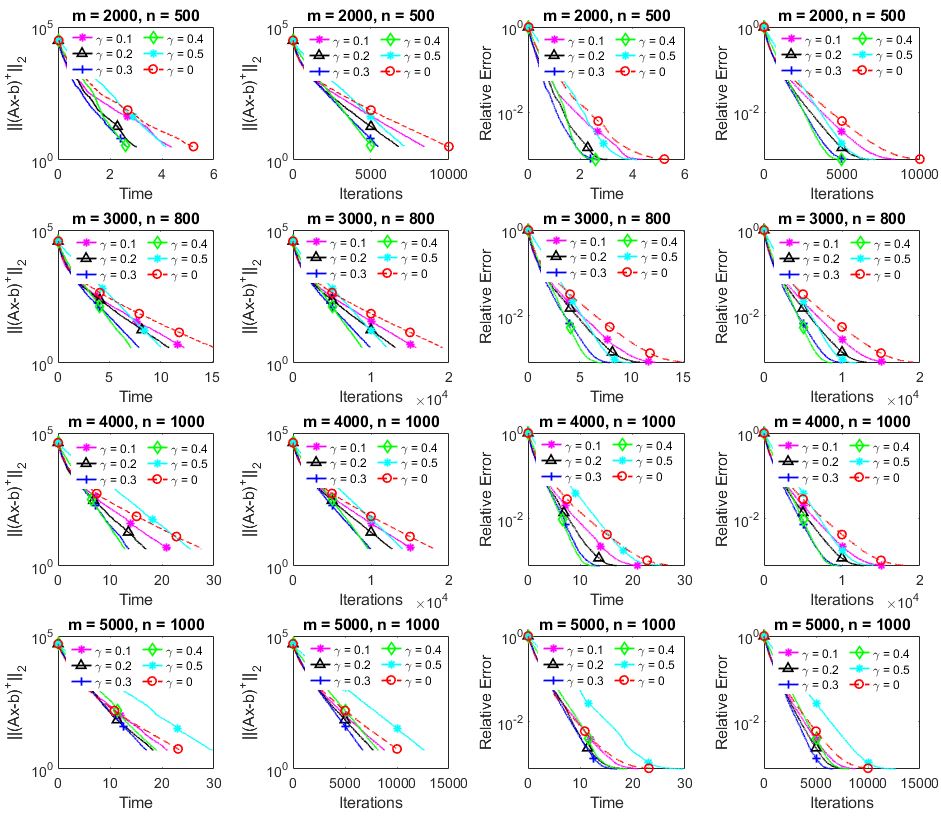}
    \caption{GK with momentum (sketch Sample size, $\tau = 50$): comparison among momentum variants on Gaussian data, left 2 panels: Positive residual error $\|\left(Ax-b\right)^+\|_2$ vs time and No. of iterations, right 2 panels: relative error $\|x_k-x_{int}\|_B/\|x_0-x_{int}\|_B$ vs time and No. of iterations.}
    \label{fig:20}
\end{figure}
\begin{figure}[htbp]
\centering
    \includegraphics[scale = 0.74]{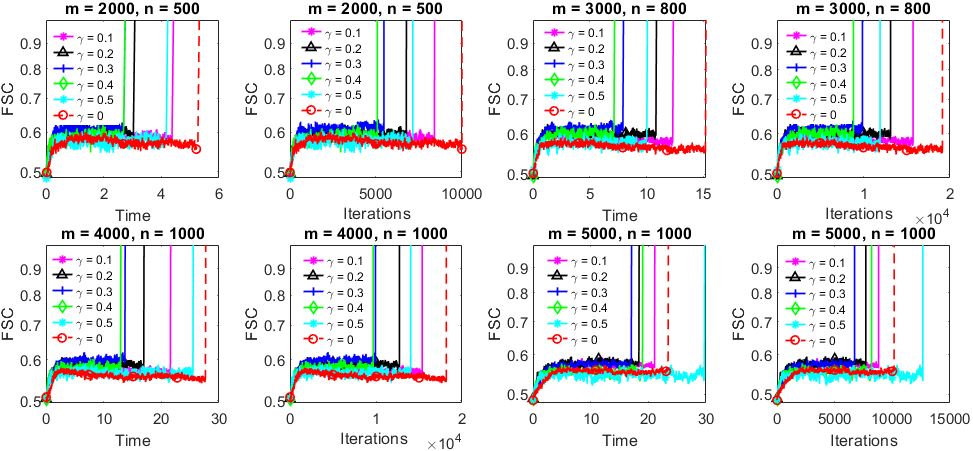}
    \caption{GK with momentum (sketch Sample size, $\tau = 50$): comparison among momentum variants on Gaussian data, FSC vs time and No. of iterations.}
    \label{fig:21}
\end{figure}

\begin{figure}[htbp]
\centering
    \includegraphics[scale = 0.76]{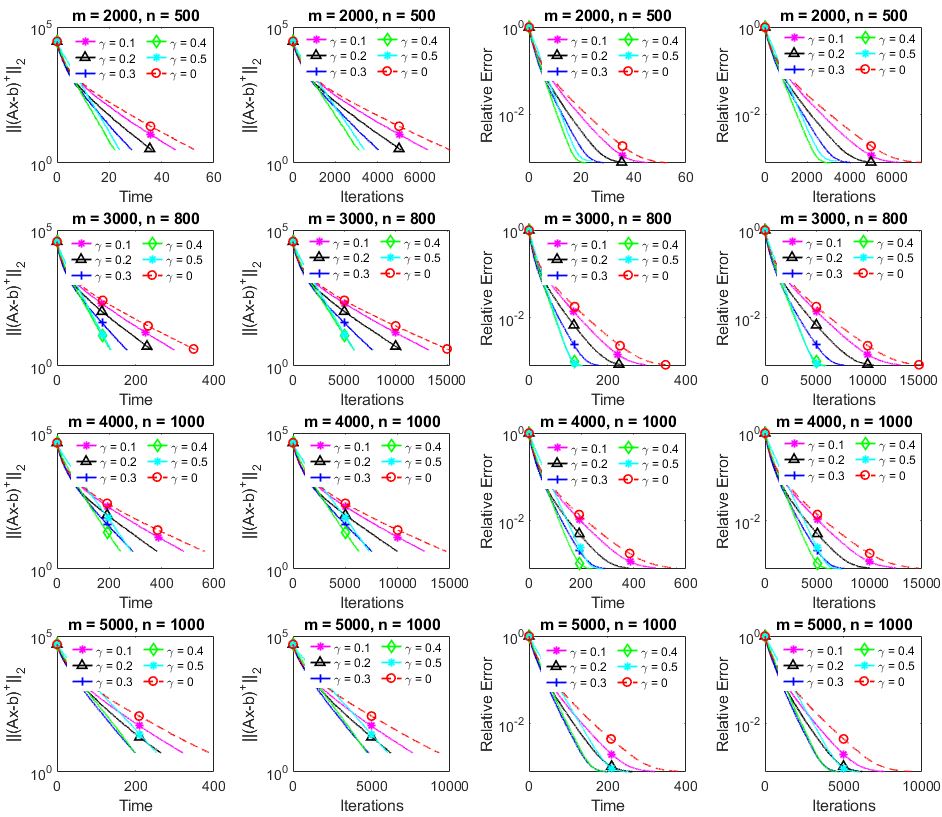}
    \caption{GK with momentum (max. distance rule, $\tau = m$): comparison among momentum variants on Gaussian data, left 2 panels: Positive residual error $\|\left(Ax-b\right)^+\|_2$ vs time and No. of iterations, right 2 panels: relative error $\|x_k-x_{int}\|_B/\|x_0-x_{int}\|_B$ vs time and No. of iterations.}
    \label{fig:22}
\end{figure}
\begin{figure}[htbp]
\centering
    \includegraphics[scale = 0.74]{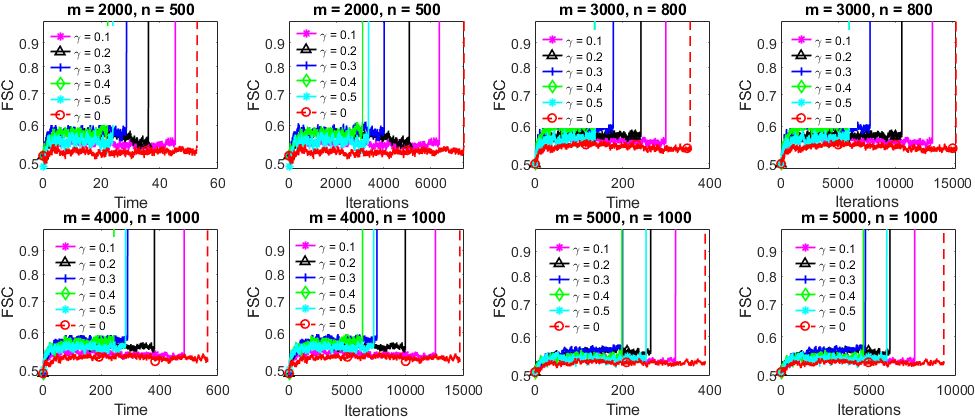}
    \caption{GK with momentum (max. distance rule, $\tau = m$): comparison among momentum variants on Gaussian data, FSC vs time and No. of iterations.}
    \label{fig:23}
\end{figure}

\begin{figure}[htbp]
\centering
    \includegraphics[scale = 0.76]{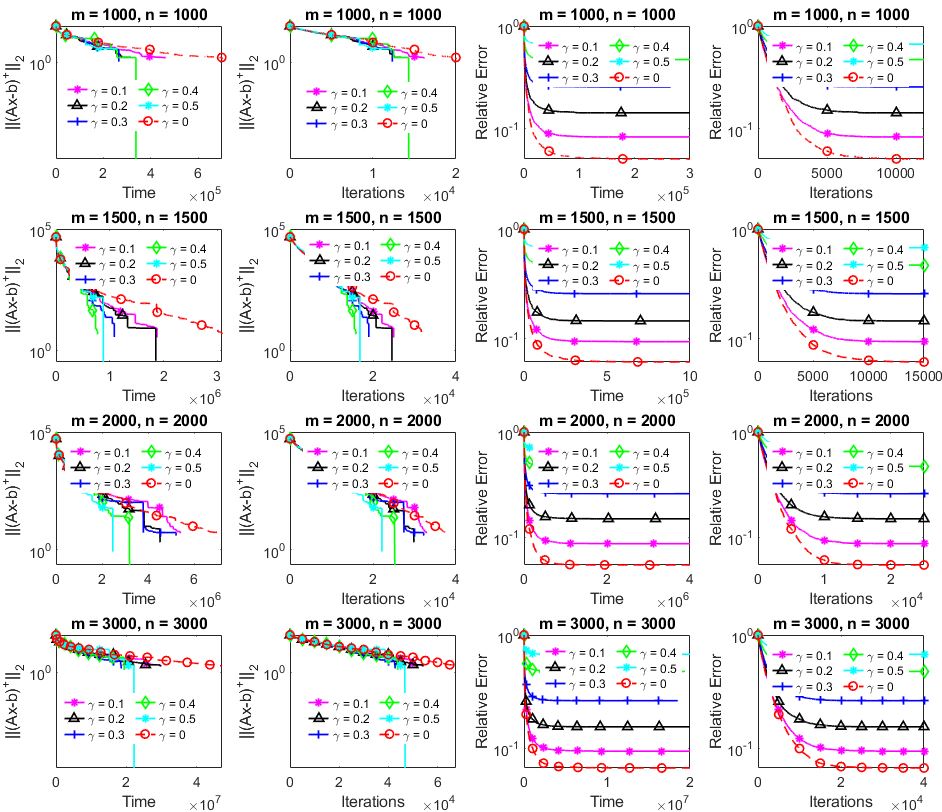}
    \caption{GCD with momentum (Uniform, sketch Sample size, $\tau = 1$): comparison among momentum variants on Gaussian data, left 2 panels: Positive residual error $\|\left(Ax-b\right)^+\|_2$ vs time and No. of iterations, right 2 panels: relative error $\|x_k-x_{int}\|_B/\|x_0-x_{int}\|_B$ vs time and No. of iterations.}
    \label{fig:24}
\end{figure}
\begin{figure}[htbp]
\centering
    \includegraphics[scale = 0.81]{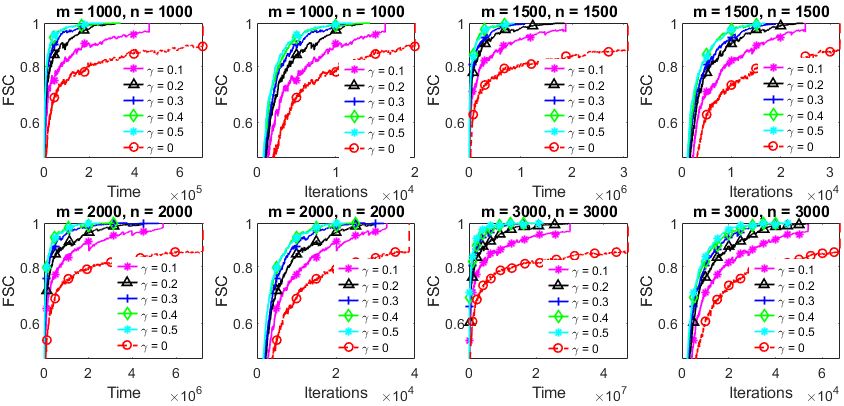}
    \caption{GCD with momentum (Uniform, sketch Sample size, $\tau = 1$): comparison among momentum variants on Gaussian data, FSC vs time and No. of iterations.}
    \label{fig:25}
\end{figure}

\begin{figure}[htbp]
\centering
    \includegraphics[scale = 0.75]{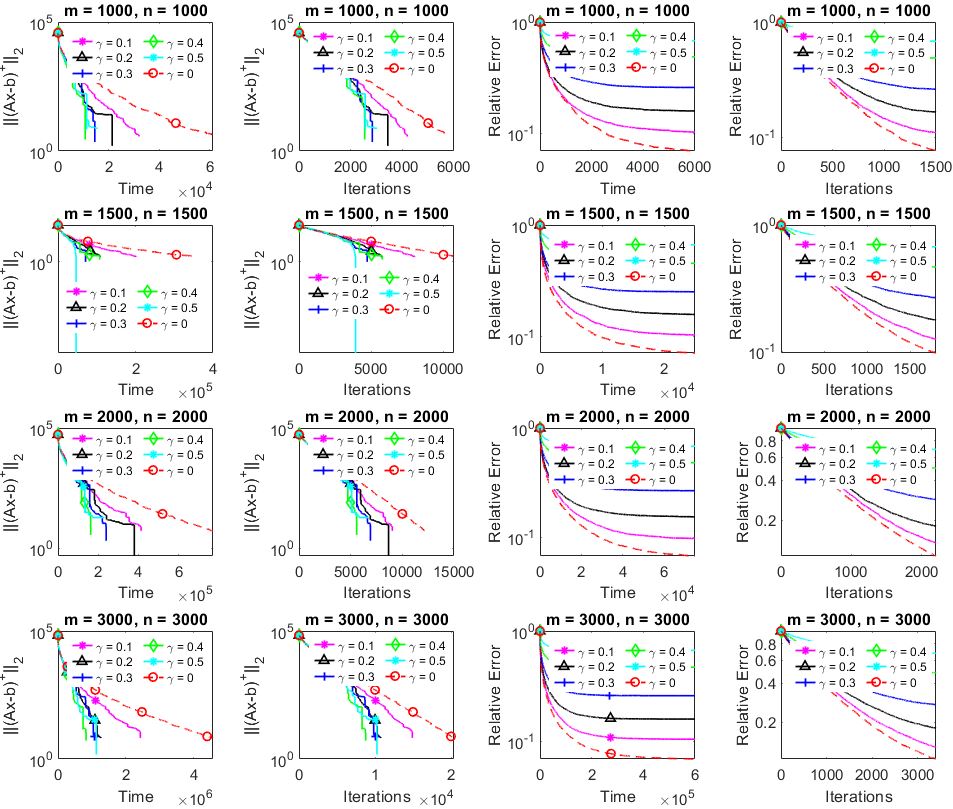}
    \caption{GCD with momentum (sketch Sample size, $\tau = 5$): comparison among momentum variants on Gaussian data, left 2 panels: Positive residual error $\|\left(Ax-b\right)^+\|_2$ vs time and No. of iterations, right 2 panels: relative error $\|x_k-x_{int}\|_B/\|x_0-x_{int}\|_B$ vs time and No. of iterations.}
    \label{fig:26}
\end{figure}
\begin{figure}[htbp]
\centering
    \includegraphics[scale = 0.77]{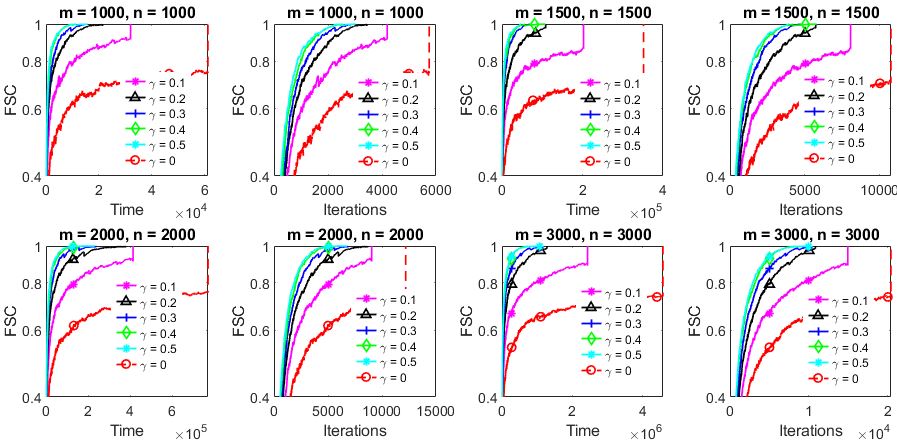}
    \caption{GCD with momentum (sketch Sample size, $\tau = 5$): comparison among momentum variants on Gaussian data, FSC vs time and No. of iterations.}
    \label{fig:27}
\end{figure}

\begin{figure}[htbp]
\centering
    \includegraphics[scale = 0.75]{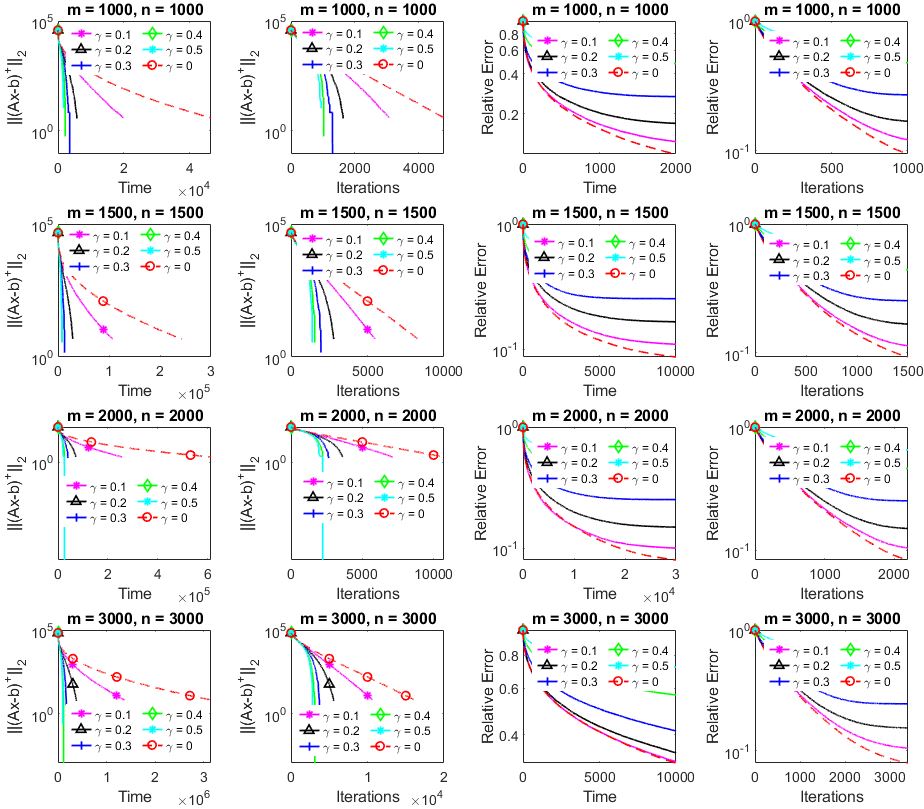}
    \caption{GCD with momentum (sketch Sample size, $\tau = 50$): comparison among momentum variants on Gaussian data, left 2 panels: Positive residual error $\|\left(Ax-b\right)^+\|_2$ vs time and No. of iterations, right 2 panels: relative error $\|x_k-x_{int}\|_B/\|x_0-x_{int}\|_B$ vs time and No. of iterations.}
    \label{fig:28}
\end{figure}
\begin{figure}[htbp]
\centering
    \includegraphics[scale = 0.76]{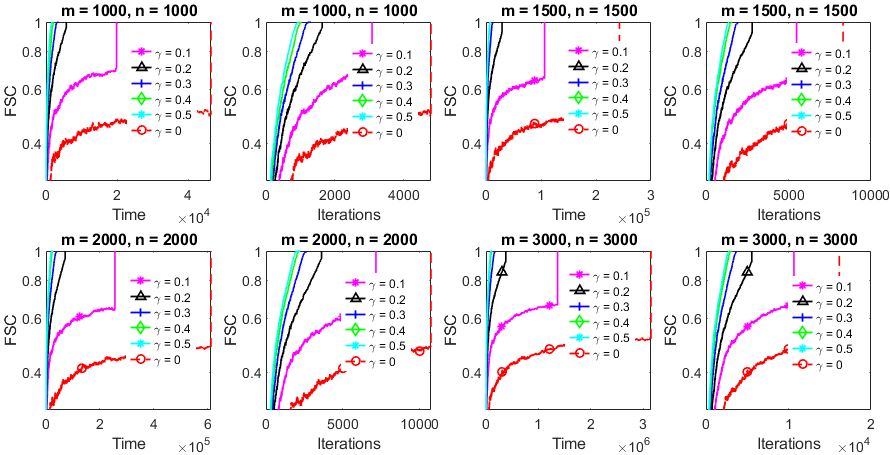}
    \caption{GCD with momentum (sketch Sample size, $\tau = 50$): comparison among momentum variants on Gaussian data, FSC vs time and No. of iterations.}
    \label{fig:29}
\end{figure}

\begin{figure}[htbp]
\centering
    \includegraphics[scale = 0.75]{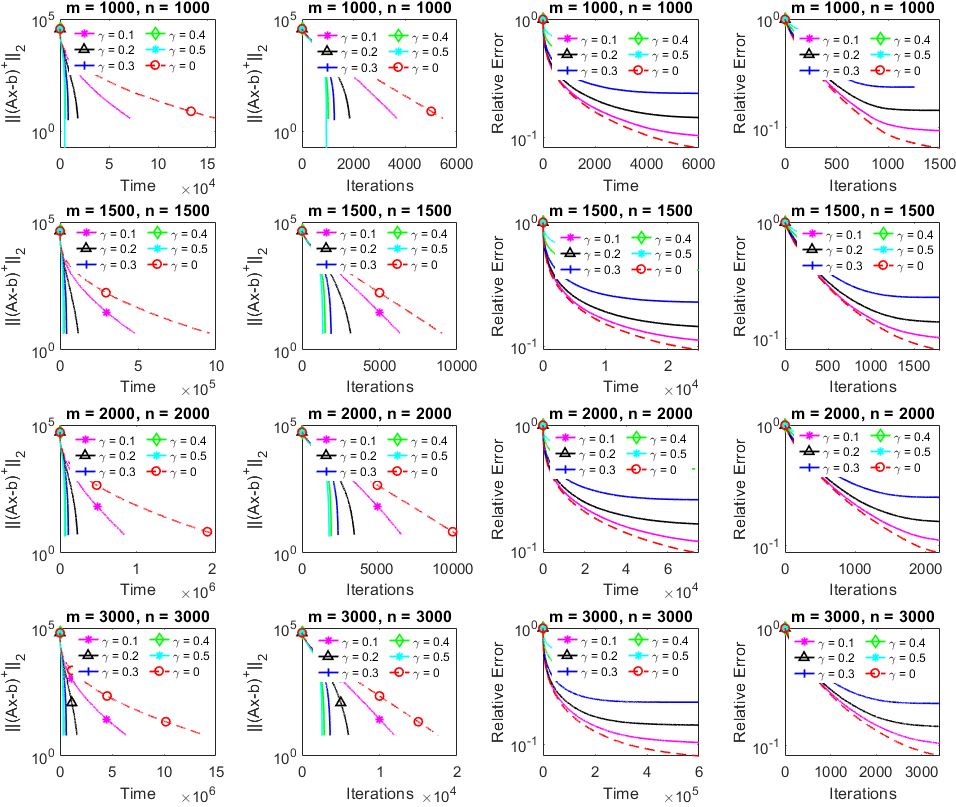}
    \caption{GCD with momentum (max. distance rule, $\tau = m$): comparison among momentum variants on Gaussian data, left 2 panels: Positive residual error $\|\left(Ax-b\right)^+\|_2$ vs time and No. of iterations, right 2 panels: relative error $\|x_k-x_{int}\|_B/\|x_0-x_{int}\|_B$ vs time and No. of iterations.}
    \label{fig:30}
\end{figure}
\begin{figure}[htbp]
\centering
    \includegraphics[scale = 0.76]{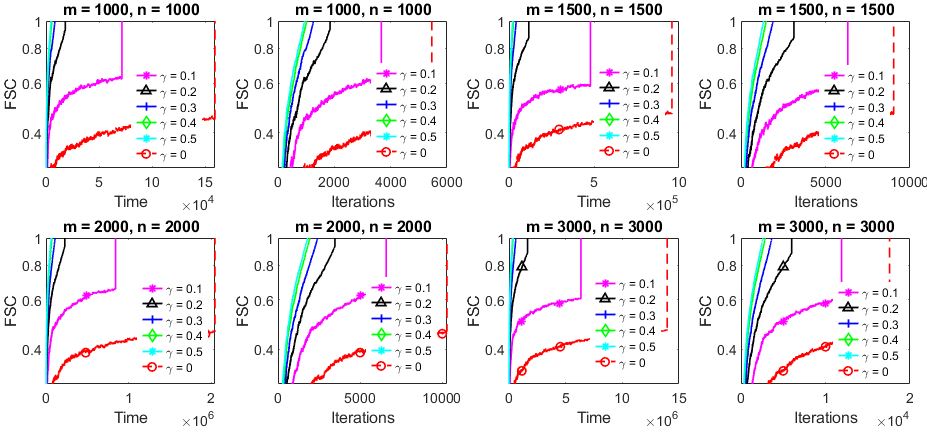}
    \caption{GCD with momentum (max. distance rule, $\tau = m$): comparison among momentum variants on Gaussian data, FSC vs time and No. of iterations.}
    \label{fig:31}
\end{figure}


\bibliographystyle{plain}
\bibliography{template}


\end{document}